%% file: main.tex
\documentclass[onefignum,onetabnum,letterpaper]{siamonline220329}

\usepackage{ifthen}
\newboolean{reportSAND}
\setboolean{reportSAND}{false}

\input{preamble/setup}
\input{preamble/setup_tikz_tensor}

\headers{Goal-Oriented Tensor Decompositions}{D. Dunlavy, E. Phipps, H. Kolla, J. Shadid, and E. Phillips}

\title{Goal-Oriented Low-Rank Tensor Decompositions for Numerical Simulation Data}

\author{
Daniel M. Dunlavy\thanks{Sandia National Laboratories (\email{\{dmdunla, etphipp, hnkolla, jnshadi\}@sandia.gov}, \email{egphill86@gmail.com})} 
\and 
Eric T. Phipps\footnotemark[1]
\and 
Hemanth Kolla\footnotemark[1]
\and 
John N. Shadid\footnotemark[1]
\and 
Edward Phillips\footnotemark[1]
}

\ifpdf
\hypersetup{ pdftitle={Goal-Oriented Tensor Decompositions}, pdfauthor={Dunlavy, Phipps, Kolla, Shadid, and Phillips} }
\fi

\begin{document}

\maketitle

\begin{abstract}
	\input{sections/abstract}
\end{abstract}



\capstartfalse
\begin{figure}[b!]
\centering
\includegraphics[width=\textwidth]{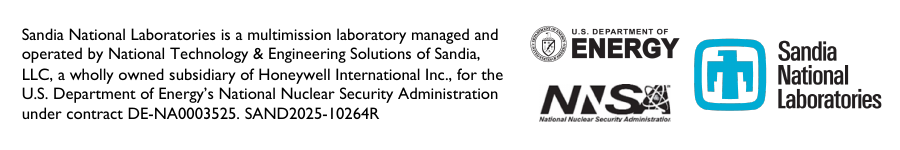}
\end{figure}
\capstarttrue

\levelstay{Introduction}
\label{sec:intro}
\input{sections/introduction}

\levelstay{Background}
\label{sec:background}
\input{sections/background}

\levelup{Goal-Oriented Tensor Decompositions}
\label{sec:goal}
\input{sections/goal_oriented}


\levelup{Application  Case Study I: Combustion}
\label{sec:exp:combustion}
\input{sections/combustion}

\levelup{Application Case Study II: Plasma Physics}
\label{sec:exp:plasma}
\input{sections/plasma}

\levelmultiup{Conclusions and Future Work}{2}
\label{sec:conc}
\input{sections/conclusion}

\section*{Acknowledgments}
\input{sections/acknowledgement}

\bibliographystyle{siamplain}
\bibliography{ref}

\appendix

\clearpage

\levelstay{Goal-oriented Gradient and Hessian}
\label{sec:go_grad_hess}
\input{sections/go_grad_hess_prec.tex}

%
%

\levelup{Plasma Physics QoI Evaluations}

\input{sections/plasma_qois}
\label{sec:app:plasma_qois}

\end{document}

%% file: preamble/setup.tex
\usepackage{amsmath,amsfonts}
\ifthenelse{\boolean{reportSAND}}{}{\usepackage{amssymb}}

\usepackage{enumitem}

\usepackage[mathscr]{eucal}

\usepackage{stmaryrd}

\usepackage{amsbsy}

\usepackage{bm}

\usepackage{braket}

\usepackage{xspace}

\usepackage{multirow}

\usepackage{xparse}

\usepackage{mathtools}

\usepackage{subcaption}

\usepackage{cleveref}
\crefname{subsection}{Section}{Sections}
\Crefname{subsection}{Section}{Sections}

\usepackage{algpseudocode}
\usepackage{algorithm}

\usepackage{lmodern}
\usepackage{listings}
\usepackage[most]{tcolorbox}

\DeclareTCBListing{macrobox}{s G{} }{IfBooleanTF={#1}{}{listing side text},title=#2,
	listing options={style=tcblatex,commentstyle=\color{red!70!black}}
}

\usepackage{booktabs}

\usepackage{coseoul}

\usepackage{fixme}
\fxsetup{
	status=draft,
	nomargin,
	inline,
	theme=color,
}
\fxsetup{marginface=\linespread{1}\footnotesize}
\FXRegisterAuthor{tk}{tke}{TK}
\FXRegisterAuthor{nj}{nje}{NJ}
\FXRegisterAuthor{ep}{epe}{EP}

\NewDocumentCommand{\Tr}{s}{\IfBooleanTF{#1}{\vphantom{\intercal}}{\intercal}}

\NewDocumentCommand{\Real} {} {\mathbb{R}}

\NewDocumentCommand{\VcRoot}{m m G{\BooleanFalse} G{\BooleanFalse}}{%
	\bm{#1\mathbf{\MakeLowercase{#2}}}%
	\IfBooleanT{#3}{^{\intercal}}%
	\IfBooleanT{#4}{^{\phantom{\intercal}}}%
}
\NewDocumentCommand{\Vc}{ O{} m t' t"} {\VcRoot{#1}{#2}{#3}{#4}}

\NewDocumentCommand{\MxRoot}{m m G{\BooleanFalse} G{\BooleanFalse}}{%
	\bm{#1\mathbf{\MakeUppercase{#2}}}%
	\IfBooleanT{#3}{^{\intercal}}%
	\IfBooleanT{#4}{^{\phantom{\intercal}}}%
}
\NewDocumentCommand{\Mx}{O{} m t' t"}{\MxRoot{#1}{#2}{#3}{#4}}

\NewDocumentCommand{\Tn}{O{} m}{\boldsymbol{#1\mathscr{\MakeUppercase{#2}}}}

\NewDocumentCommand{\Tm}{O{} m m t' t"}{\MxRoot{#1}{#2}{#4}{#5}_{(#3)}}

\NewDocumentCommand{\X} {} {\Tn{X}}

\NewDocumentCommand{\Xk} { O{k} } {\mathbf{X}_{(#1)}}

\NewDocumentCommand{\Mk} { O{k} } {\mathbf{M}_{(#1)}}

\NewDocumentCommand{\Bk} { O{k} } {\mathbf{B}_{(#1)}}

\DeclareDocumentCommand{\xi}{ s } 
{
	\IfBooleanTF{#1}
	{x_{i_1 \dots i_d}}
	{x_{i}}
}

\def\lastiter{\bar}

\NewDocumentCommand{\ttm}{O{k}}{\mathbin{\times_{\mspace{-2mu}#1}}}

\NewDocumentCommand{\Xt} {s} {\IfBooleanTF{#1}{\Mx{X}_t}{\Tn{X}_t}}
\NewDocumentCommand{\Xh} {s} {\IfBooleanTF{#1}{\Mx{X}_h}{\Tn{X}_h}}

\NewDocumentCommand{\yit}{s}{\IfBooleanT{#1}{\tilde}y_{it}}

\NewDocumentCommand{\yiht}{s}{\IfBooleanT{#1}{\tilde}y_{ih}}
\NewDocumentCommand{\cit}{s}{\IfBooleanT{#1}{\tilde}c_{it}}

\NewDocumentCommand{\Wt} {s} {\IfBooleanTF{#1}{\Mx{W}_t}{\Tn{W}_t}}
\NewDocumentCommand{\Wh} {s} {\IfBooleanTF{#1}{\Mx{W}_h}{\Tn{W}_h}}

\NewDocumentCommand{\M} {} {\Tn{M}}
\NewDocumentCommand{\Mt} {s} {\IfBooleanTF{#1}{\Mx{M}_t}{\Tn{M}_t}}
\NewDocumentCommand{\Mtold} {s} {\IfBooleanTF{#1}{\Mx[\lastiter]{M}_t}{\Tn[\lastiter]{M}_t}}
\NewDocumentCommand{\Mh} {s} {\IfBooleanTF{#1}{\Mx{M}_h}{\Tn{M}_h}}

\NewDocumentCommand{\mind}{o}{\boldsymbol{\mathbf{i}}}

\NewDocumentCommand{\A}{t'}{\mathbf{A}\IfBooleanT{#1}{^{\mkern -4mu \intercal}}}

\NewDocumentCommand{\Ak} { G{k} } {%
	\Mx{A}^{\mkern -4mu (#1)}
}
\NewDocumentCommand{\Akt} { G{k}G{t} } {%
	\Mx{A}^{\mkern -4mu (#1)}_{#2}
}
\NewDocumentCommand{\Aktt} { G{k} } {%
	\Tn{A}^{\mkern -4mu (#1)}
}

\NewDocumentCommand{\akj}{O{k} G{j}}{%
	\Vc{a}^{\mkern -3mu {(#1)}}_{#2}
}

\NewDocumentCommand{\Akold} { G{k} } {%
	\Mx[\lastiter]{A}^{\mkern -4mu (#1)}
}

\NewDocumentCommand{\lj}{G{j} t"}{w_{#1}\IfBooleanT{#2}{^{\phantom{(1)}}\!\!}}
\NewDocumentCommand{\lvec}{}{\Vc{w}}

\NewDocumentCommand{\st}{G{t} t'}{\Vc{s}_{#1} \IfBooleanT{#2}{^{\intercal}}}

\NewDocumentCommand{\HWt}{ G{t} }{\mathcal{H}_{#1}}

\NewDocumentCommand{\Y}{s}{\Tn[\IfBooleanT{#1}{\tilde}]{Y}}
\NewDocumentCommand{\Yt}{s}{\Tn[\IfBooleanT{#1}{\tilde}]{Y}_t}
\NewDocumentCommand{\Yh}{s}{\Tn[\IfBooleanT{#1}{\tilde}]{Y}_h}
\NewDocumentCommand{\Yht}{s}{\Tn[\IfBooleanT{#1}{\tilde}]{Y}_{h}}
\NewDocumentCommand{\Yk}{s}{\Mx[\IfBooleanT{#1}{\tilde}]{Y}_{\mkern -3mu (k)}}
\NewDocumentCommand{\Yn} { O{} G{n} t' t"} {\MxRoot{#1}{Y}{#3}{#4}_{#2}}
\NewDocumentCommand{\Ytk}{s}{(\Tn[\IfBooleanT{#1}{\tilde}]{Y}_t)_{\mkern -3mu (k)}}
\NewDocumentCommand{\Yhk}{s}{(\Yh)_{(k)}}
\NewDocumentCommand{\Yhtk}{s}{(\Yht)_{(k)}}
\NewDocumentCommand{\Yv}{s}{\Vc[\IfBooleanT{#1}{\tilde}]{Y}}

\NewDocumentCommand{\KT} { s } {
	\llbracket 
	\IfBooleanTF{#1}{\lvec;}{}
	\Ak{1}, \dots,  \Ak{d} \rrbracket
}

\NewDocumentCommand{\KTt} { G{\st} } {
	\llbracket 
	#1;
	\Ak{1}, \dots,  \Ak{d} \rrbracket
}
\NewDocumentCommand{\KTtt} { G{\st}G{t} } {
	\llbracket 
	#1;
	\Akt{1}{#2}, \dots,  \Akt{d}{#2} \rrbracket
}

\NewDocumentCommand{\KTdef}{O{} O{} m !g}{%
	#1\llbracket \IfValueTF{#4}{#3;#4}{#3} #2\rrbracket%
}
\NewDocumentCommand{\KTauto}{m !g}{\KT[\left][\right]{#1}{#2}}
\NewDocumentCommand{\KTbig}{m !g}{\KT[\bigl][\bigr]{#1}{#2}}
\NewDocumentCommand{\KTbigg}{m !g}{\KT[\biggl][\biggr]{#1}{#2}}
\NewDocumentCommand{\KTBig}{m !g}{\KT[\Bigl][\Bigr]{#1}{#2}}
\NewDocumentCommand{\KTBigg}{m !g}{\KT[\Biggl][\Biggr]{#1}{#2}}

\NewDocumentCommand{\Gk} { G{k} } {%
	\Mx{G}^{(#1)}
}

\NewDocumentCommand{\Zk} { G{k} t'} {\MxRoot{}{Z}{#2}_{#1}}

\NewDocumentCommand{\Zkold} { G{k} t'} {\MxRoot{}{Z}{#2}_{#1}}

\NewDocumentCommand{\Z} { t' } {\MxRoot{}{Z}{#1}}

\NewDocumentCommand{\I} {} {\mathcal{I}}

\NewDocumentCommand{\FPD} { s m m } {
	\IfBooleanTF{#1}
	{\partial #2 / \partial #3}
	{\frac{\partial #2}{\partial #3}}
}

\DeclareMathOperator{\vc}{vec}

\NewDocumentCommand{\K} {} {\Tn{K}}

\newcommand{\pluseq}{\mathrel{+}=}


%% file: preamble/setup_tikz_tensor.tex
\usepackage{pgfplots}
\usetikzlibrary{matrix}
\usetikzlibrary{positioning}
\usetikzlibrary{calc}
\usetikzlibrary{3d}
\usetikzlibrary{shadows}
\usetikzlibrary{math}
\usetikzlibrary{backgrounds}
\usetikzlibrary{external}
\usetikzlibrary{topaths}
\usetikzlibrary{pgfplots.colorbrewer}
\usetikzlibrary{decorations.shapes,decorations.text}

\usepackage{xcolor}

\definecolor{ttbackground}{gray}{0.95}
\definecolor{tthighlight}{gray}{0.7} 
\definecolor{ttannotate}{gray}{0.4}
\definecolor{ttdark}{gray}{0.3}

\definecolor{c1}{HTML}{e31a1c} 
\definecolor{c1a}{HTML}{fb9a99} 
\definecolor{c2}{HTML}{1f78b4} 
\definecolor{c2a}{HTML}{a6cee3} 
\definecolor{c3}{HTML}{33a02c} 
\definecolor{c3a}{HTML}{b2df8a} 
\definecolor{c4}{HTML}{6a3d9a} 
\definecolor{c4a}{HTML}{cab2d6} 
\definecolor{c5}{HTML}{ff7f00} 
\definecolor{c5a}{HTML}{fdbf6f} 
\definecolor{c9}{HTML}{ffff99} 
\definecolor{c7}{HTML}{b15928} 
\definecolor{c8}{HTML}{f781bf} 
\definecolor{c6}{HTML}{999999} 

\definecolor{t1}{HTML}{263B91} 
\definecolor{t2}{HTML}{006A39} 
\definecolor{t3}{HTML}{AF2930} 
\definecolor{t4}{HTML}{F16722} 
\definecolor{t5}{HTML}{00AEC8} 
\definecolor{t6}{HTML}{39B34A} 

\pgfplotsset{compat=1.18}

\tikzset{
	default colors/.style={
		tfaces=c4a,tedges=c4,
		rf1={draw=c1a,fill=c1},
		rf2={draw=c2a,fill=c2},
		rf3={draw=c3a,fill=c3},    
		ttf1={draw=c1,fill=c1a},
		ttf2={draw=c2,fill=c2a},
		ttf3={draw=c3,fill=c3a},    
	},
	label/.style={},
	f1 label/.style={label,inner xsep=0.15em,inner ysep=0em,anchor=south west},
	f2 label/.style={label,inner sep=0.1em,anchor=west},
	f3 label/.style={label,inner xsep=0.1em,inner ysep=0em,anchor=north west},
	weight label/.style={label,inner sep=0.15em,anchor=east},
}

\makeatletter

\def\tikz@lib@matrixpic@draw#1--#2\pgf@stop{%
	\begin{scope}[join=bevel,x={(\vxx,\vxy)},y={(\vyx,\vyy)}]
		\coordinate (-south west) at (0,0);
		\coordinate (-north west) at (0,#1);
		\coordinate (-north east) at (#2,#1);
		\coordinate (-south east) at (#2,0);
		\coordinate (-center) at (#2/2,#1/2);
		\coordinate (-west) at (0,#1/2);
		\coordinate (-east) at (#2,#1/2);
		\coordinate (-north) at (#2/2,#1);
		\coordinate (-south) at (#2/2,0);
		\path[matrixpic,draw=none] (-south west) rectangle (-north east);
		\begin{scope}[shift=(-south west)]
			\path[matrixpic/grid] (-south west) grid (-north east);
		\end{scope}
		\path[matrixpic,fill=none] (-south west) rectangle (-north east);
	\end{scope}  
}

\tikzset{
	pics/matrixpic/.style = {
		setup code = \tikz@lib@tensor@setup,
		background code = \tikz@lib@matrixpic@draw#1\pgf@stop
	},
	pics/matrixpic/.default={2--2},
	matrixpic/.style={draw,fill=white,fill opacity=0.75},
	matrixpic/grid/.style={step=1},
}

\makeatother

\tikzset{
	mstyle/.style={matrixpic/.append style={#1}},
	mgrids/.style={matrixpic/grid/.append style={draw,#1}},
	mfill/.style={matrixpic/.append style={fill=#1}},
	mdraw/.style={matrixpic/.append style={draw=#1}},
}

\makeatletter
\def\tikz@lib@rotensor@get#1{\pgfkeysvalueof{/tikz/rotensor/#1}}

\def\tikz@lib@rotensor@draw#1--#2--#3\pgf@stop{%
	\begin{scope}[join=bevel,x={(\vxx,\vxy)},y={(\vyx,\vyy)},z={(\vzx,\vzy)}]
		\def\offset{\tikz@lib@rotensor@get{offset}}
		\def\width{\tikz@lib@rotensor@get{width}}
		\def\xmult{\tikz@lib@rotensor@get{xmult}}
		\def\ymult{\tikz@lib@rotensor@get{ymult}}
		\def\zmult{\tikz@lib@rotensor@get{zmult}}
		\pgfmathsetmacro{\halfwidth}{\width/2}
		\coordinate (boffset) at (\xmult*\offset,0,0);
		\coordinate (aoffset) at (0,\ymult*\offset,0);
		\coordinate (coffset) at (0,0,\zmult*\offset);
		\path (0,#1,0) coordinate (-origin);
		\path (-origin) ++ (aoffset) coordinate (-mode one start)
		coordinate (-f1-north west) ++ (\width,0) coordinate (-f1-north east);
		\path (-origin) ++ (boffset) coordinate (-mode two start)
		++ (\width,0) coordinate (-f2-north west) ++ (0,-\width) coordinate (-f2-south west);
		\path (-origin) ++ (coffset) coordinate (-mode three start)
		coordinate (-f3-south west) ++ (\width,0) coordinate (-f3-south east);
		\path (0,0,0) coordinate (-mode one end) coordinate (-f1-south west) ++ (\width,0) coordinate (-f1-south east);
		\path (#2,#1,0) coordinate (-mode two end) coordinate (-f2-north east) ++ (0,-\width) coordinate (-f2-south east);
		\path (0,#1,-#3) coordinate (-mode three end) coordinate (-f3-north west) ++ (\width,0) coordinate (-f3-north east);
		\coordinate (-bbox south west) at (-mode one end);
		\coordinate (-bbox south east) at (#2,0,0);
		\path let \p1=(-mode three end),\p2=(-mode two end) in (0,\y1) coordinate (-bbox north west)
		(\x2,\y1) coordinate (-bbox north east);
		\foreach \mode/\foo in {f3-/three,f1-/one,f2-/two}
		\path[rotensor/factors,rotensor/factor \foo]
		(-\mode north west) -- (-\mode north east) -- (-\mode south east) -- (-\mode south west) -- cycle;
		\foreach \foo in {f3-,f1-,f2-,bbox }
		{
			\coordinate (-\foo south) at ($(-\foo south west)!0.5!(-\foo south east)$);
			\coordinate (-\foo north) at ($(-\foo north west)!0.5!(-\foo north east)$);
			\coordinate (-\foo east) at ($(-\foo south east)!0.5!(-\foo north east)$);
			\coordinate (-\foo west) at ($(-\foo south west)!0.5!(-\foo north west)$);
			\coordinate (-\foo center) at ($(-\foo south)!0.5!(-\foo north)$);
		}
		\coordinate (-f1 label) at (-f1-south east);
		\coordinate (-f2 label) at (-f2-east);
		\coordinate (-f3 label) at (-f3-north east);
		\coordinate (-weight label) at (-origin);
		\foreach \loc in {center,south west,south,south east,east,north east,north,north west,west}
		\coordinate (-\loc) at (-bbox \loc);
	\end{scope}
}

\tikzset{
	pics/rotensor/.style = {
		setup code = \tikz@lib@tensor@setup,
		background code = \tikz@lib@rotensor@draw#1\pgf@stop
	},
	pics/rotensor/.default={1--1--1},
	rotensor/.is family,
	rotensor,
	factors/.style={draw,fill=white},
	factor one/.style={},
	factor two/.style={},
	factor three/.style={},
	offset/.initial=0.1,
	width/.initial=0.2,
	xmult/.initial=1,
	ymult/.initial=-1,
	zmult/.initial=-2,
}

\makeatother

\tikzset{
	redges/.style={rotensor/factors/.append style={draw=#1}},
	redges/.default={black},
	rfaces/.style={rotensor/factors/.append style={fill=#1}},
	rfaces/.default={white},
	rfactors/.style={rotensor/factors/.append style={#1}},
	rfactors/.default={},
	rf1/.style={rotensor/factor one/.style={#1}},
	rf1/.default={},
	rf2/.style={rotensor/factor two/.style={#1}},
	rf2/.default={},
	rf3/.style={rotensor/factor three/.style={#1}},
	rf3/.default={},
	rwidth/.style={rotensor/width=#1},
	rwidth/.default={0.2},
	roffset/.style={rotensor/offset=#1},
	roffset/.default=0.1,
}

\makeatletter
\def\tikz@lib@ttensor@get#1{\pgfkeysvalueof{/tikz/ttensor/#1}}

\def\tikz@lib@ttensor@draw#1--#2--#3--#4--#5--#6\pgf@stop{%
	\begin{scope}[join=bevel,x={(\vxx,\vxy)},y={(\vyx,\vyy)},z={(\vzx,\vzy)}]
		\coordinate (uoffset) at (\tikz@lib@ttensor@get{xmult}*\tikz@lib@ttensor@get{uoffset},0,0);
		\coordinate (voffset) at (\tikz@lib@ttensor@get{xmult}*\tikz@lib@ttensor@get{voffset},0,0);
		\coordinate (woffset) at (0,0,\tikz@lib@ttensor@get{zmult}*\tikz@lib@ttensor@get{woffset});
		\pgfmathsetmacro{\topheight}{-#6*\tikz@lib@tensor@get{zscale}*sin(\tikz@lib@tensor@get{zangle}}
		\pgfmathsetmacro{\ystart}{#2 > 0 ? max(0,#4-#1)+max(0,#5-#4-\topheight) : max(0,#4-#1)}
		\pgfmathsetmacro{\uwidth}{#1 > 0 ? #4 : 0}
		\path (0,\ystart) coordinate (-factor one south west)
		-- ++(0,#1) coordinate (-factor one north west)
		-- ++(\uwidth,0) coordinate (-factor one north east)
		-- ++(0,-#1) coordinate (-factor one south east)
		-- cycle;
		\coordinate (-factor one center) at ($ (0,\ystart) + (#4/2,#1/2) $);
		\pic[ttensor/core] (-core) at ($ (0,\ystart) + (\uwidth,#1) + (uoffset) - (0,#4) $) {tensor=#4--#5--#6};
		\pgfmathsetmacro{\vtheight}{#2 > 0 ? #5 : 0}
		\path
		($(-core-bbox north east) + (voffset) $) coordinate (-factor two north west)
		-- ++(#2,0) coordinate (-factor two north east)
		-- ++(0,-\vtheight) coordinate (-factor two south east)
		-- ++(-#2,0) coordinate (-factor two south west)
		-- cycle;
		\coordinate (-factor two center) at ($ (-factor two south west) + (#2/2,#5/2) $);
		\pgfmathsetmacro{\wwidth}{#3 > 0 ? #6 : 0}
		\coordinate (-factor three south west) at ($ (-core-back north west) + (woffset) $);
		\path[canvas is xz plane at y=#4]
		(-factor three south west)
		-- ++(\wwidth,0) coordinate (-factor three south east)
		-- ++(0,-#3) coordinate (-factor three north east)
		-- ++(-\wwidth,0) coordinate (-factor three north west)
		-- cycle;
		\path[canvas is xz plane at y=#4] (-factor three south west) -- ++(#6/2,-#3/2) coordinate (-factor three center);
		\foreach \nbr in {one ,two ,three }
		{
			\coordinate (-factor \nbr north) at ($(-factor \nbr north west)!0.5!(-factor \nbr north east)$);
			\coordinate (-factor \nbr south) at ($(-factor \nbr south west)!0.5!(-factor \nbr south east)$);
			\coordinate (-factor \nbr east) at ($(-factor \nbr north east)!0.5!(-factor \nbr south east)$);
			\coordinate (-factor \nbr west) at ($(-factor \nbr north west)!0.5!(-factor \nbr south west)$);
		}
		\foreach \fold/\fnew in {factor one /f1-,factor three /f3-,factor two /f2-}
		\foreach \loc in {north west,north,north east,east,south east,south,south west,west,center}
		\coordinate (-\fnew \loc) at (-\fold \loc);
		\coordinate (-bbox south west) at (0,0);
		\path let \p1=(-factor three north east),\p2=(-factor two north east), \n1 = {max(\x1,\x2)} in
		(0,\y1) coordinate (-bbox north west)
		(\n1,0) coordinate (-bbox south east)
		(\n1,\y1) coordinate (-bbox north east);
		\coordinate (-bbox south) at ($(-bbox south west)!0.5!(-bbox south east)$);
		\coordinate (-bbox north) at ($(-bbox north west)!0.5!(-bbox north east)$);
		\coordinate (-bbox east) at ($(-bbox south east)!0.5!(-bbox north east)$);
		\coordinate (-bbox west) at ($(-bbox south west)!0.5!(-bbox north west)$);
		\foreach \loc in {south west,south,south east,east,north east,north,north west,west}
		\coordinate (-\loc) at (-bbox \loc);
		\path[ttensor/factors,ttensor/factor one,shift=(-f1-south west),draw=none] (-f1-south west) rectangle (-f1-north east);
		\path[ttensor/factors,ttensor/factor two,shift=(-f2-south west),draw=none] (-f2-south west) rectangle (-f2-north east);
		\path[ttensor/factors,ttensor/factor three,shift=(-f3-south west),canvas is xz plane at y=0,draw=none] (-f3-south west) rectangle (-f3-north east);    
		\path[ttensor/factor grids,ttensor/f1 grid,shift=(-f1-south west)] (-f1-south west) grid (-f1-north east);
		\path[ttensor/factor grids,ttensor/f2 grid,shift=(-f2-south west)] (-f2-south west) grid (-f2-north east);
		\path[ttensor/factor grids,ttensor/f3 grid,shift=(-f3-south west),canvas is xz plane at y=0] (-f3-south west) grid (-f3-north east);    
		\path[ttensor/factors,ttensor/factor one,shift=(-f1-south west),fill=none] (-f1-south west) rectangle (-f1-north east);
		\path[ttensor/factors,ttensor/factor two,shift=(-f2-south west),fill=none] (-f2-south west) rectangle (-f2-north east);
		\path[ttensor/factors,ttensor/factor three,shift=(-f3-south west),canvas is xz plane at y=0,fill=none] (-f3-south west) rectangle (-f3-north east);    
	\end{scope}  
}

\tikzset{
	pics/ttensor/.style = {
		setup code = \tikz@lib@tensor@setup,
		background code = \tikz@lib@ttensor@draw#1\pgf@stop
	},
	pics/ttensor/.default={2--2--2--1--1--1},
	ttensor/.is family,
	ttensor,
	factors/.style={draw,fill=white, fill opacity=0.75},
	core/.style={},
	factor one/.style={},
	factor two/.style={},
	factor three/.style={},
	uoffset/.initial=0.25,
	voffset/.initial=0.25,
	woffset/.initial=0.25,
	xmult/.initial=1,
	ymult/.initial=-1,
	zmult/.initial=-2,
	factor grids/.style={step=1},
	f1 grid/.style={},
	f2 grid/.style={},
	f3 grid/.style={},
}

\makeatother

\tikzset{
	ttcore/.style={ttensor/core/.append style={#1}},
	ttfactors/.style={ttensor/factors/.append style={#1}},
	ttf1/.style={ttensor/factor one/.append style={#1}},
	ttf2/.style={ttensor/factor two/.append style={#1}},
	ttf3/.style={ttensor/factor three/.append style={#1}},
	ttgrids/.style={ttensor/factor grids/.append style={draw,#1}},
}

\makeatletter
\def\tikz@lib@tensor@get#1{\pgfkeysvalueof{/tikz/tensor/#1}}

\def\tikz@lib@tensor@setup{%
	\pgfmathsetlengthmacro{\vxx}%
	{\tikz@lib@tensor@get{scale}*\tikz@lib@tensor@get{xscale}*cos(\tikz@lib@tensor@get{xangle})*1cm}
	\pgfmathsetlengthmacro{\vxy}%
	{\tikz@lib@tensor@get{scale}*\tikz@lib@tensor@get{xscale}*sin(\tikz@lib@tensor@get{xangle})*1cm}
	\pgfmathsetlengthmacro{\vyx}%
	{\tikz@lib@tensor@get{scale}*\tikz@lib@tensor@get{yscale}*cos(\tikz@lib@tensor@get{yangle})*1cm}
	\pgfmathsetlengthmacro{\vyy}%
	{\tikz@lib@tensor@get{scale}*\tikz@lib@tensor@get{yscale}*sin(\tikz@lib@tensor@get{yangle})*1cm}
	\pgfmathsetlengthmacro{\vzx}%
	{\tikz@lib@tensor@get{scale}*\tikz@lib@tensor@get{zscale}*cos(\tikz@lib@tensor@get{zangle})*1cm}
	\pgfmathsetlengthmacro{\vzy}%
	{\tikz@lib@tensor@get{scale}*\tikz@lib@tensor@get{zscale}*sin(\tikz@lib@tensor@get{zangle})*1cm}
}

\def\tikz@lib@tensor@draw#1--#2--#3\pgf@stop{%
	\begin{scope}[join=bevel,x={(\vxx,\vxy)},y={(\vyx,\vyy)},z={(\vzx,\vzy)}]
		\tikzset{
			canvas front face/.style={canvas is xy plane at z=0},
			canvas back face/.style={canvas is xy plane at z=-#3},
			canvas top face/.style={canvas is xz plane at y=#1},
			canvas bottom face/.style={canvas is xz plane at y=0},
			canvas right face/.style={canvas is yz plane at x=#2},
			canvas left face/.style={canvas is yz plane at x=0},
		}
		\path (0,0,0) coordinate (-front south west) coordinate (-left south west) coordinate (-bottom south west);
		\path (#2,0,0) coordinate (-front south east) coordinate (-right south west) coordinate (-bottom south east);
		\path (0,#1,0) coordinate (-front north west) coordinate (-left north west) coordinate (-top south west); 
		\path (#2,#1,0) coordinate (-front north east) coordinate (-right north west) coordinate (-top south east);
		\path (0,0,-#3) coordinate (-back south west) coordinate (-left south east) coordinate (-bottom north west);
		\path (#2,0,-#3) coordinate (-back south east) coordinate (-right south east) coordinate (-bottom north east);
		\path (0,#1,-#3) coordinate (-back north west) coordinate (-left north east) coordinate (-top north west);
		\path (#2,#1,-#3) coordinate (-back north east) coordinate (-right north east) coordinate (-top north east);
		\path ($(-front north west)!0.5!(-front north east)$) coordinate (-front north) coordinate (-top south);
		\path ($(-front south west)!0.5!(-front south east)$) coordinate (-front south) coordinate (-bottom south);
		\path ($(-front north west)!0.5!(-front south west)$) coordinate (-front west) coordinate (-left west);
		\path ($(-front north east)!0.5!(-front south east)$) coordinate (-front east) coordinate (-right west);
		\path ($(-back north west)!0.5!(-back north east)$) coordinate (-back north) coordinate (-top north);
		\path ($(-back south west)!0.5!(-back south east)$) coordinate (-back south) coordinate (-bottom north);
		\path ($(-back north west)!0.5!(-back south west)$) coordinate (-back west) coordinate (-left east);
		\path ($(-back north east)!0.5!(-back south east)$) coordinate (-back east) coordinate (-right east);
		\path ($(-front north west)!0.5!(-back north west)$) coordinate (-top west) coordinate (-left north);                            
		\path ($(-front north east)!0.5!(-back north east)$) coordinate (-top east) coordinate (-right north);
		\path ($(-front south west)!0.5!(-back south west)$) coordinate (-left south) coordinate (-bottom west); 
		\path ($(-front south east)!0.5!(-back south east)$) coordinate (-right south) coordinate (-bottom east);
		\coordinate (-left center) at ($(-left east)!0.5!(-left west)$);
		\coordinate (-right center) at ($(-right east)!0.5!(-right west)$);
		\coordinate (-top center) at ($(-top east)!0.5!(-top west)$);
		\coordinate (-bottom center) at ($(-bottom east)!0.5!(-bottom west)$);
		\coordinate (-front center) at ($(-front east)!0.5!(-front west)$);
		\coordinate (-back center) at ($(-back east)!0.5!(-back west)$);
		\coordinate (-center) at ($(-left center)!0.5!(-right center)$);
		\path (-front south west) coordinate (-bbox south west);
		\path (-back north east) coordinate (-bbox north east);
		\path let \p1=(-front south west), \p2=(-back north west) in (\x1, \y2) coordinate (-bbox north west);
		\path let \p1=(-back east), \p2=(-front south) in (\x1, \y2) coordinate (-bbox south east);
		\path ($ (-bbox north west)!0.5!(-bbox north east) $) coordinate (-bbox north);
		\path ($ (-bbox south west)!0.5!(-bbox south east) $) coordinate (-bbox south);
		\path ($ (-bbox north west)!0.5!(-bbox south west) $) coordinate (-bbox west);
		\path ($ (-bbox north east)!0.5!(-bbox south east) $) coordinate (-bbox east);
		\path ($ (-bbox north)!0.5!(-bbox south) $) coordinate (-bbox center);
		\path let \p1=(-front north),\p2=(-back east) in (\x2,\y1) coordinate (-bbox east proj);
		\path[tensor/outline] (-front south west) -- (-front north west)
		-- (-back north west) -- (-back north east)
		-- (-back south east) -- (-front south east) -- cycle;
		\foreach \ff in {back ,left ,bottom }{ 
			\path[canvas \ff face,tensor/hidden faces,tensor/\ff face] (-\ff south west) rectangle (-\ff north east);
		}
		\path[tensor/hidden edges]
		(-back north west) -- (-back south west) -- (-back south east)
		(-back south west) -- (-front south west);
		\foreach \ff in {top ,right ,front }{ 
			\path[canvas \ff face,tensor/faces,tensor/\ff face] (-\ff south west) rectangle (-\ff north east);
		}
		\foreach \ff in {top ,right ,front }{ 
			\path[canvas \ff face,tensor/grids,tensor/\ff grid] (-\ff south west) grid (-\ff north east);
		}
		\path[fill=none,tensor/edges] (-front south west) -- (-front north west)
		-- (-back north west) -- (-back north east)
		-- (-back south east) -- (-front south east) -- cycle 
		(-front north west) -- (-front north east) -- (-back north east) 
		(-front north east) -- (-front south east);
		\foreach \loc in {center,south west,south,south east,east,north east,north,north west,west,east proj}
		\coordinate (-\loc) at (-bbox \loc);
		
	\end{scope}
}

\tikzset{
	pics/tensor/.style = {
		setup code = \tikz@lib@tensor@setup,
		background code = \tikz@lib@tensor@draw#1\pgf@stop
	},
	pics/tensor/.default={1--1--1},
	tensor/.is family,
	tensor,
	faces/.style={fill=white, fill opacity=0.75},
	hidden faces/.style={},
	grids/.style={},
	outline/.style={},
	front face/.style={},
	front grid/.style={},
	right face/.style={},
	right grid/.style={},
	top face/.style={},
	top grid/.style={},
	back face/.style={},
	bottom face/.style={},
	left face/.style={},
	edges/.style={draw},
	hidden edges/.style={},
	scale/.initial=1,
	xangle/.initial=0,
	yangle/.initial=90,
	zangle/.initial=225,
	xscale/.initial=1,
	yscale/.initial=1,
	zscale/.initial=0.5445
}

\makeatother

\tikzset{
	tfaces/.style={tensor/faces/.style={fill opacity=0.75,fill=#1}},
	tfaces/.default={white},
	tedges/.style={tensor/edges/.style={draw=#1}},
	tedges/.default={},
	thidden/.style={tensor/hidden edges/.style={draw=#1}},
	thidden/.default={},
	tgrids/.style={tensor/grids/.style={draw=#1}},
	tgrids/.default={},
	tshadow/.style={tensor/outline/.style={fill=white,drop shadow}},
}

\tikzset{
	tscale/.style={tensor/scale=#1},
	tscale/.default=1,
	frontal/.style={canvas is xy plane at z=#1},
	frontal/.default=0,
	horizontal/.style={canvas is zx plane at y=#1,rotate=90,transform shape},
	horizontal/.default=0,
	lateral/.style={canvas is yz plane at x=#1,rotate=-90,transform shape},
	lateral/.default=0,
}

\tikzset{
	columnfibers/.style={tgrids=#1,
		tensor/front grid/.append style={ystep=100},
		tensor/right grid/.append style={xstep=100}},
	columnfibers/.default=black,
	rowfibers/.style={tgrids=#1,
		tensor/front grid/.append style={xstep=100},
		tensor/top grid/.append style={xstep=100}},
	rowfibers/.default=black,
	tubefibers/.style={tgrids=#1,
		tensor/right grid/.append style={ystep=100},
		tensor/top grid/.append style={ystep=100}},
	tubefibers/.default=black,
}



%% file: sections/abstract.tex
We introduce a new low-dimensional model of high-dimensional numerical simulation data based on low-rank tensor  decompositions. Our new model aims to minimize  differences between the model data and simulation data as well as \emph{functions of the model data} and \emph{functions of the simulation data}. This novel  approach to dimensionality reduction of  simulation data provides a means of directly incorporating quantities of interests and  invariants associated with conservation principles associated with the simulation data into the low-dimensional  model, thus enabling more accurate analysis of the simulation without requiring access to the full set of high-dimensional data. Computational results of applying this approach to two standard low-rank tensor decompositions of data arising from simulation of combustion and plasma physics are presented.

%% file: sections/introduction.tex
High-dimensional numerical simulation data is ubiquitous across scientific computing disciplines, including plasma physics, fluids, earth systems, and mechanics.  Often such data is generated on high-performance computing (HPC) resources and is too large to be stored, analyzed, or moved as a whole. To facilitate efficient post-simulation analysis of the data, scientists often use samples and/or low-dimensional models of the high-dimensional data to avoid these challenges. 
Here we present a new dimension reduction method for numerical simulation data that models both the high-dimensional data and functions of the high-dimensional data. 

Numerical simulation data is naturally represented as a tensor, or multi-dimensional array, consisting of the value of each simulation variable at each point in space and time.  Low-rank tensor decompositions---akin to low-rank matrix decompositions of two-dimensional data---provide low-dimensional models of tensor data and have been used effectively in a variety of data analysis tasks, including data compression, surrogate modeling, pattern identification, and anomaly detection. However, existing tensor decomposition methods target simple element-wise error metrics between the high-dimensional data and low-dimensional model of the data, often failing to faithfully represent important physics quantities of interest (QoIs) or invariants arising from conservation principles.  In this work, we introduce \emph{goal-oriented  low-rank tensor decompositions}, which incorporate QoIs directly into the optimization problems used in fitting the low-rank models to high-dimensional tensor data. Specifically, we derive two such goal-oriented decompositions---extending the common Canonical Polyadic (CP) and Tucker models---and explore how well these new low-rank models better capture several standard QoIs typically used in analyzing data produced via combustion and plasma physics simulations.

The purpose of this paper is to introduce the new goal-oriented low-rank tensor decompositions and demonstrate that they lead to improvements of low-dimensional modeling of high-dimensional numerical simulation data. As such, the target audience for this work includes both tensor decomposition experts and computational scientists.

%% file: sections/background.tex

In this section, we introduce tensors and two of the standard low-rank tensor decompositions that form the foundation of our new goal-oriented low-rank decompositions. We also discuss the connections of our approach to related work on constrained low-rank tensor decompositions.

Tensors, or multi-dimensional arrays, are a powerful means of representing relationships in multiway data~\cite{KoBa09}. For example, in the context of scientific simulations presented in this work, tensors represent the spatio-temporal solution of coupled systems of differential equations.  Analysis of high-dimensional tensor data is usually facilitated through low-rank tensor decompositions, which approximate a given tensor in terms of simpler, structured tensors defined with fewer parameters than the number of elements in the high-dimensional data.  In this work, we focus on two common tensor decompositions, the Canonical Polyadic (CP, also known as CANDECOMP/PARAFAC) and Tucker decompositions, summarized below. 

The \emph{order} of a tensor is the number of \textit{dimensions} and each tensor dimension is called a
\textit{mode}. Following standard notation~\cite{KoBa09}, lowercase letters denote scalars (tensors of order zero), e.g., $x$; bold lowercase letters denote vectors (tensors of order one), e.g., $\Vc{v}$; and bold uppercase letters denote matrices (tensors of order two), e.g., $\Mx{A}$. Tensors of order three and
higher are denoted with bold capital script letters, e.g., $\Tn{X}$.
We use multi-index notation to indicate tensor elements, i.e., $x_{\mind} \equiv \xi*$ denotes the entry $ \mind=(i_1,\dots,i_d)\in \I \equiv \set{1, \dots, I_1} \otimes \dots \otimes \set{1, \dots, I_d}$ of a $d$-way tensor $\X\in\Real^{I_1 \times \cdots \times I_d}$.   Finally, given $\X\in\Real^{I_1 \times \cdots \times I_d}$, we define $\Xk[n]$ for $n=1,\dots,d$ to be the mode-$n$ matricization of $\X$, which is a matrix of size $I_n \times (I_1\cdots I_{n-1} I_{n+1} \cdots I_d)$ whose rows are given by a linearization of the mode-$n$ slices.

\leveldown{Canonical Polyadic (CP) Decompositions}
\label{sec:cp}

\begin{figure}
	\centering
	\begin{subfigure}[c]{0.45\linewidth}
	    \centering
	    \includegraphics[width=\textwidth]{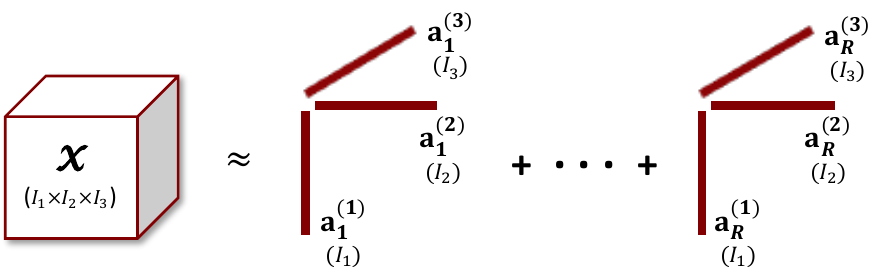}
	    \caption{}
	    \label{fig:CP-decomp}
    \end{subfigure}
	\qquad
	\begin{subfigure}[c]{0.45\linewidth}
		\centering
		\includegraphics[width=\textwidth]{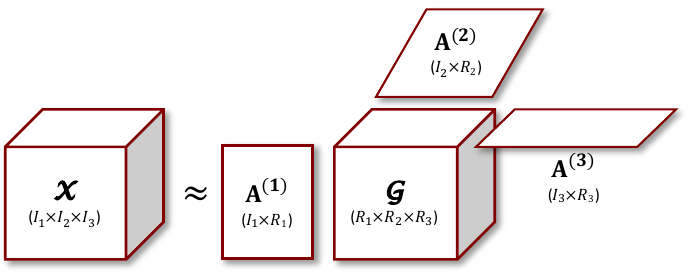}
		\caption{}
		\label{fig:HOSVD}
	\end{subfigure}
	\caption{(\subref{fig:CP-decomp}) Schematic of the canonical polyadic (CP) decomposition for a general, non-symmetric,
		third-order tensor. The input tensor is approximated as a sum of vector outer
		products.
		(\subref{fig:HOSVD}) Schematic of the Tucker decomposition for a third-order tensor.
		Mathematically, the decomposition can be written as $\Tn{X} \approx \Tn{M} = \Tn{G} 
		\times_{1} \Ak{1} \times_{2} \Ak{2} {\times_{3}} \Ak{3} $.}
\end{figure}

For a given $d$-way tensor $\X \in \Real^{I_1 \times \cdots \times I_d}$, the CP decomposition~\cite{Hi27a,Hi27,Ca44,Ca52,CaCh70,Ha70}
attempts to find a good approximating low-rank model tensor $\M$ of the form 
\begin{equation}\label{eq:CP_model}
\X \approx \M  = \sum_{r=1}^R  \akj[1] \circ \akj[2] \circ \dots \circ \akj[d],
\end{equation}
where $\akj$ is a column vector of size $I_k$,
$\circ$ represents the vector outer product, and $R$ is the approximate rank.
See~\cref{fig:CP-decomp} for a graphical depiction.
The column vectors for each mode $k$ are often collected into a matrix
$\Ak = [ \akj{1} \, \cdots  \; \akj{R} ]$
of size $I_k\times R$ called a factor matrix. 
Given factor matrices $\set{ \Ak{1}, \dots, \Ak{d} }$, 
we use the notation $\M=\llbracket \Ak{1},\dots,\Ak{d} \rrbracket$~\cite{BaKo07}. 
For standard CP decompositions, $\M$ is computed by solving the  following optimization problem,
\begin{equation}\label{eq:CP_problem}
\min_{\M} \;\; f\left(\X, \M \right) = \min_{\M} \;\;  \| \X - \M \|_F^2 \quad \mbox{s.t.} \quad \M = \llbracket \Ak{1},\dots,\Ak{d} \rrbracket,
\end{equation}
where $\| \X - \M\|_F^2 = \sum_{\mind \in \I} (x_{\mind} - m_{\mind})^2$ denotes the tensor Frobenius (sum-of-squares) norm and $\I$, defined as above, denotes the set of all multi-indices of the tensor elements.
Many approaches have been developed for efficiently solving \cref{eq:CP_problem} that are scalable to large, sparse or dense tensors including alternating least-squares (ALS)~\cite{CaCh70,Ha70},
gradient-based optimization~\cite{TTB_CPOPT}, nonlinear least-squares~\cite{Pa97}, and damped Gauss-Newton~\cite{PhTiCi13a,singh2020comparison} approaches, among others.  

\levelstay{Tucker Decompositions}
\label{sec:tucker}

The Tucker decomposition~\cite{Tu66} attempts to approximate a given tensor $\X \in \Real^{I_1 \times \cdots \times I_d}$ in the form
\begin{equation}\label{eq:tucker}
\X \approx \M = \Tn{G}\times_{1}\Ak{1}\times_{2}\dots\times_d\Ak{d},
\end{equation}
where $\Tn{G}\in\Real^{R_1\times\dots\times R_d}$ is called the core tensor and the columns of each factor matrix $\Ak{n}\in\Real^{I_n\times R_n}$, $n=1,\dots,d$ are usually orthonormal.  See~\cref{fig:HOSVD} for an example in three dimensions.  Here the notation $\Y = \Tn{X}\times_n\Mx{A}$ is called the $n$-mode product and is defined in terms of matricized tensors as $\Mx{Y}_{(n)} = \Mx{A}\Xk[n]$.  As in the CP decomposition, we use the shorthand notation $\M = \llbracket \Tn{G};\Ak{1},\dots,\Ak{d}\rrbracket$. 
 If $R_n \ll I_n$, $n=1,\dots,d$, then the memory required to store $\M$ is significantly smaller than $\X$, and hence Tucker decompositions have been used effectively for data compression~\cite{BaKlKo20,AuBaKo16}.  Each factor matrix $\Ak{n}$ is an approximate basis for the column space  corresponding to the mode-$n$ matricization of $\X$; thus, the Tucker decomposition is often considered to be a form of generalization of the matrix SVD or PCA to higher dimensions~\cite{DeDeVa00}.  
 
 Similar to the problem for CP decompositions in \cref{eq:CP_problem}, for standard Tucker decompositions, $\M$ is  computed by solving the following optimization problem
 \begin{equation}\label{eq:Tucker_problem}
 	\min_{\M} \;\; f\left(\X, \M \right) = \min_{\M} \;\;  \| \X - \M \|_F^2 \quad \mbox{s.t.} \quad \M = \llbracket \Tn{G};\Ak{1},\dots,\Ak{d}\rrbracket.
 \end{equation}
The only difference in the optimization problems between the CP and Tucker decompositions is  the specific form of the low-rank tensor $\M$.
Many approaches exist for computing Tucker decompositions, including the Higher-Order SVD (HOSVD)~\cite{DeDeVa00}, sequentially truncated HOSVD (ST-HOSVD)~\cite{VaVaMe12}, and Higher-Order Orthogonal Iteration (HOOI)~\cite{DeDeVa00a}.
For a given set of Tucker ranks, both the HOSVD and ST-HOSVD algorithms compute quasi-optimal solutions, with reconstruction error within a factor of $\sqrt d$ of the optimal approximation as measured by a sum-of-squares (Frobenius) error \cite{Hackbusch12}.
For a given approximation error tolerance, both algorithms can select Tucker ranks in order to guarantee the error tolerance is satisfied with respect to the Frobenius norm, providing a continuous trade-off between data reduction and approximation error.

\levelstay{Constrained Tensor Decompositions}

In general, when computing low-rank tensor decompositions via numerical optimization approaches, some  loss function involving the tensor data $\X$ and low-rank model $\M$ is minimized---e.g., $f$ in \cref{eq:CP_problem} or \cref{eq:Tucker_problem}---subject to the specific low-rank structure of the model $\M$---e.g., either the CP or Tucker structure described above. In addition to the constraint of the low-rank model, more general  constraints can be included in the optimization problem to  reflect additional information known about the data.

Much prior work exists  on incorporating constraints into CP and Tucker decompositions~\cite{tour_constraints}, including non-negativity (and more general bound) constraints~\cite{BHK18, cp_nn_constraints, CoVa25}, linear constraints~\cite{snect, parafac_constrained}, constraints on the rank of the decomposition~\cite{cp_rank_constraints}, and constraints presented in the form of coupled matrix-tensor factorizations~\cite{acar_cmtf, coupled}. However, such constraints amount to either bound or linear constraints, and the methods used to solve the associated constrained optimization problems are often not applicable to general nonlinear constraints. Most related to the work presented here is recent work that includes constraints on low-rank matrix and tensor decomposition factors via regularization penalties. For example, Huang, et al. introduced a general method for incorporating constraints through penalties on the factors and illustrated the use of non-negativity constraints, nuclear-norm regularization (for tensor completion), and sparsity-inducing regularization (for dictionary learning)~\cite{HuSiLi}. More recently, that work was extended to support any proximable regularizer, broadening the types of constraints that could be placed on the factors~\cite{RoScCaAd22}. As we show in Section~\ref{sec:goal}, our approach presented here also incorporates constraints as penalties to the decomposition optimization problems in \cref{eq:CP_problem} or \cref{eq:Tucker_problem}, but we illustrate the use on more general nonlinear constraints than this previous work.

\levelstay{Related Work}

As mentioned above, low-rank approximation through tensor decomposition has applications to many different computational and analysis tasks, including data compression \cite{BaKlKo20,DeCJV2024}, reduced-order modeling \cite{GhahremaniB2024}, surrogate modeling \cite{KizhakkinanDLVRR2023,Ion_Phd_thesis_2024}, and anomaly detection \cite{FanaeetG2016,LiYuan2023}, and in many of these applications, including QoIs into the low-rank approximation could potentially improve the utility of the decomposition.
Typically, only the error of primary state variables is minimized while computing the decomposition, and it has been observed that the errors in QoIs, which are derived from the primary variables, can be orders of magnitude higher \cite{Kolla2020}.
In this work, we don't focus on any particular application of tensor decomposition, but instead focus on the algorithmic, computational, and numerical consequences of computing the goal-oriented decomposition itself.  Furthermore, we recognize that tensor methods are not the only approach for solving these problems and numerous other methods for the above applications have been developed, some also incorporating goal-oriented techniques.
In particular, the literature on goal-oriented techniques for reduced order modeling of scientific simulations is vast, e.g.,
\cite{bui2007goal,fang2017efficient,borggaard2016goal,farrando2022goal,cheng2016semi}.
Furthermore, neural networks (NN) and neural operators (NO) have recently emerged as data-driven surrogate models of scientific systems, and `physics-informed' variants of these that attempt to preserve physics in addition to minimizing state prediction error have been proposed, such as physics-informed neural networks (PINNs) \cite{RaissiPK_PINN_2019} and physics-informed neural operators (PINOs) \cite{Li_etal_PINO_2024}. 
Finally, several data compression approaches that can preserve quantities of interest have recently been investigated, including MGARD~\cite{AiTuWhKl19}, which can preserve linear quantities of interest, as well as~\cite{LeGoCh22}, which proposes a hybrid approach based on tensor decompositions, autoencoders, MGARD for achieving a prescribed error guarantee, and QoI preservation as a post-processing step after reconstruction.  Our approach differs in that it produces a reduced/compressed model that attempts to directly preserve the relevant QoIs (including nonlinear QoIs) without post-processing.

%% file: sections/goal_oriented.tex

\leveldown{Problem Formulation}

As described above, traditional tensor decompositions can often be formulated as the solution to optimization problem such as 
\begin{equation}\label{eq:traditional}
\min_{\M} f(\X, \M) \quad \mbox{s.t. $\M$ is of CP or Tucker form,}
\end{equation}
where $f$ is chosen based on the statistical model of the data tensor $\X$.  In this work, we use $f(\X, \M) = \|\X-\M\|_F^2$, i.e., sum-of-squares loss, which is appropriate when the error tensor $\X-\M$ is assumed to be normally distributed (however, many other choices are possible).  We further assume the data tensor $\X$ represents the spatio-temporal evolution of some multi-physics system on a fixed Cartesian grid, consisting of one or more spatial modes (our examples below contain either two or three spatial dimensions)\footnote{This includes the modeling of systems involving more complex geometries and unstructured meshes by collapsing all spatial dimensions into a single mode given by all of the vertices in the mesh.}, a variable mode representing the different physics variables present in the system (e.g., density, temperature, momentum, etc.), and a temporal mode.  We model quantities-of-interest (QoIs), whether they be invariants that are satisfied by the data or merely quantities that are derived from the data, as scalar valued functions $g_{q}(\X_t)$  for a given set of time points $t\in\mathbb{T}_q$, where $\X_t$ represents a temporal slice of $\X$ at time $t$ (i.e., $\X_t = \X(:,:,:,:,t)$ for a tensor with three spatial modes) and $\mathbb{T}_q \subseteq \{1,\dots,\tau\}$ is the set of indices in the temporal mode of the tensor used when computing the $q$th QoI.  To preserve these QoIs in the tensor decomposition, one would like to constrain the tensor decomposition to enforce $g_{q}(\M_t) = g_{q}(\X_t)$.  However, unless the data $\X$ is exactly low rank, it is unlikely these constraints can be satisfied exactly and still obtain reasonable levels of reduction and some accuracy in the data reconstruction so that $\M\approx\X$.  We therefore attempt to preserve them through the following penalty formulation,
\begin{equation}\label{eq:goal_oriented}
\min_{\M} \alpha_0 f(\X, \M) + \sum_{q=1}^{Q} \alpha_{q} \sum_{t\in\mathbb{T}_q} \left(g_{q}(\X_t)-g_{q}(\M_t)\right)^2 \; ,
\end{equation}
where $Q$ is the number of QoIs and $\alpha_0$,\dots,$\alpha_{Q}$, are user-chosen weighting coefficients.  Note that the quadratic penalty in~\cref{eq:goal_oriented} could be replaced with other penalty functions; e.g., an augmented Lagrangian penalty~\cite{NocedalWright}. While not explored here, the use of alternate penalty functions to incorporate the QoIs into the low-rank decompositions could be explored through future work.
We next discuss multiple derivative-based optimization approaches for solving \cref{eq:goal_oriented} followed by several practical considerations that must be addressed for any optimization-based approach.

\levelstay{Optimization Approaches}
\label{sec:optimization}
\input{sections/optimization}

\levelup{Solution Procedure}
\label{sec:sol_procedure}

Many multiphysics problems involve solution variables that exhibit a wide range of scales that can cause numerical difficulties when solving~\cref{eq:goal_oriented} if not handled appropriately.  To obtain more accurate low-rank models of this data, we first perform centering and scaling of each variable slice in $\X$ to produce a scaled tensor $\Tn{\tilde{X}}$ where, e.g., for a 5-way data tensor $\X\in\Real^{I_1\times I_2 \times I_3 \times I_4 \times \tau}$,
\begin{equation}
  \Tn{\tilde{X}}(i_1,i_2,i_3,v,t) = \frac{\X(i_1,i_2,i_3,v,t) - \Vc{\mu}(v)}{\Vc{\sigma}(v)}.
\end{equation}
The computation of $\Vc{\mu}$ and $\Vc{\sigma}$ are problem dependent, with common choices of the variable mean for $\Vc{\mu}$ and standard deviation for $\Vc{\sigma}$.  For traditional tensor decompositions, one can compute a corresponding low-rank model $\Tn{\tilde{M}}$ from $\Tn{\tilde{X}}$ and then transform back to unscaled variables by unscaling the corresponding variable factor matrix, e.g., for CP:
\begin{equation}\label{eq:scaling_transformation}
  \M = S(\Tn{\tilde{M}}) = S(\llbracket \Mx{A}, \Mx{B}, \Mx{C}, \Mx{D}, \Mx{E} \rrbracket) \equiv \llbracket \Mx{A}, \Mx{B}, \Mx{C}, \mbox{diag}(\Vc{\sigma})\Mx{D}, \Mx{E} \rrbracket + \llbracket \Vc{1}_{I_1}, \Vc{1}_{I_2}, \Vc{1}_{I_3}, \Vc{\mu}, \Vc{1}_{\tau} \rrbracket
\end{equation}
where $\Vc{1}_n$ is a length-$n$ vector of 1's and $S$ is the function that maps the scaled data and model to the unscaled data and model. However, QoI evaluations often only make sense when evaluated on unscaled tensors, and thus we modify \cref{eq:goal_oriented} to explicitly unscale data before QoI evaluation:
\begin{equation}\label{eq:scaled_goal_oriented}
\min_{\Tn{\tilde{M}}} f_{go}(\Tn{\tilde{X}},\Tn{\tilde{M}}) \;\; \mbox{where} \;\; f_{go}(\Tn{\tilde{X}},\Tn{\tilde{M}}) \equiv \alpha_0 f(\Tn{\tilde{X}}, \Tn{\tilde{M}}) + \sum_{q=1}^{Q} \alpha_{q} \sum_{t\in\mathbb{T}_q} \left(g_{q}(S(\Tn{\tilde{X}}_t))-g_{q}(S(\Tn{\tilde{M}}_t))\right)^2,
\end{equation}
which requires modifying the derivative formulas presented in \cref{sec:go_grad_hess} to incorporate \cref{eq:scaling_transformation} when computing the goal derivatives. 

Once the scaling vectors have been determined, the next steps in finding an approximate solution to \cref{eq:scaled_goal_oriented} are obtaining a good initial guess and choosing the weighting coefficients.  For the former, we first obtain an approximate solution to 
\begin{equation}
  \min_{\Tn{\tilde{M}}} f(\Tn{\tilde{X}},\Tn{\tilde{M}})
\end{equation}
using a traditional approach such as CP-ALS when $\Tn{\tilde{M}}$ is of CP form and ST-HOSVD when $\Tn{\tilde{M}}$ is of Tucker form, to provide the initial guess $\Tn{\tilde{M}}^0$.  We then fix the CP/Tucker rank of $\Tn{\tilde{M}}$ based on the corresponding rank of $\Tn{\tilde{M}}^0$ and choose the weighting coefficients as
\begin{equation}
\begin{aligned}
  \alpha_0 &= \frac{1}{Q+1} \frac{1}{f(\Tn{\tilde{X}},\Tn{\tilde{M}}^0)}, \\
  \alpha_{q} &= \frac{1}{Q+1} \frac{1}{\sum_{t\in\mathbb{T}_q}\left(g_{q}(S(\Tn{\tilde{X}}_t))-g_{q}(S(\Tn{\tilde{M}}^0_t))\right)^2}, \;\; q=1,\dots,Q,
\end{aligned}
\end{equation}
which ensures the tensor term and each goal-term contributes equally to the objective function. 

When presenting results of the various low-rank decompositions in~\Cref{sec:exp:combustion} and~\Cref{sec:exp:plasma}, we report relative tensor reconstruction error using the unscaled data and model $\|\Tn{X}-\Tn{M}\|_F/\|\Tn{X}\|_F$, scaled data and model $\|\Tn{\tilde{X}}-\Tn{\tilde{M}}\|_F/\|\Tn{\tilde{X}}\|_F$,  or both.

\levelstay{Non-uniqueness and Stopping Criteria}
\label{sec:stop_crit}
It is well-known that CP and Tucker models obtained from the solution of optimization problems such as \cref{eq:traditional} are not locally unique and instead lie along a manifold of equivalent solutions~\cite{TensorBook}.  For CP, this non-uniqueness arises from the CP factor scaling ambiguity whereby any factor matrix column may be scaled by a nonzero value and not change the resulting tensor reconstruction as long as the same column of the factor matrices for the other dimensions are scaled appropriately so that the product of scales is unity.  For Tucker, the situation is even worse since any factor matrix may be multiplied by a nonsingular matrix and not change the tensor reconstruction as long as the core is multiplied by the inverse of that matrix in the corresponding dimension.  Since the QoI functions in \cref{eq:goal_oriented} or \cref{eq:scaled_goal_oriented} rely on tensor reconstructions, these same non-uniqueness properties are inherited within the goal-oriented formulation.  This lack of uniqueness can cause problems when monitoring the norm of the gradient of the objective function to determine stopping criteria, because it often does not converge monotonically to zero and instead exhibits a saw-tooth pattern (while the objective function is constant along the manifold of equivalent solutions at a point, and the gradient is therefore orthogonal to the tangent space of the manifold at that point, steps with a nonzero component along the tangent space can cause the norm of the gradient to increase since the manifold is curved).  Thus determining appropriate stopping criteria is somewhat challenging.
However, in this work we are not aiming to find the CP/Tucker model that best minimizes the QoI error, but rather to reduce this error over what is obtained in a traditional, non-goal-oriented low-rank decomposition formulation.  Thus our approach is to run the optimization method for a fixed number of iterations.  In particular, because of the quick convergence of the trust-region Newton optimization method, we found 5 optimization iterations to be sufficient for achieving significant goal reduction in most of our experiments (however most of the experiments below used 20 iterations to ensure the results were near a local minimum).  Furthermore, with the choice of weighting coefficients described above, the objective function $f_{go}$ evaluated at the initial guess satisfies $f_{go}(\Tn{\tilde{X}},\Tn{\tilde{M}}^0)=1$, and if a solution $\Tn{\tilde{M}}^\ast$ is found that exactly preserves the QoIs such that $g_{q}(S(\Tn{\tilde{X}}_t)) = g_{q}(S(\Tn{\tilde{M}}_t^\ast))$ 
while not changing the tensor loss (i.e., $f(\Tn{\tilde{X}},\Tn{\tilde{M}}^\ast) = f(\Tn{\tilde{X}},\Tn{\tilde{M}}^0)$), then $f_{go}(\Tn{\tilde{X}},\Tn{\tilde{M}}^\ast)=1/(Q+1)$.  These provide approximate (though not strict) upper and lower bounds on the objective function $f_{go}$, and by comparing $f_{go}$ to $1/(Q+1)$ throughout the optimization process, one can gauge how well the solution to \cref{eq:goal_oriented} has reduced the QoI error.

We note that the non-uniqueness challenges of these optimization formulations can be rectified by incorporating additional constraints into the problem formulation.  For CP this is relatively straightforward by constraining each factor matrix column to unit-norm and introducing an additional $R$ weight variables (where $R$ is the CP rank) into the optimization problem.  For Tucker, however, formulating explicit algebraic constraints that eliminate the above ambiguity is more challenging.  Instead, one can formulate both CP and Tucker on Riemannian manifolds and leverage existing work on optimization methods over Riemannian manifolds (e.g.,~\cite{genrtr}).  In the Tucker case, this is done by constraining the columns of the factor matrices to be orthonormal and recognizing the resulting core is not a free variable, but rather the projection of the data onto spaces spanned by the factor matrix columns in each dimension~\cite{ElSa09}.  By eliminating the core, one can then pose and solve the Tucker minimization over a product of Grassmann manifolds~\cite{ElSa09,DeHo16}.  However, extending such a method to the goal-oriented Tucker formulation is non-trivial since the constraint between the core and factor matrices is no longer satisfied and therefore can't be formulated as a product Riemannian manifold.  While not explored here, this manifold optimization approach could be explored through future work.

\levelstay{Software Implementation}
\label{sec:software}
For the numerical experiments below, we implemented the goal-oriented formulation in \cref{eq:scaled_goal_oriented} in MATLAB, leveraging the Tensor Toolbox for MATLAB~\cite{matlab_ttb_dense,BaKo07} dense tensor, CP, and Tucker data structures as well as the Matricized Tensor Times Khatri-Rao Product (MTTKRP) and Tensor Times Matrix (TTM) operations needed for the gradient and Hessian calculations described in \cref{sec:go_grad_hess}.  As described in \cref{sec:tr_newton} we used the Manopt trust region-Newton as the optimization solver.  While effective for investigating the feasibility of this approach, this serial, Matlab-based implementation was found to be quite limiting in the size of simulation data sets that could be studied.  In particular, the approach for MTTKRP in the Tensor Toolbox for MATLAB for dense tensors involves explicitly forming the Khatri-Rao product, which for large data sets and large CP ranks is memory prohibitive.  

\levelstay{Computational Cost}
\label{sec:cost}
The goal-oriented formulation in \cref{eq:scaled_goal_oriented} adds an additional $Q$ terms to the objective function defining the CP/Tucker model, so the computational cost of the goal-oriented approach is necessarily greater than a traditional approach, even one based on similar optimization methods such as trust-region Newton/Gauss-Newton.  Comparing \cref{eq:gocp_jac_trans_vec} with \cref{eq:cp_grad} and \cref{eq:gotucker_jac_trans_vec} with \cref{eq:tucker_grad}, we see the cost of computing the gradient of the additional QoI terms is similar to the gradient cost for the Frobenius loss term, once each QoI derivative tensor $\Tn{Z}$ has been computed.  Note that these tensors can be combined across multiple QoIs to compute the goal-oriented gradient by applying \cref{eq:gocp_jac_trans_vec}/\cref{eq:gotucker_jac_trans_vec} once instead of $Q$ times, and so the gradient cost of the goal-oriented approach is about twice that of a traditional, non-goal-oriented approach (not including the cost to calculate the QoI derivative tensors).  Furthermore, comparing \cref{eq:gotucker_jac_trans_vec} and \cref{eq:gotucker_jac_vec} with \cref{eq:tucker_grad} and \cref{eq:tucker_jac_vec}, we see the cost of Hessian-vector products is again similar between the QoI terms and the Frobenius loss term for the Tucker method.  However, by comparing \cref{eq:gocp_jac_trans_vec} and \cref{eq:gocp_jac_vec} with \cref{eq:cp_hess_vec}, we see the special structure of Frobenius loss does allow for more efficient Hessian-vector products for the traditional CP model than can be computed for the goal-oriented approach, since it requires $d$ CP model reconstructions.  Moreover, the cost of computing the derivative tensor $\Tn{Z}$ for each QoI can be large.  Each derivative tensor is the same size and shape as the original data $\Tn{X}$, and in most applications, the cost of computing it should be proportional to the total number of entries in $\Tn{X}$.  However, the finite element integrations defining the QoIs derivatives for the plasma physics examples shown below are challenging to vectorize within our Matlab-based implementation, making the goal-oriented approach substantially more expensive than a traditional CP or Tucker decomposition method.  We would expect the cost to be more competitive with a (Newton-based) traditional optimization method in an HPC environment where such integrals can be implemented more efficiently.


%% file: sections/optimization.tex

In this section, we describe two derivative-based optimization methods proposed to solve the optimization problems defining the goal-oriented CP and Tucker decompositions above.  The first is the Limited-Memory Broyden-Fletcher-Goldfarb-Shanno (L-BFGS) quasi-Newton algorithm commonly used for tensor decomposition methods, followed by trust-region Newton methods employing truncated Conjugate-Gradient inner solvers and Gauss-Newton Hessian approximations.

\leveldown{Quasi-Newton}

The L-BFGS algorithm has been widely used for tensor decomposition, including CP~\cite{TTB_CPOPT}, generalized CP~\cite{HoKoDu2020}, and Tucker~\cite{HeChLi23} decompositions, as it provides superlinear convergence rates but only requires specification of the gradient of the objective function.  The algorithm works by constructing secant approximations of the Hessian matrix over multiple steps of the optimization algorithm, which can be implemented by low-rank updates to the Hessian and its inverse (see~\cite{NocedalWright} for more details).  Our experiments specifically employ the implementation provided by the L-BFGS-B code~\cite{ZhByLuNo1997} available in the Tensor Toolbox for MATLAB.  Formulas for computing the needed gradients for the goal-oriented Tucker and CP decompositions are described in \cref{sec:go_grad_hess}.  

\levelstay{Trust-Region Newton}
\label{sec:tr_newton}

As will be demonstrated in \cref{sec:results}, our numerical experiments indicated the above L-BFGS method had difficulty in achieving significant reduction in the goal terms compared to a traditional tensor decomposition, particularly for the Tucker method.  However, Newton-type methods such as trust region-Newton, Gauss-Newton and Levenberg-Marquardt have proven to be effective for tensor decomposition in other contexts, e.g.,~\cite{PhTiCi2011,PhTiCi13a,RaBa20,VaVeLa21}.  So we also explored applying such methods to the goal-oriented formulation, in particular focusing on trust-region Newton methods.  The trust-region Newton method approximates the objective function at each step by a local quadratic model constructed from the objective function gradient and Hessian.  The quadratic model is minimized in the trust-region via Newton's method.  The implementation used here was provided by the Manopt~\cite{manopt} MATLAB package for optimization on Riemannian manifolds (which would allow inclusion of constraints such as optimization over Grassmann manifolds for Tucker decompositions as described later in \cref{sec:stop_crit}, but those capabilities are not explored here).  In Euclidean space, the Manopt trust-region method~\cite{genrtr} is equivalent to the method described in~\cite{NocedalWright}, using a truncated conjugate gradient (tCG) linear solve for the inexact Newton step.  This provides much faster quadratic convergence over the L-BFGS method described above, but requires implementing (approximate) Hessian-vector products for use in CG.  
Since the optimization problems defining the goal-oriented CP and Tucker decompositions are a form of nonlinear least-squares, we found the Gauss-Newton Hessian approximation to be effective, which only requires implementation of the first derivative of the goal functions.  Formulas for computing Gauss-Newton Hessian-vector products are derived in \cref{sec:go_grad_hess} for both CP and Tucker models.
The truncated CG solves are also improved through the use of a preconditioner.  In this work, we employed simple block-diagonal preconditioning in the CP case, only considering the Frobenius norm term in the objective function, where the diagonal blocks are explicitly constructed and inverted through Cholesky decompositions.  In the Tucker case, even just forming these blocks was too memory intensive, so we resorted to diagonal preconditioning.   Note that while it can be shown that each goal-oriented term introduces a rank-$|\mathbb{T}_q|$ correction to the diagonal blocks arising from the Frobenius norm term, effective incorporation of these terms in the preconditioner will be studied in future work.  Construction of the diagonal blocks and their inverses is standard in the CP and Tucker literature (e.g.,~\cite{PhTiCi2011,TensorBook}), and they are summarized in \cref{sec:fro_grad_hess_prec} for reference.

%% file: sections/combustion.tex

\leveldown{Background} 
\label{sec:hcci_background}

We use a dataset of direct numerical simulation (DNS) of turbulent combustion in conditions representing
a homogeneous charge compression ignition (HCCI) engine. The simulation \cite{BhagatwalaCL2014} considers
autoignition of a turbulent ethanol-air mixture in a 2D spatial domain, performed with the DNS code S3D
\cite{Chen_S3D} which solves a system of PDEs governing conservation of mass, momentum, energy, and chemical
species mass fractions described in \cref{tab:Combustion_System}. S3D uses an explicit eighth-order 
finite difference scheme for spatial derivatives and explicit fourth-order Runge-Kutta scheme for temporal 
derivatives. A structured grid rectangular Cartesian domain with uniform grid spacing and uniform time step 
are also typical choices (although not strictly necessary). 

%
\begin{table}[tp]
\begin{center}
\begin{tabular}{|c|c|}
\hline
\tiny \parbox{1.5in}{\scriptsize Mass} &  \scriptsize \parbox{4.0in} {\[ \frac{\partial \rho }{\partial t} + \nabla \cdot (\rho {\bf u}) = 0 \] } \\ \hline
\tiny \parbox{1.5in}{\scriptsize Momentum} & \scriptsize\parbox{4.0in} {\[ \frac{\partial {{\bf m}}} {\partial t} + \nabla \cdot \left[ {{\bf m} \otimes {\bf u}} + ( P{\bf I} + {{\bf \Pi}})\right]  = {\bf 0};~{\bf m} \equiv
 \rho {\bf u} \]} \\ \hline
\tiny \parbox{1.5in}{\scriptsize Total Energy } &  \scriptsize\parbox{4.0in} {\[ \frac{\partial (\rho e_t)}{\partial t} + \nabla \cdot \left[ \rho {\bf u} e_t + {\bf q} \right] -  {{\bf \Pi} : {\nabla {\bf u}}}  = {\bf 0};~ e_t \equiv 
e + \frac{ {\bf u} \cdot {\bf u}}{2} \]} \\  \hline
\tiny \parbox{1.5in}{\scriptsize Species Mass Fractions } &  \scriptsize\parbox{4.0in} {\[ \frac{\partial \rho {\bf y}}{\partial t} + \nabla \cdot \left[ {{\bf m} \otimes {\bf s}} + {\bf d} \right] - {\boldsymbol \dot{\omega}} = {\bf 0} \]} \\  \hline 
\end{tabular}
\end{center}
\caption{System of PDEs for low Mach number turbulent combustion solved by S3D. $\rho$ is density, $\bf u$ is 
the velocity vector, $\bf m$ is momentum vector, $P$ is pressure, $\bf \Pi$ is the viscous shear stress tensor, 
$e_t$ and $e$ are total and internal energies respectively, $\bf q$ is the heat flux vector, $\bf s$ is the 
vector of chemical species mass fractions, $\bf d$ is the species diffusion flux, and $\boldsymbol \dot{\omega}$ is 
the vector of species chemical production rates. In addition to the PDEs, a set of constitutive laws complete
the description of the computational model: thermodynamic relationships compute temperature, $T$, as a function
of internal energy, $e$; pressure is prescribed using an equation of state, $P = P(\rho {\bf s}, T)$; molecular
transport laws relate $\bf \Pi$ to $\nabla {\bf u}$, $\bf q$ to $\nabla T$, and $\bf d$ to $\nabla {\bf s}$,
respectively; finite-rate chemical kinetics and Arrhenius equation provide species production rates ${\boldsymbol 
\dot{\omega}}  = {\boldsymbol \dot{\omega}} (\rho {\bf s}, P, T)$. }
\label{tab:Combustion_System}
\end{table}
The HCCI simulation is performed in a doubly-periodic $x-y$ spatial domain, discretized with an equi-spaced
$672\times672$ grid. The data contains, at each grid point, the solution primary variables that include velocities
($u_x, u_y$), temperature, pressure, and 28 chemical species ($\rho {\bf s}$) that capture the finite-rate chemical 
kinetics of ethanol combustion. A total of 626 snapshots in time are considered. 
Accordingly, the HCCI data tensor is denoted as $\X \in \Real^{I_1 \times I_2 \times I_3 \times \tau}$, where

\begin{itemize}
        \item $I_1$ represents the spatial discretization in $x$ ($I_1 = 672$);
        \item $I_2$ represents the spatial discretization in $y$ ($I_2 = 672$);
        \item $I_3$ represents the variables, including chemical species ($I_3 = 32$); and
        \item $\tau$ represents the temporal snapshots ($\tau = 626$).
\end{itemize}

Following~\cite{Kolla2020}, the HCCI data are provided in dimensionless units by dividing each variable slice by the maximum of magnitude of that variable across all spatial grid and time points, and therefore explicit scaling/unscaling operations as described in \cref{sec:sol_procedure} are not required.

\levelstay{Quantities of Interest}


The tensor entries for the combustion data represent solution values at the discrete mesh points which are
uniformly spaced, which simplifies the expressions for the QoIs and the associated operations on the tensors.
Any integral over the spatial domain $\Omega$ is equivalent to a summation over the first two (or three for
3D cases) spatial modes times a constant factor that denotes an elemental volume: 
\[
\int_{\Omega} () ~d\Omega = \sum_{i_1=1}^{|I_1|}\sum_{i_2=1}^{|I_2|} () ~\Delta x \Delta y. 
\]
Likewise, any integral over time is equivalent to a summation over the last tensor mode. QoIs that are
functions of the solution variables imply an operation over the variables mode at each point in space/time.

\leveldown{Conservation of Mass} 
Since the HCCI simulation is performed over a doubly periodic spatial domain, it is a closed system that does not exchange mass
with the surroundings. Hence, the mass, computed as volume integral of density, is conserved in time\footnote{
Strictly speaking, a finite difference discretization is not conservative, but we expect the discretization
errors to be much smaller than errors due to the low-rank tensor representation.}. Note that density itself 
is not stored in the data files, rather it is implied by the species mass fractions vector, via the relationship
\[
\sum_{i_3=1}^{28} s_{i_3} = 1 \implies \rho = \sum_{i_3=1}^{28} \rho s_{i_3}.
\]
Accordingly, the point-wise density of the HCCI data tensor is
\begin{equation}\label{eq:density_pointwise}
        \Tn{D}_{\X}(i_1, i_2, t) = \sum_{{i_3}=1}^{28} \X(i_1, i_2, i_3, t) \; .
\end{equation}
We denote the full density tensor as $\Tn{D}_{\X} \in \Real^{I_1 \times I_2 \times \tau}$, and define our QoI
functional for mass conservation at each time step as the volume integral of density
\begin{equation}\label{eq:density_integral_error_per_time}
   \int_{\Omega} \rho ~d\Omega  \equiv  g_{1}(\X_t) = \sum_{i_1=1}^{I_1}\sum_{i_2=1}^{I_2} \Tn{D}_{\X}(i_1,i_2,t) \;,
   ~~t = 1,2,\dots, \tau,
\end{equation}
where the constant factor of $\Delta x \Delta y$ has been dropped for simplicity. 


\levelstay{Conservation of Kinetic Energy} The point-wise velocities in the $x$-direction and $y$-direction are
the last two quantities in the point-wise variables vector, and hence slices of the HCCI data tensor,
\begin{equation}\label{eq:velocity_x_pointwise}
        \Tn{U}^{(x)}_{\X}(i_1,i_2,t) = \X(i_1, i_2, 31, t) \; ,
\end{equation}
and
\begin{equation}\label{eq:velocity_y_pointwise}
        \Tn{U}^{(y)}_{\X}(i_1,i_2,t) = \X(i_1, i_2, 32, t) \;, 
\end{equation}
respectively. Kinetic energy is an important quantity in turbulent flows, and the point-wise kinetic energy
is given by
\begin{equation}\label{eq:kinetic_energy_pointwise}
        \K_{\X}(i_1,i_2,t) = \Tn{D}_{\X}(i_1,i_2,t) \left[\left(\Tn{U}^{(x)}_{\X}(i_1,i_2,t)\right)^2 + \left(\Tn{U}^{(y)}_{\X}(i_1,i_2,t)\right)^2 \right] \; .
\end{equation}
We denote the full kinetic energy tensor as $\K_{\X} \in \Real^{I_1 \times I_2 \times \tau}$, which is defined as
\begin{equation}\label{eq:kinetic_energy_pointwise2}
        \K_{\X} = \Tn{D}_{\X} \ast \left[\Tn{U}^{(x)}_{\X}  \ast  \Tn{U}^{(x)}_{\X} + \Tn{U}^{(y)}_{\X}  \ast  \Tn{U}^{(y)}_{\X}\right] \; ,
\end{equation}
where $\Tn{A} \ast \Tn{B}$ is the Hadamard (or element-wise) product of tensors $\Tn{A}$ and $\Tn{B}$~\cite{QiLu17}.
We now define our QoI functional for kinetic energy as
\begin{equation}\label{eq:kinetic_energy_integral_error_per_time}
        g_{2}(\X_t) = \sum_{i_1=1}^{I_1}\sum_{i_2=1}^{I_2} \Tn{K}_{\X}(i_1,i_2,t) \; ,
        ~~t = 1,2,\dots, \tau,
\end{equation}
representing conservation at each timestep.


\levelup{Results} 

In this section we present results of computing low-rank tensor decompositions on a subset of the HCCI data tensor. We focus on a subset of time steps of the simulation data where there exist interesting dynamics, specifically time steps 360 through 409 during an ignition phase. Thus, the HCCI data tensor is $\X\in\Real^{672 \times 672 \times 32 \times 50}$. As noted above in \cref{sec:hcci_background}, 
explicit scaling/unscaling of each variable slice is not necessary for this data.
Otherwise, we use the solution strategy described in \cref{sec:sol_procedure}, computing a traditional CP/Tucker decomposition via CP-ALS/ST-HOSVD using a supplied CP rank and ST-HOSVD error tolerance, respectively, which serves as the initial guess for the goal-oriented decomposition.  For CP-ALS, we set a tight stopping tolerance of $10^{-9}$ on the change in the tensor fit and a maximum of 100 iterations to ensure the relative tensor reconstruction error for the unscaled tensor data and the CP model is dominated by rank-truncation error rather than solver error.  We then choose the weighting coefficients as described above and run 20 iterations of the trust region-Newton optimization method as described in \cref{sec:stop_crit}.

\begin{figure*}[t]
	\centering
	\begin{subfigure}[b]{0.49\textwidth}
		\includegraphics[width=\textwidth]{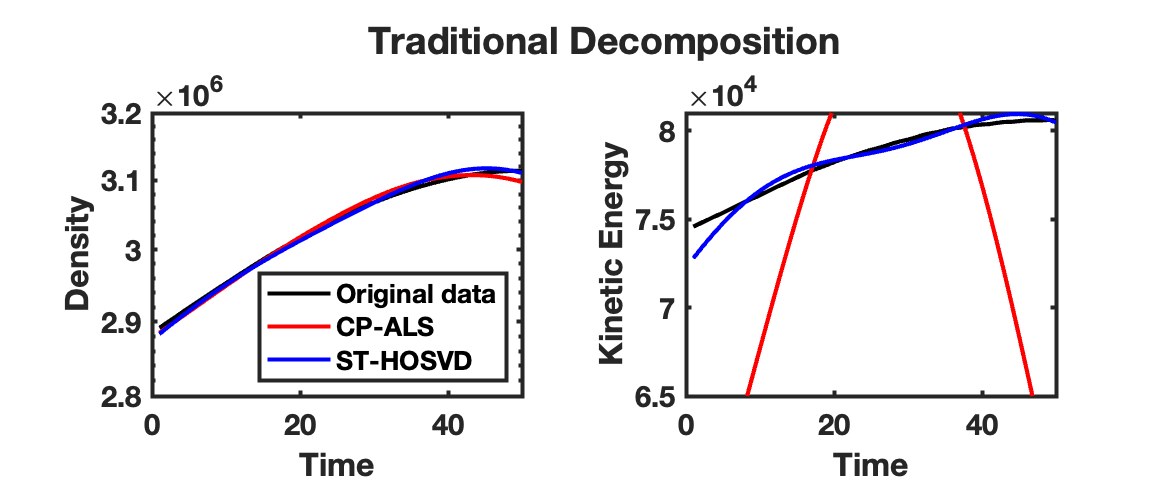}
		\caption{Initial density and kinetic energy QoIs}
		\label{fig:hcci_qois_initial}
	\end{subfigure}
	\hfill
	\begin{subfigure}[b]{0.49\textwidth}
		\includegraphics[width=\textwidth]{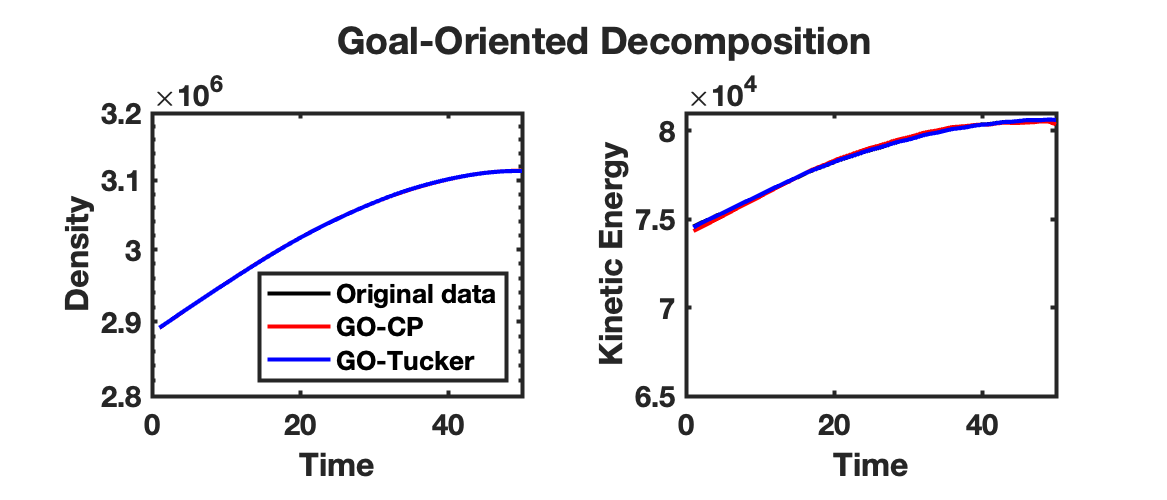}
		\caption{Final density and kinetic energy QoIs}
		\label{fig:hcci_qois_final}
	\end{subfigure}
	\caption{Traditional (CP-ALS/ST-HOSVD) and goal-oriented (GO-CP/GO-Tucker) decompositions of the HCCI data for CP rank of 70 and initial ST-HOSVD truncation tolerance of $0.1$ (corresponding to a core tensor of size $46\times 40\times 7 \times 3$), both resulting in compression ratios of about 7,000x.}
	\label{fig:hcci_qois}
\end{figure*}

The density and kinetic energy QoIs described above are displayed for the HCCI data tensor and reconstruction obtained from a rank-70 CP-ALS decomposition and ST-HOSVD decomposition with truncation tolerance of $0.1$ in~\cref{fig:hcci_qois_initial}. For the density QoI, there are small but clearly visible differences observed between the data and traditional CP-ALS and ST-HOSVD decompositions. However, for the kinetic energy QoI, the differences are much more pronounced, especially for the CP-ALS decomposition. Such results are what motivated the development of the goal-oriented decompositions: when computing decompositions that simply aim to match the reconstructed low-rank model to the data---i.e., as in traditional tensor decompositions---it is often the case that QoIs associated with the data are not well modeled. In comparison, computing the goal-oriented decompositions---both GO-CP and GO-Tucker---leads to much more accurate modeling of the density and kinetic QoIs, as shown in~\cref{fig:hcci_qois_final}. Furthermore, this increase in QoI modeling accuracy comes at minimal increase in the relative tensor reconstruction error---about $0.04\%$ for GO-CP and $0.09\%$ for GO-Tucker.

\begin{figure*}[t]
	\centering
	\includegraphics[width=\textwidth]{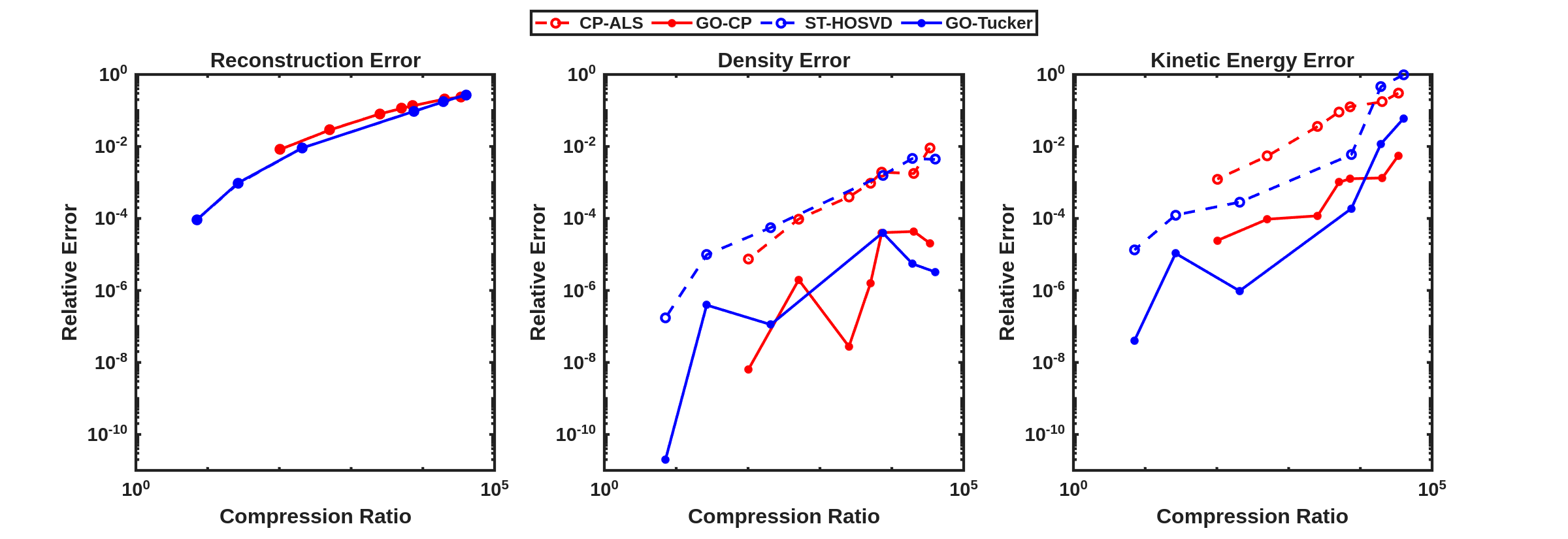}\\
	\caption{Relative tensor reconstruction and QoI errors for traditional (CP-ALS/ST-HOSVD) and goal-oriented (GO-CP/GO-Tucker) decompositions for the HCCI problem.  Significant improvement in the relative QoI errors using the goal-oriented decompositions are observed (solid versus dashed lines in the center and right plots) with little observed change in the relative tensor reconstruction errors (left plot).}
	\label{fig:hcci_cp_tucker_goal_error}
\end{figure*}

We next consider the performance of the goal-oriented approach compared to the traditional CP-ALS and ST-HOSVD approaches over a range of CP ranks and initial ST-HOSVD truncation tolerances.  \Cref{fig:hcci_cp_tucker_goal_error} displays the relative tensor reconstruction error as well as the relative QoI error for the density and kinetic energy QoIs as functions of compression ratio of the CP/Tucker reduced models, for both traditional and goal-oriented methods.  We make several observations:
\begin{itemize}
	\item The Tucker relative tensor reconstruction error is (slightly) less than the CP relative tensor reconstruction error for the same compression ratio, as is typically expected.
	\item The relative tensor reconstruction errors are essentially identical between the traditional and goal-oriented approaches, for both CP and Tucker. 
	\item In all compression ratios tested in these experiment, the relative QoI errors were improved using the goal-oriented approach. For the density QoI, we see around 2--4 orders of magnitude improvement in both CP and Tucker decompositions, while for the kinetic energy QoI, we see around 1--3 orders of magnitude improvement. 
\end{itemize}
	
\begin{figure}[h!]
	\centering
	\begin{subfigure}[b]{0.3\textwidth} 
		\includegraphics[width=\textwidth]{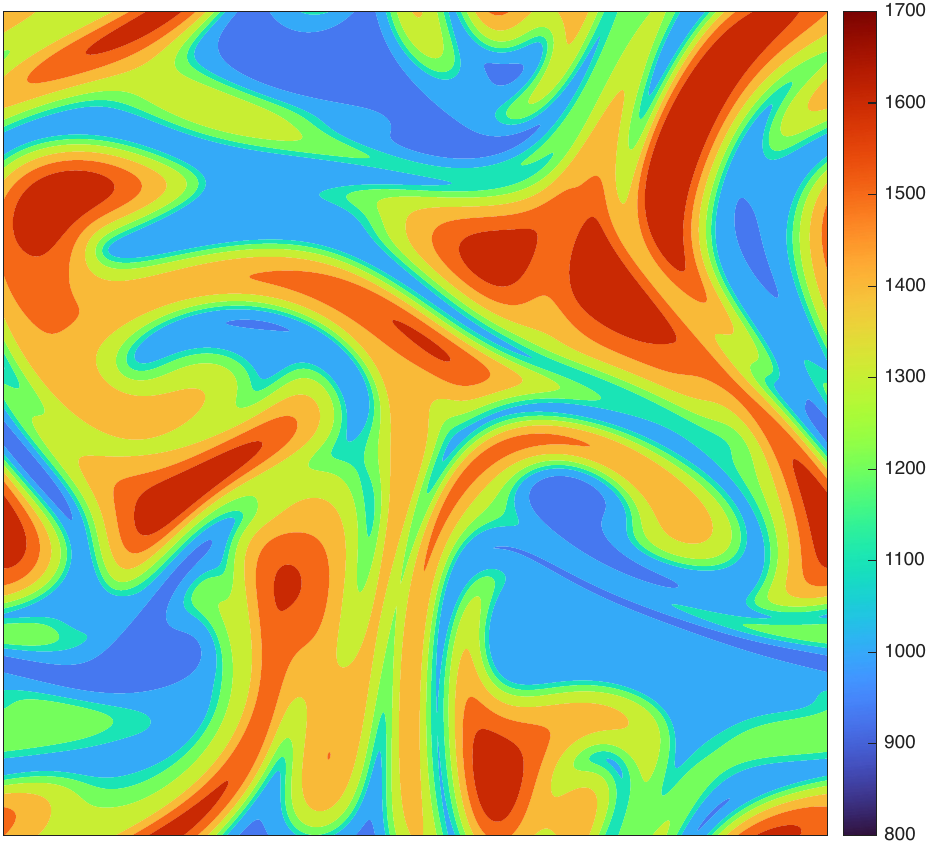}
		\caption{HCCI Data}
		\label{fig:hcci_data}
	\end{subfigure}
	\begin{subfigure}[b]{0.275\textwidth}
		\includegraphics[width=\textwidth]{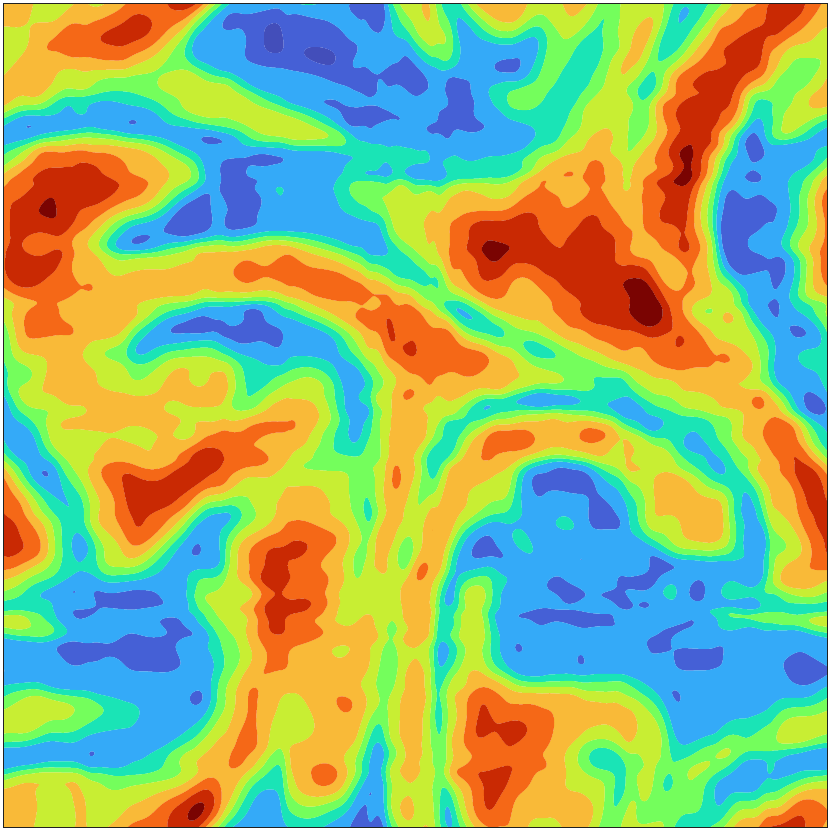}
		\caption{CP-ALS}
		\label{fig:hcci_cp}
	\end{subfigure}
	\begin{subfigure}[b]{0.275\textwidth}
		\includegraphics[width=\textwidth]{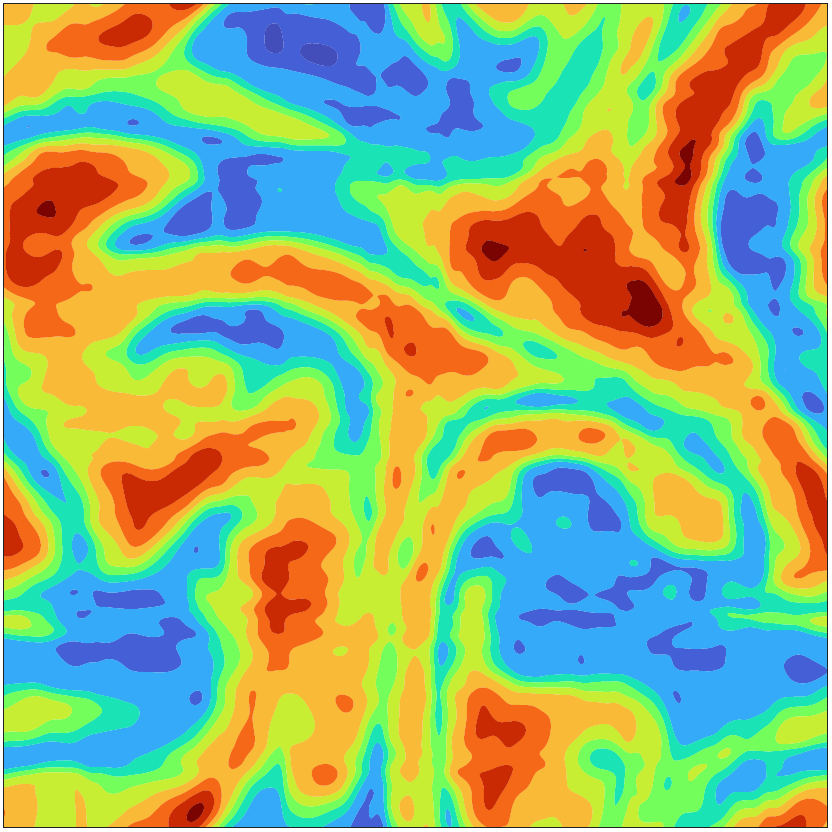}
		\caption{GO-CP}
		\label{fig:hcci_gocp}
	\end{subfigure} \\
	\begin{subfigure}[b]{0.575\textwidth} 
		\raggedleft 
		\captionsetup{justification=raggedleft,singlelinecheck=false} 
		\includegraphics[height=0.478\textwidth]{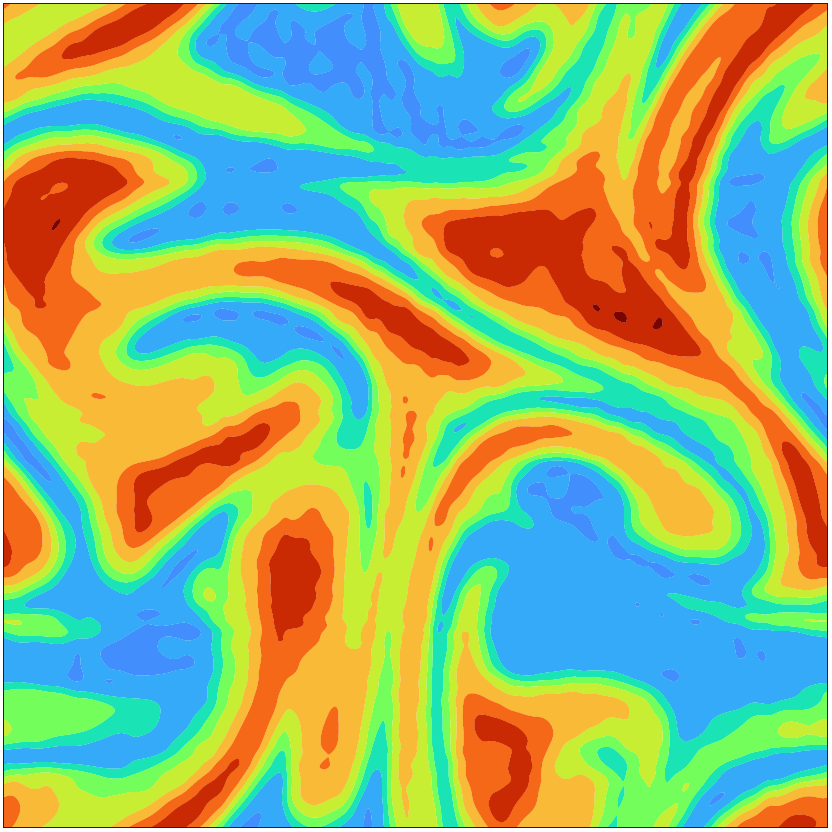} 
		\caption{ST-HOSVD\hspace*{0.1\textwidth}} 
		\label{fig:hcci_hosvd}
	\end{subfigure}
	\begin{subfigure}[b]{0.275\textwidth}
		\includegraphics[width=\textwidth]{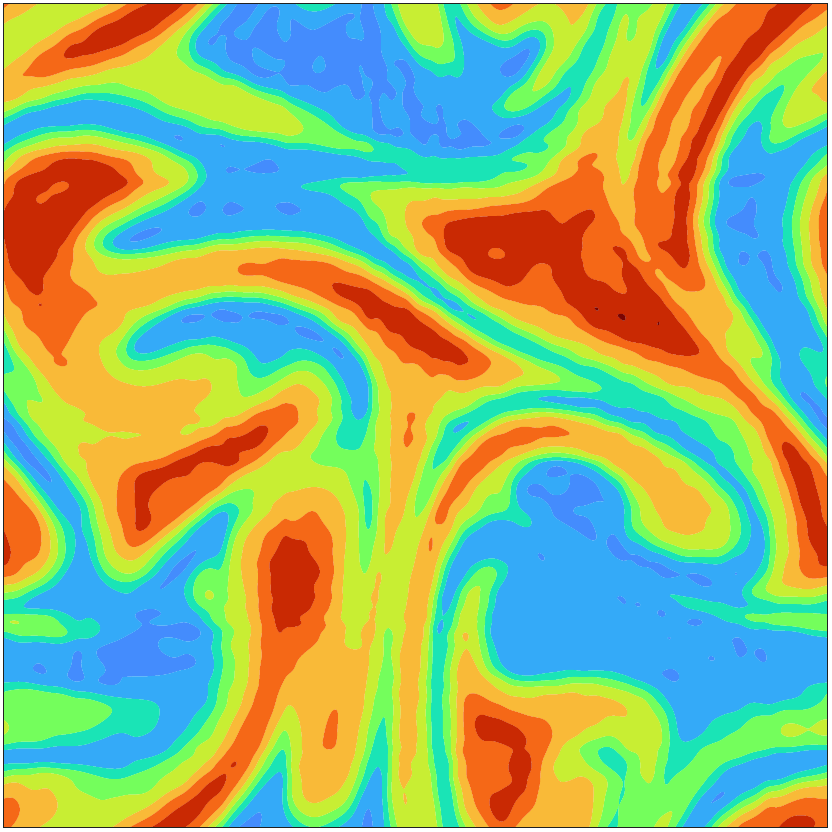}
		\caption{GO-Tucker}
		\label{fig:hcci_gotucker}
	\end{subfigure}
	\caption{Visual comparison of the (\subref{fig:hcci_data}) HCCI data and reconstructions from decompositions computed using (\subref{fig:hcci_cp}) CP-ALS, (\subref{fig:hcci_gocp}) GO-CP, (\subref{fig:hcci_hosvd}) ST-HOSVD, and (\subref{fig:hcci_gotucker}) GO-Tucker for the temperature variable over all spatial grids points at time step 401 out of 626. The decompositions use either a CP rank of 70 and initial ST-HOSVD truncation tolerance of $0.1$ (corresponding to a core tensor of size $46\times 40\times 7 \times 3$), resulting in compression ratios of about 7,000x for both the CP and Tucker decompositions.}
	\label{fig:hcci_temp}
\end{figure}

To illustrate the impact of the compression gained from computing the low-rank tensor decompositions of the HCCI data, \cref{fig:hcci_data} presents a slice of data for the temperature variable over all spatial grid points at time step 401 in the simulation. In \cref{fig:hcci_cp} and \cref{fig:hcci_gocp}, we see the corresponding reconstructions using the CP-ALS and GO-CP decompositions, respectively, noting there are no clearly visible differences, despite the relative QoI errors of GO-CP being dramatic improvements over those of CP-ALS. Similarly, in \cref{fig:hcci_hosvd} and \cref{fig:hcci_gotucker}, we see the corresponding reconstructions using the ST-HOSVD and GO-Tucker decompositions, respectively, noting also here are no clearly visible differences despite the dramatic improvements in the relative QoI errors using GO-Tucker.

%% file: sections/plasma.tex

\leveldown{Background}


The computational modeling of complex plasma physics systems is of critical importance in science and advanced technology.
A heavily used base-level model is resistive magnetohydrodynamics (MHD)  \cite{GoedbloedPoedts2004,bittencourt2013fundamentals}.  
This model is useful for describing the macroscopic dynamics of conducting fluids in the presence of electromagnetic fields and is 
often employed to study aspects of astrophysical phenomena (e.g., stellar interiors, solar 
flares), important science and technology applications (e.g., tokamak and stellarator devices, alternate pulsed fusion concepts), 
and basic plasma physics phenomena (e.g., magnetic reconnection, hydromagnetic instabilities)  \cite{GoedbloedPoedts2004,bittencourt2013fundamentals}.
The mathematical basis for the continuum modeling of these systems is the solution of the 
governing partial differential equations (PDEs) describing conservation of mass, momentum, and energy, augmented by an evolution equation for the magnetic field.
The magnetic field leads to Lorentz forces in the momentum equation, ${\bf j} \times {\bf B} = (\nabla \times {\bf B} / \mu_0) \times {\bf B}$, and a Joule heating energy dissipation term, $\eta \|{\bf j}\|^2$, in the energy equation. 
An outline of these equations is provided in \cref{eq: MHDSystem}. 
\begin{table}[t]
\begin{center}
\begin{tabular}{|c|c|}
\hline
\tiny \parbox{1.0in}{\scriptsize Mass} &  \scriptsize \parbox{4.25in} {\[ \frac{\partial \rho }{\partial t} + \nabla \cdot (\rho {\bf u}) = 0 \] } \\ \hline
\tiny \parbox{1.0in}{\scriptsize Momentum} & \scriptsize\parbox{4.25in} {\[ \frac{\partial {{\bf m}}} {\partial t} + \nabla \cdot \left[ {{\bf m} \otimes {\bf u}} + ( P{\bf I} + {{\bf \Pi}})\right] - {\bf j} \times {\bf B}  = {\bf 0};~{\bf m} \equiv
 \rho {\bf u} \]} \\ \hline
\tiny \parbox{1.0in}{\scriptsize Energy } &  \scriptsize\parbox{4.25in} {\[ \frac{\partial (\rho e)}{\partial t} + \nabla \cdot \left[ \rho {\bf u} e + {\bf q} \right] -  {{\bf \Pi} : {\nabla {\bf u}}} - \eta \|{\bf j}\|^2  = {\bf 0} \]} \\  \hline
\tiny \parbox{1.0in}{\scriptsize Magnetic Field } &  \scriptsize\parbox{4.25in} {\[ \frac{\partial {\bf B}}{\partial t} + \nabla \cdot \left[ {\bf u} \otimes {\bf B} - { \bf B} \otimes {\bf u} - \frac{\eta}{\mu_0} \left(\nabla {\bf B} - (\nabla {\bf B})^T \right) + \psi{\bf I}  \right] = {\bf 0};~  \nabla \cdot \mathbf{B} = {\bf 0}\]} \\  \hline \hline
\tiny \parbox{1.0in}{\scriptsize Plasma QoIs} &  \scriptsize \parbox{4.25in} {\[ {  IE = \int_\Omega (\rho e) d\Omega;~KE = \int_\Omega \frac{\| {\bf m} \|^2}{2 \rho} d\Omega;~ ME = \int_\Omega \frac{\| {\bf B} \|^2}{2\mu_0}  d\Omega;~M = \int_\Omega \|{\bf m}\|^2 d\Omega;  }\]} \\ \hline 
\end{tabular}
\end{center}
\caption{Low Mach Number Compressible Resistive MHD Plasma Model and QoIs.}
\label{eq: MHDSystem}
\end{table}

In the computational MHD simulations that are used in the results that follow, a fully-implicit variational multiscale (VMS)  
unstructured finite element (FE) discretization of the governing MHD system is employed. The fully-implicit time integration, based on RK methods, allow efficient solution of longer-time scale dynamics \cite{ShadidetalResistiveMHD2016,BonillaShadidetalVMSMHDMCF2023}
and the unstructured FE spatial discretization allows for complex geometries (e.g. ITER \cite{BonillaShadidetalVMSMHDMCF2023}) and local grid refinement. Simulations of transient, large, 3D simulations can  generate solutions
with $O(10M)$ -- $O(1B)$ elements and $O(10^3)$ -- $O(10^4)$ time-steps resulting in FE databases from $O(100GB)$ -- $O(10TB)$ in size.

\levelstay{Quantities of Interest}


In the pursuit of fundamental scientific understanding of complex plasma physics systems, accurate representations of the state (solution variables $(\rho, {\bf m}, (\rho e), {\bf B})$) as well 
a number of scientific QoIs are important. As an example, magnetic reconnection is a fundamental process whereby a magnetic field is structurally altered via some dissipation mechanism,
resulting in a rapid conversion of magnetic field energy into plasma energy and significant plasma transport/flow. 
Magnetic reconnection dominates the dynamics of many space and laboratory plasmas, and is at the root of phenomena such as solar flares, coronal mass ejections, 
and major disruptions in MCF energy (MFE) experiments \cite{GoedbloedPoedts2004,biskamp1994magnetic}.  In this context, there are a number of important QoIs
that include the internal plasma energy (IE), magnetic energy (ME), kinetic energy (KE), the total energy, $TE = IE + KE + ME$, along with the total momentum, $\sqrt{M}$ (see definitions in \cref{eq: MHDSystem}).
In the evolution of the state variables and QoIs that are considered here, this transformation of energy components (ME, KE) in magnetic reconnection will be evident\footnote{It should be noted that due to energy transfer through the boundary there is a small decrease in the total energy in example problems that are presented.}. Commonly in the evaluation of reconnection phenomena, an understanding of the growth rate of 
these transformations is of interest, and for example, the slope of $\sqrt{M(t)}$ provides such a metric. Alternatively, using $\sqrt{ME(t)}$, one could compute a decay rate for the transformation of magnetic energy.


\begin{figure} [t]
\centering
\includegraphics[width=.6\textwidth]{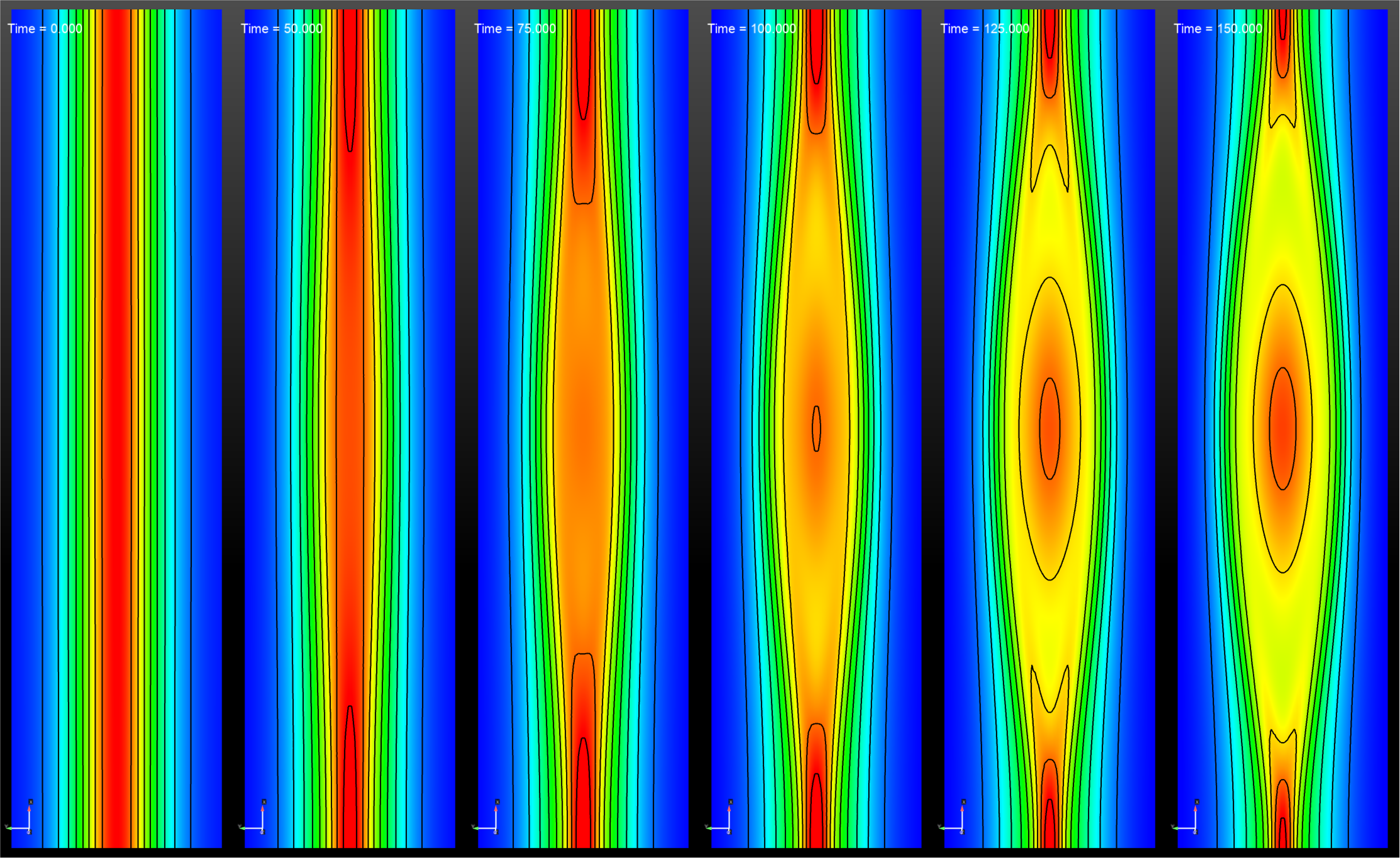}
\caption{Tearing mode evolution of unstable 1D current sheet to 2D magnetic island. 
Colored contours and isolines of the current $J_z$ are shown at times $t~=~0, 50, 75, 100, 125,150$.}
\label{fig:2dtearing}
\end{figure}

In this work, we attempt to preserve IE, KE, ME as quantities of interest which then implicitly preserve TE.  
Together with the momentum QoI, this results in four QoI functions that are approximated using the finite element discretization scheme described above.   Calculation of these QoIs is straightforward but substantially more complicated than the combustion QoIs described above, so descriptions of their evaluations are provided in~\cref{sec:app:plasma_qois}.

\levelstay{Results}
\label{sec:results}

\leveldown{2D Compressible Tearing Mode}

\begin{figure*}[t]
	\centering
	\begin{subfigure}[b]{0.49\textwidth}
		\includegraphics[width=\textwidth]{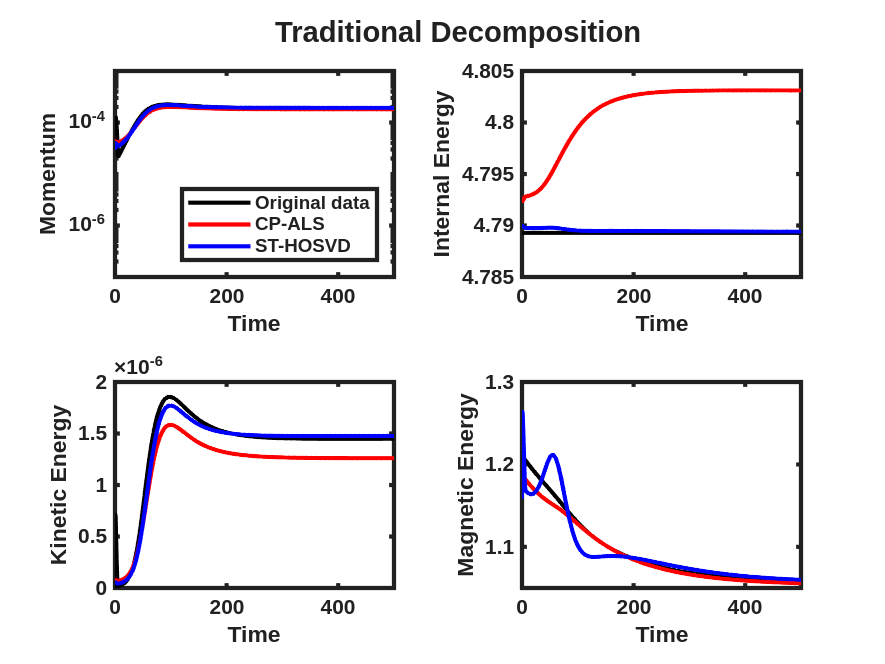}
		\caption{Initial momentum and energy QoIs}
		\label{fig:tearing_sheet_qois_initial}
	\end{subfigure}
	\hfill
	\begin{subfigure}[b]{0.49\textwidth}
		\includegraphics[width=\textwidth]{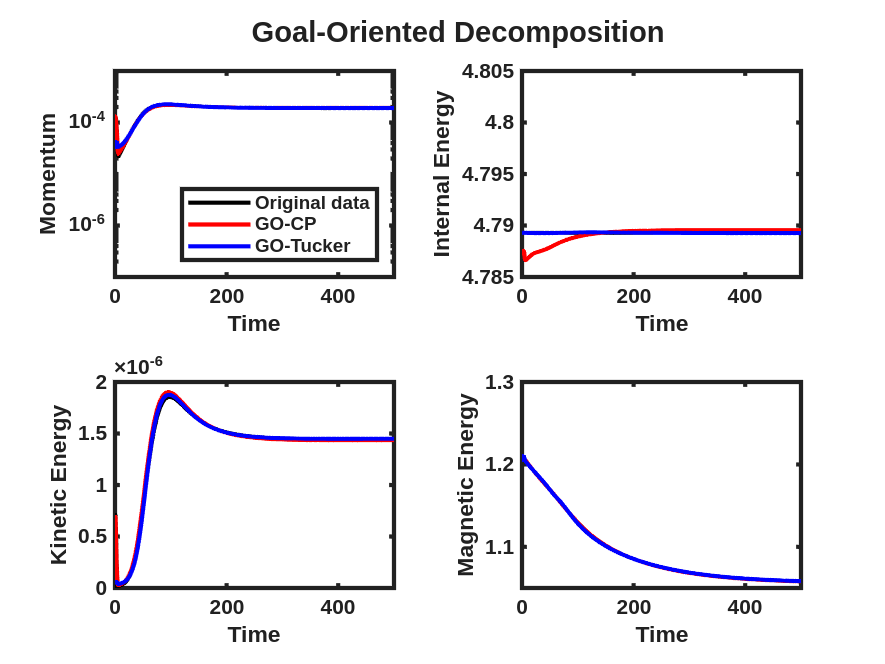}
		\caption{Final momentum and energy QoIs}
		\label{fig:tearing_sheet_qois_final}
	\end{subfigure}
	\caption{Traditional and goal-oriented CP and Tucker decomposition of the tearing sheet data for CP rank of 7 and initial ST-HOSVD truncation tolerance of $0.1$ (corresponding to a core tensor of size $8\times 9\times 5 \times 3$), both resulting in compression ratios of about 60,000x.}
	\label{fig:tearing_sheet_cp_qois}
\end{figure*}

The magnetic reconnection  problem is a 2D low flow-Mach number compressible tearing mode simulation that
 follows the  unstable evolution of a thin current sheet formed by a sheared magnetic field (Harris sheet) within an initially stationary velocity field.
 The domain is a rectangle $[0,4]\times[0,1]$ discretized on a uniform grid resulting in 501 and 201 grid points in the $x$ and $y$ directions, respectively.  The computation is for a Lundquist number of $10^3$ and consists of 500 time steps. Details of the full problem setup can be found in~\cite{chacon-pop-08-3dmhd,BonillaShadidetalVMSMHDMCF2023}.
The thin current sheet becomes unstable and forms a magnetic island structure as presented in~\cref{fig:2dtearing} with the $x$-axis oriented in the vertical direction.
The rate of the  decay of the $L_2$ norm of the magnetic energy and the growth of the $L_2$ norm of momentum provide information of the time-scale of the linear  phase of the
instability~\cite{chacon-pop-08-3dmhd}.  These rates correspond to exponential fits of the slopes of these QoIs in~\cref{fig:tearing_sheet_cp_qois}.

We now consider goal-oriented CP and Tucker decomposition of the 2-D tearing mode plasma physics data stored as a 4-way tensor $\X\in\Real^{401\times 201 \times 12 \times 501}$ consisting of each of the 12 simulation variables $(\rho, T, {\bf m}, {\bf u}, (\rho e), {\bf B})$ at each grid and time point.  We use the solution strategy described in \cref{sec:sol_procedure} where we first center each variable slice by its mean and scale it by its standard deviation and then compute a traditional CP/Tucker decomposition via CP-ALS/ST-HOSVD using a supplied CP rank and ST-HOSVD error tolerance, respectively, which serves as the initial guess for the goal-oriented decomposition.  As for the HCCI problem, we set a CP-ALS stopping tolerance of $10^{-9}$.  Since the momentum values are near zero for the first few time steps, we only include time steps in the momentum QoI starting from the fifth time step---i.e., $\mathbb{T}_q = \{5,...,\tau\}$ for the momentum QoI.  We then choose the weighting coefficients as described above and run 20 iterations of the trust region-Newton optimization method as described in \cref{sec:stop_crit}.


The momentum and energy QoIs described above are displayed for the data tensor and reconstruction obtained from a rank-7 CP-ALS decomposition and ST-HOSVD decomposition with truncation tolerance of $0.1$ in~\cref{fig:tearing_sheet_qois_initial}, where small but clearly visible differences are observed in these quantities, even though the relative tensor reconstruction error of the unscaled data and model is approximately $6\cdot 10^{-4}$ for CP-ALS and $1\cdot 10^{-4}$ for ST-HOSVD.  These same QoIs are then displayed after the goal-oriented decomposition as described above using the trust-region Newton optimization method in~\cref{fig:tearing_sheet_qois_final}, where these visual differences between the QoIs for the data and the low-rank models have largely disappeared.  Furthermore, the relative tensor reconstruction error is again about $6\cdot 10^{-4}$ and $1\cdot 10^{-4}$, respectively, indicating substantial reduction in relative QoI error with no increase in the relative tensor reconstruction error.


\begin{figure*}[t]
	\centering
	\includegraphics[width=\textwidth]{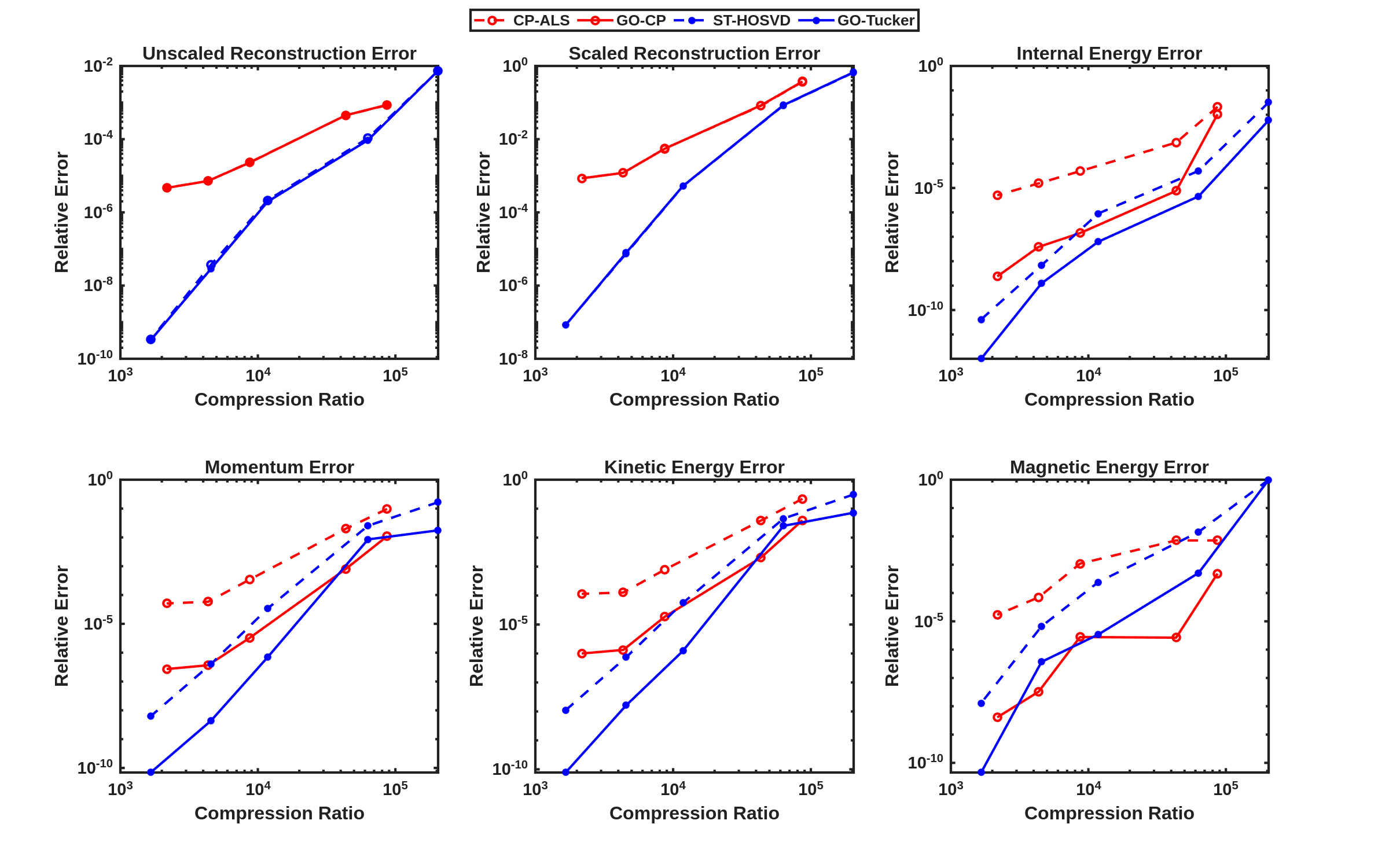}
	\caption{Relative tensor reconstruction and QoI errors for traditional and goal-oriented CP and Tucker decompositions, including both scaled and unscaled data and models for the 2D tearing mode problem.  Significant improvement in the relative QoI errors are observed with little observed change in the relative tensor reconstruction errors.}
	\label{fig:tearing_sheet_cp_tucker_goal_error}
\end{figure*}

We now consider the performance of the goal-oriented approach compared to the traditional CP/Tucker approach over a range of CP ranks and initial ST-HOSVD truncation tolerances.  \Cref{fig:tearing_sheet_cp_tucker_goal_error} displays the relative tensor reconstruction error (for both the scaled and unscaled data and models) as well as the relative QoI error for the momentum and energy QoIs as functions of compression ratio of the CP/Tucker reduced models, for both traditional and goal-oriented methods.  We make several observations:
\begin{itemize}
  \item The Tucker relative tensor reconstruction error is typically less than the CP relative tensor reconstruction error for the same compression ratio, as is typically expected.
  \item The relative QoI errors for Tucker is also typically less than for CP.
  \item The relative tensor reconstruction errors for the traditional and goal-oriented approaches are essentially identical, for both CP and Tucker.  This is true for both the scaled and unscaled data and models, indicating the goal improvements are not due to increased relative tensor reconstruction error of small variables relative to larger variables.
  \item Around 2--3 orders of magnitude improvement in relative QoI error for all QoIs is typically observed for the goal-oriented CP method, whereas goal-oriented Tucker results in about 1--2 orders of magnitude improvement.
\end{itemize}

\begin{figure*}[t]
	\centering
	\includegraphics[width=0.6\textwidth]{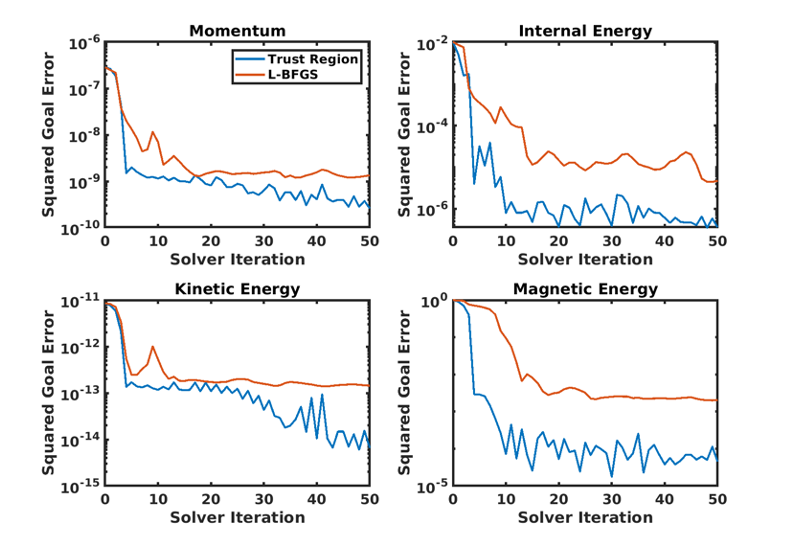}
	\caption{Convergence history for the trust-region Newton and L-BFGS optimization methods for goal-oriented Tucker. for the 2D tearing mode problem.}
	\label{fig:opt_conv_iter}
\end{figure*}

The above calculations all relied on use of the trust-region Newton optimization method with Gauss-Newton Hessian approximation and approximate block-diagonal preconditioner.  We found this approach to be much more reliable than the L-BFGS method, often finding more accurate solutions (in terms of QoI error) in significantly fewer iterations.  For example, \cref{fig:opt_conv_iter} displays the convergence history for the trust-region Newton and L-BFGS methods for the above goal-oriented Tucker experiment with HOSVD truncation tolerance of $0.1$, where we see lower achieved QoI error for the trust-region Newton method, yielding a good solution in just a few optimization iterations.

\begin{figure} [htbp]
\centering
\includegraphics[width=\textwidth]{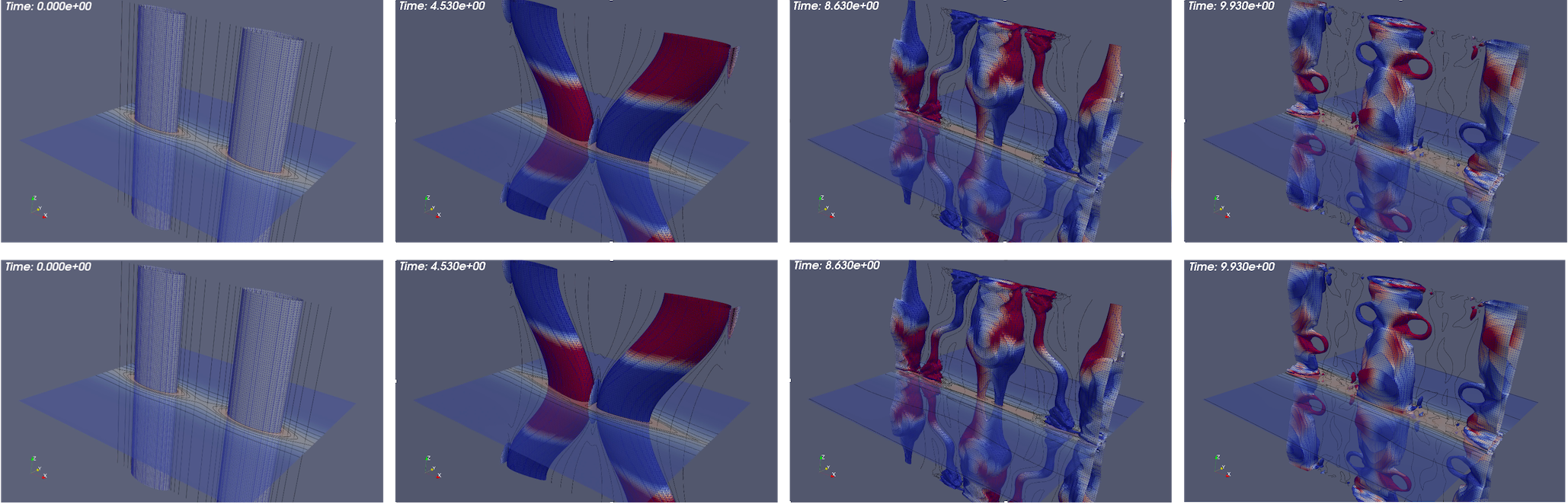}
\caption{Evolution of unstable 3D current flux tubes with magnetic reconnection.  Top row FE solution/data and bottom ST-HOSVD tensor reconstruction at 10x reduction.
Density iso-surface at 1.5 colored by $\bm{u}_x$, and slice plane colored by density at $t~=~0, 4.530, 8.630, 9.930$.}
\label{fig:3D_IC_HOSVD_10x}
\end{figure}

\begin{figure} [htbp]
\centering
\includegraphics[width=\textwidth]{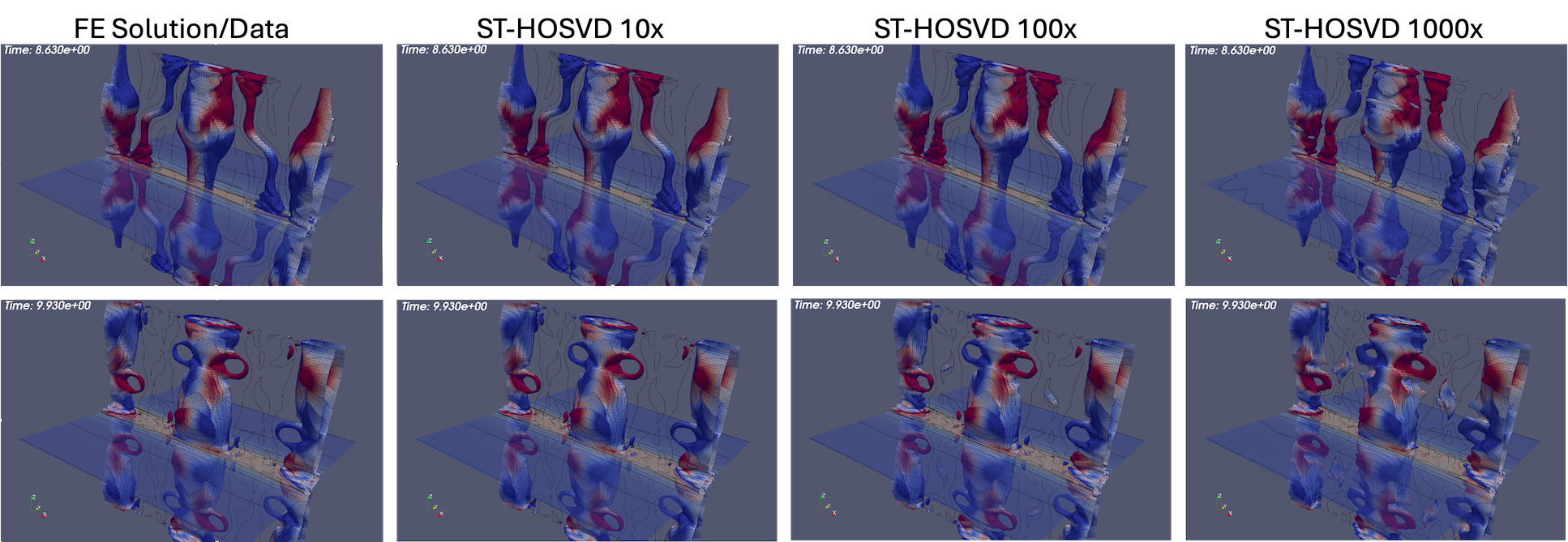}
\caption{Comparison of FE solution/data and ST-HOSVD reconstructions for compressions at $t~=8.630, 9.930$.}
\label{fig:3D_IC_Compare_HOSVD}
\end{figure}

\levelstay{3D Island Coalescence}
The island coalescence problem follows the unstable evolution of two 3D current flux tubes (which in the cross plane become islands) embedded in a  sheared magnetic 
field Harris sheet and is described in detail in \cite{ShadidetalResistiveMHD2016}. The computational simulations presented here are for a higher 
Lundquist number that exhibit more significant dynamics in the evolution and reconnection \cite{chacon-pop-08-3dmhd,ShadidetalResistiveMHD2016}.
This example allows a demonstration of compression and evaluation of the evolution of appropriate QoIs, but also the more complex dynamics
facilitate qualitative discussion of the reconstructions using low-rank tensor decompositions at various compression ratios shown in \cref{fig:3D_IC_HOSVD_10x} and \cref{fig:3D_IC_Compare_HOSVD}. 
Visualization of an iso-surface (in this case density) is a challenging qualitative metric for comparison of the reconstructions as small changes in the density distribution in space will modify the structure of the iso-surface at a given value.
\cref{fig:3D_IC_HOSVD_10x} compares the initial condition and three discrete times in the evolution of the reconnecting flux tubes for the original FE solution and 
ST-HOSVD decomposition at approximately 10x compression. Clearly evident is the quality of the ST-HOSVD reconstruction at 10x reduction.  

 In~\cref{fig:3D_IC_Compare_HOSVD}, the ST-HOSVD reconstruction for up to three orders of magnitude of reduction is visualized.
 Here it is evident that at the significant reduction of 100x, the reconstruction appears very representative of the original data for visualizing the structure of the iso-surface.
 Additionally, even up to a 1000x reduction, the major structural features of the iso-surface resemble reasonably well the original FE data with the existence of the 
 loop structures on the three main island iso-surfaces at time $t = 9.930$. 

\begin{figure*}[htbp]
	\centering
	\hfill
	\begin{subfigure}[b]{0.3\textwidth}
		\includegraphics[width=\textwidth]{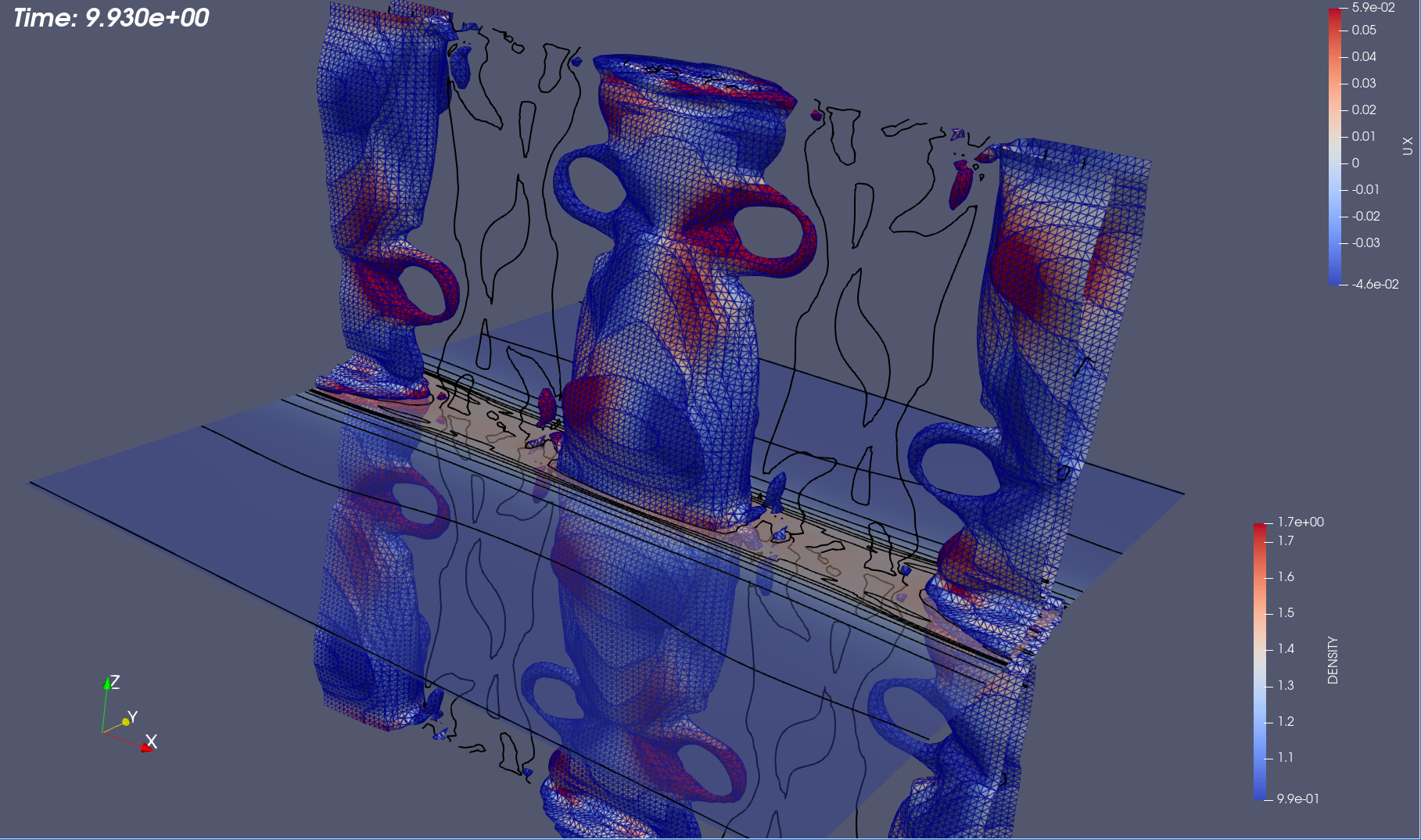}
		\caption{FE Solution/Data.}
		\label{fig:3D_IC_Reference_FE}
	\end{subfigure}
	\hfill
	\begin{subfigure}[b]{0.315\textwidth}
		\includegraphics[width=\textwidth]{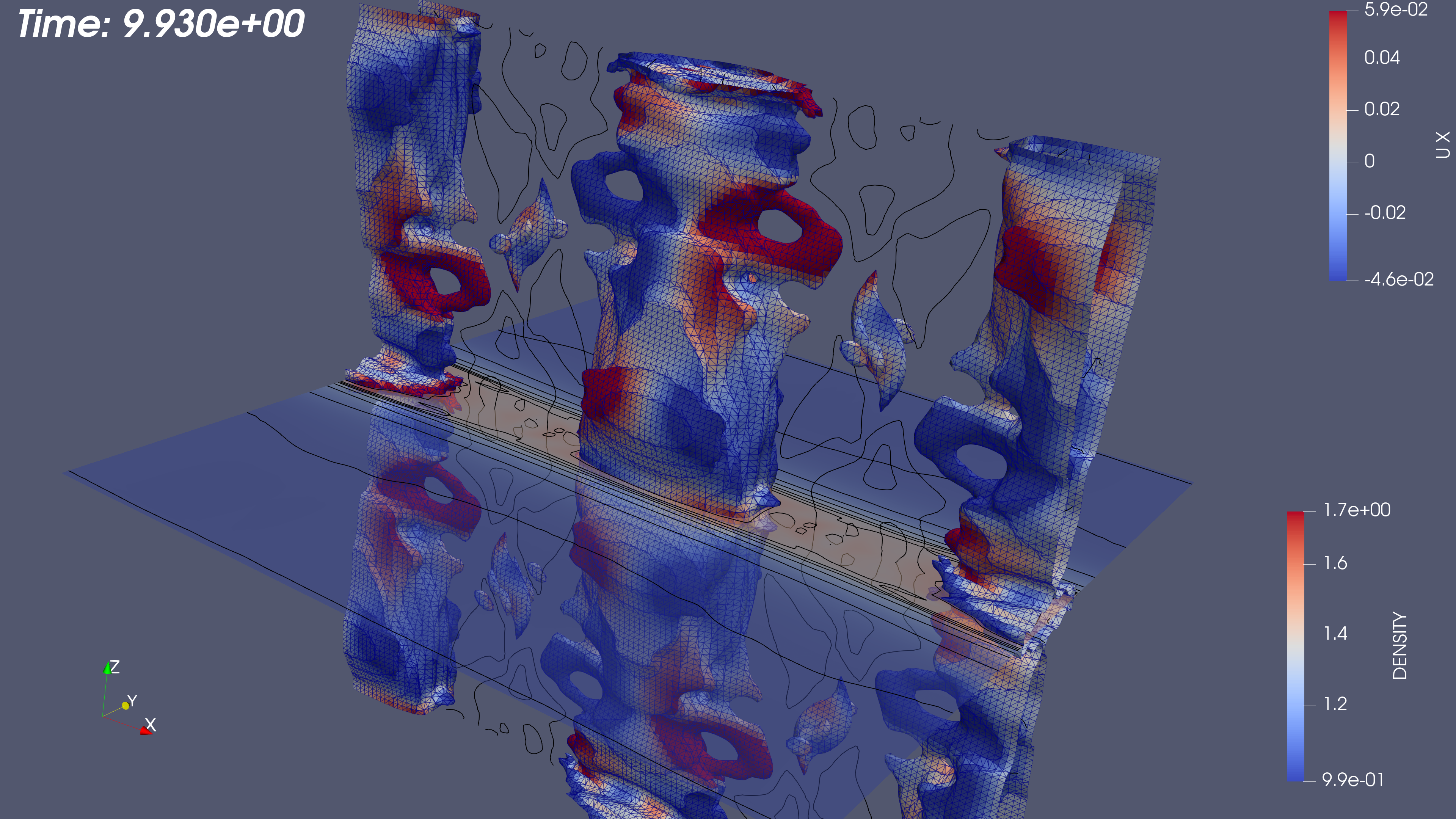}
		\caption{ST-HOSVD full-database.}
		\label{fig:3D_IC_HOSVD_full_database}
	\end{subfigure}
    \hfill
	\begin{subfigure}[b]{0.32\textwidth}
		\includegraphics[width=\textwidth]{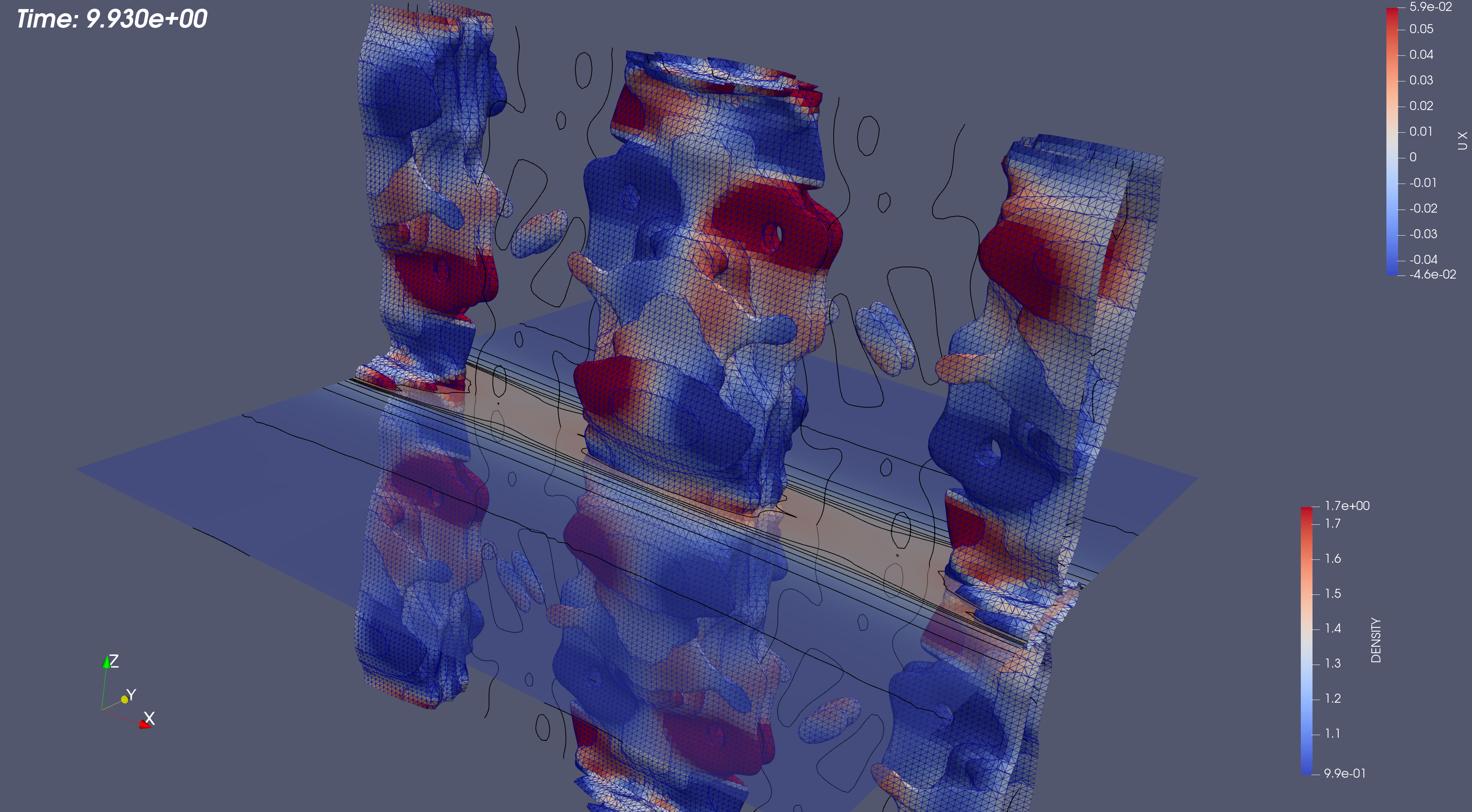}
		\caption{ST-HOSVD coarse-database.}
		\label{fig:fig:3D_IC_HOSVD_coarse_database}
	\end{subfigure}
    \hfill
	\caption{Degradation of reconstruction quality due to coarsening the 3D IC database to every third time step for the density iso-surface at 1.5 colored by $\bm{u}_x$, and slice plane colored by density at t = 9.930. This coarse-database is used for
	the study of the goal-oriented methods.}
	\label{fig:3D_IC_EffectCoarseDatabase_Reconstruction}
\end{figure*}
%
%

Next, a comparison of the reconstructions of the four plasma physics QoIs along with the quality of the iso-surface reconstruction for traditional and goal-oriented decompositions with varying levels of reduction is studied. 
The setup of the goal-oriented approach is identical to the tearing sheet problem described above, except the default stopping tolerance of $10^{-4}$ for CP-ALS is used for the initial guess for goal-oriented CP.
In this study, the original island coalescence data set produced by the simulation code was coarsened in time by selecting every third time-step to overcome memory and runtime limitations with our current serial, Matlab-based software implementation\footnote{These limitations are expected to be removed in a future software implementation.}. It should be noted that, as expected, the temporal coarsening somewhat adversely affects the general quality of the reconstructions over what would have been obtained with the original data set for the same level of reduction, for both the traditional and goal-oriented approaches. This degradation is illustrated in \cref{fig:3D_IC_EffectCoarseDatabase_Reconstruction} where the reference FE data, the full-database ST-HOSVD and coarsened-database reconstructions at a reduction of 1,000x are presented.

\begin{figure*}[ht!]
	\centering
	\hfill
	\begin{subfigure}[b]{0.49\textwidth}
		\includegraphics[width=\textwidth]{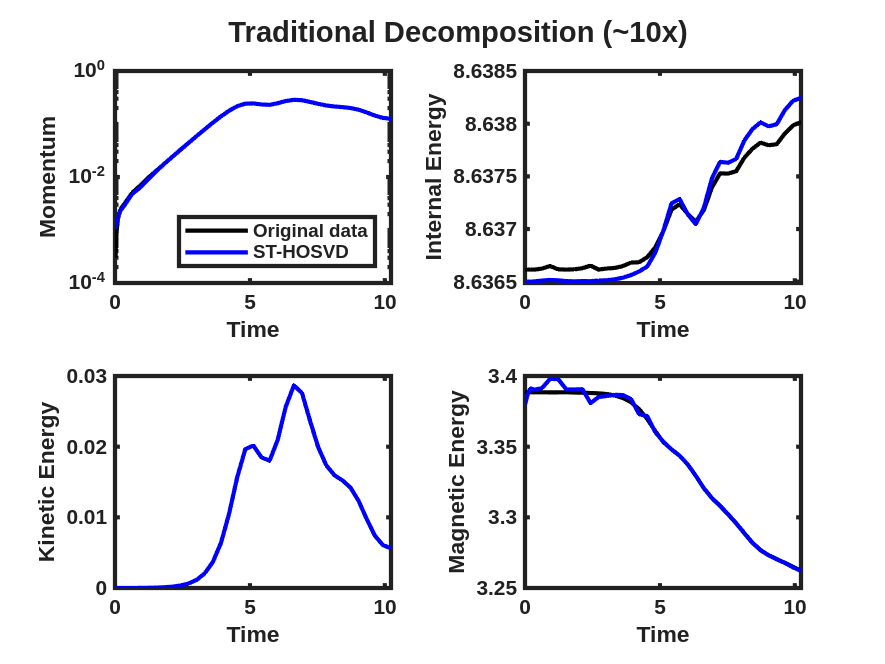}
		\caption{Relative tensor reconstruction errors of $0.00087$/$0.0083$ for unscaled/scaled ST-HOSVD.}
		\label{fig:island_coalescence_qois_initial_10}
	\end{subfigure}
	\hfill
	\begin{subfigure}[b]{0.49\textwidth}
		\includegraphics[width=\textwidth]{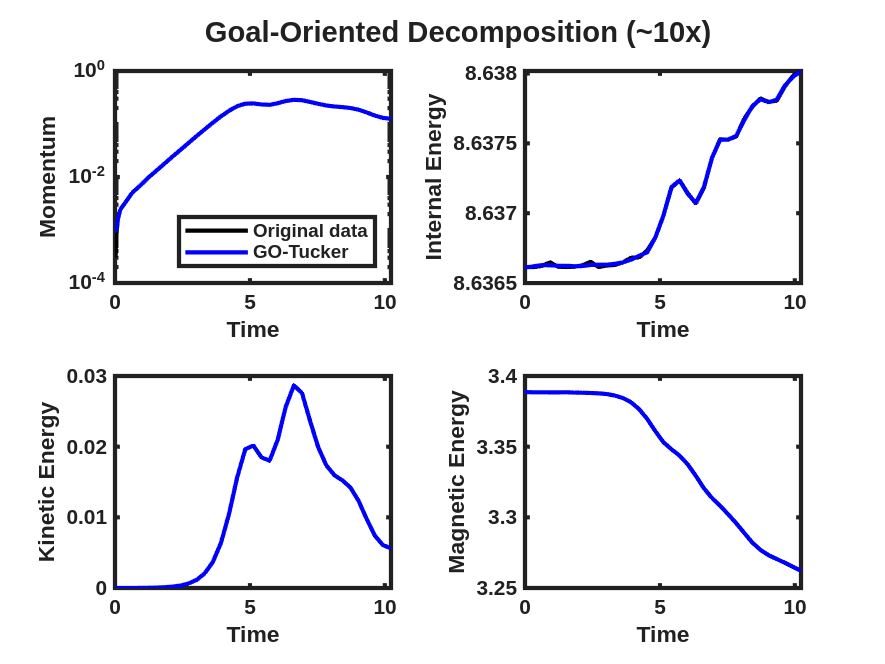}
		\caption{Relative tensor reconstruction errors of $0.00085$/$0.0084$ for unscaled/scaled GO-Tucker.}
		\label{fig:island_coalescence_qois_final_10}
	\end{subfigure}
	\hfill
	\begin{subfigure}[b]{0.49\textwidth}
		\includegraphics[width=\textwidth]{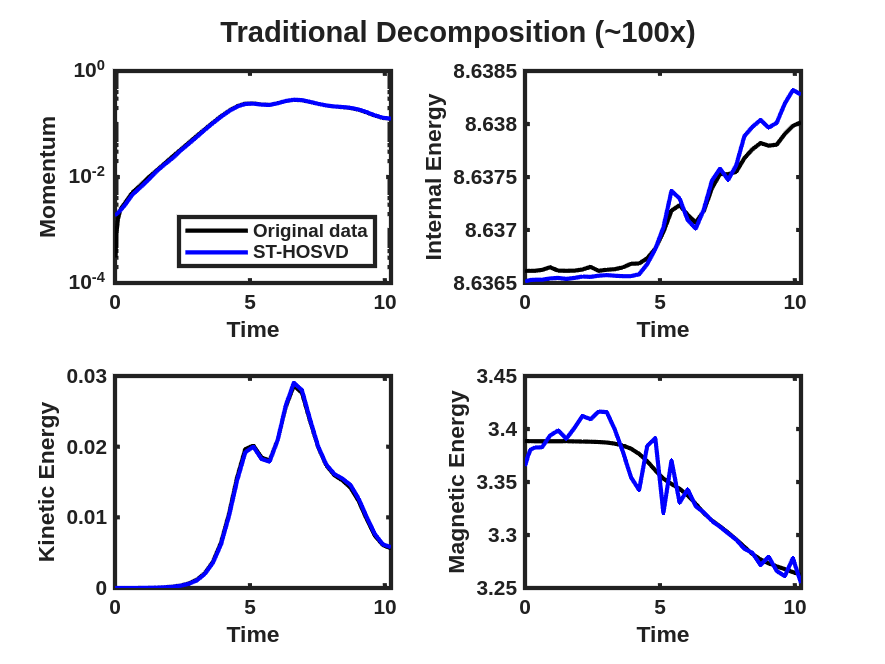}
		\caption{Relative tensor reconstruction errors of $0.0062$/$0.062$ for unscaled/scaled ST-HOSVD.}
		\label{fig:island_coalescence_qois_initial_100}
	\end{subfigure}
	\begin{subfigure}[b]{0.49\textwidth}
		\includegraphics[width=\textwidth]{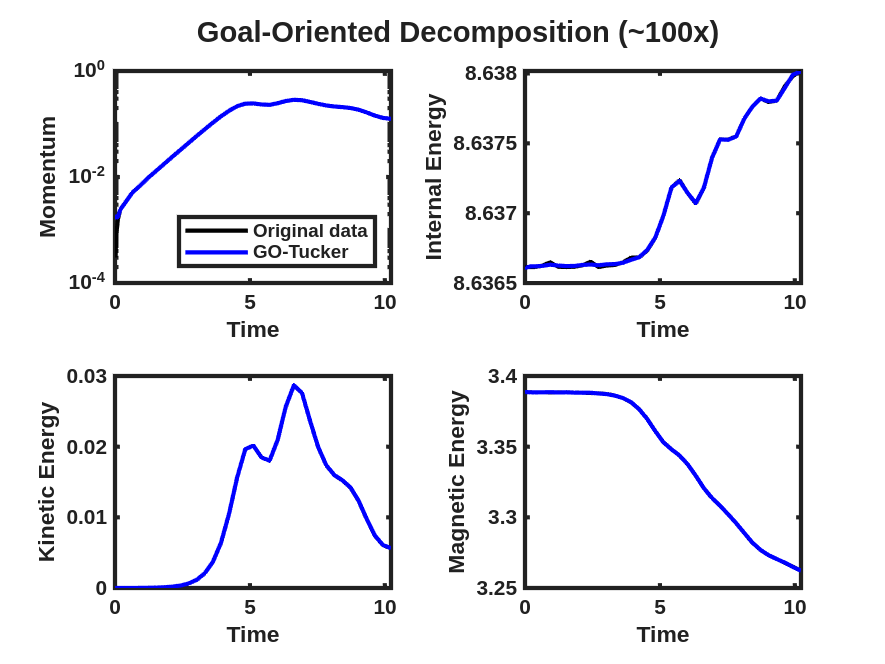}
		\caption{Relative tensor reconstruction errors of $0.0063$/$0.063$ for unscaled/scaled GO-Tucker.}
		\label{fig:island_coalescence_qois_final_100}
	\end{subfigure}
    \hfill
	\caption{Traditional (left column) and goal-oriented (right column) Tucker decomposition of the island coalescence data for with varying scaled initial ST-HOSVD tolerances ($0.01$, $0.07$), yielding corresponding compression ratios of approximately 10x and 100x,  respectively (CP decompositions for 10x and 100x compressions require too much memory due to the serial MATLAB implementation and thus are not provided).  Even with the relatively high accuracy of these decompositions, visual differences in the QoIs are observed for the ST-HOSVD decomposition which are eliminated with the goal-oriented approach while the unscaled and scaled relative tensor reconstruction errors are nearly the same between the approaches.}
	\label{fig:island_coalescence_qois1}
\end{figure*}

\begin{figure*}[ht!]
	\centering
	\hfill
	\begin{subfigure}[b]{0.49\textwidth}
		\includegraphics[width=\textwidth]{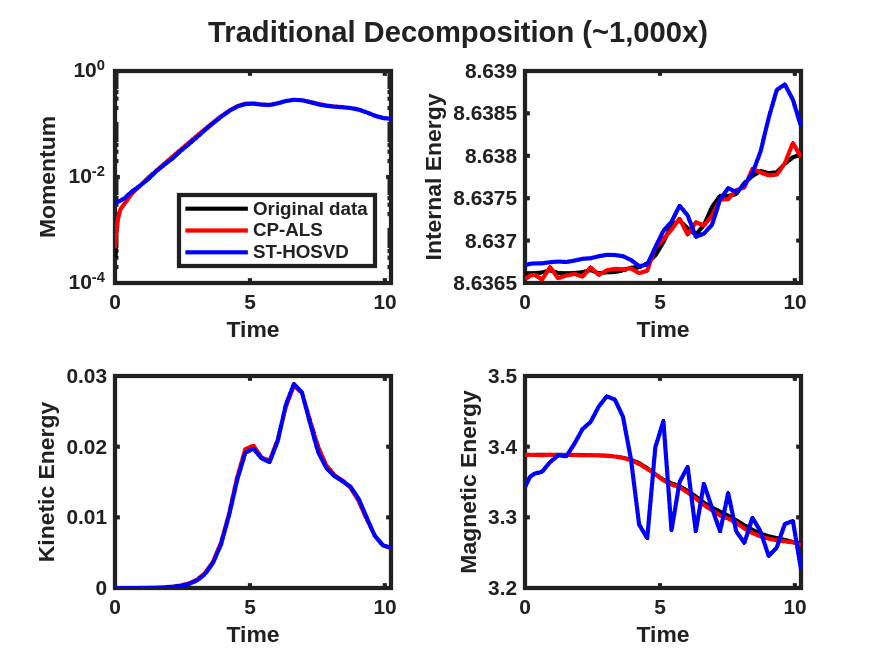}
		\caption{Relative tensor reconstruction errors of $0.011$/$0.053$ for unscaled/scaled CP-ALS and $0.012$/$0.13$ for unscaled/scaled ST-HOSVD.}
		\label{fig:island_coalescence_qois_initial_1000}
	\end{subfigure}
	\hfill
	\begin{subfigure}[b]{0.49\textwidth}
		\includegraphics[width=\textwidth]{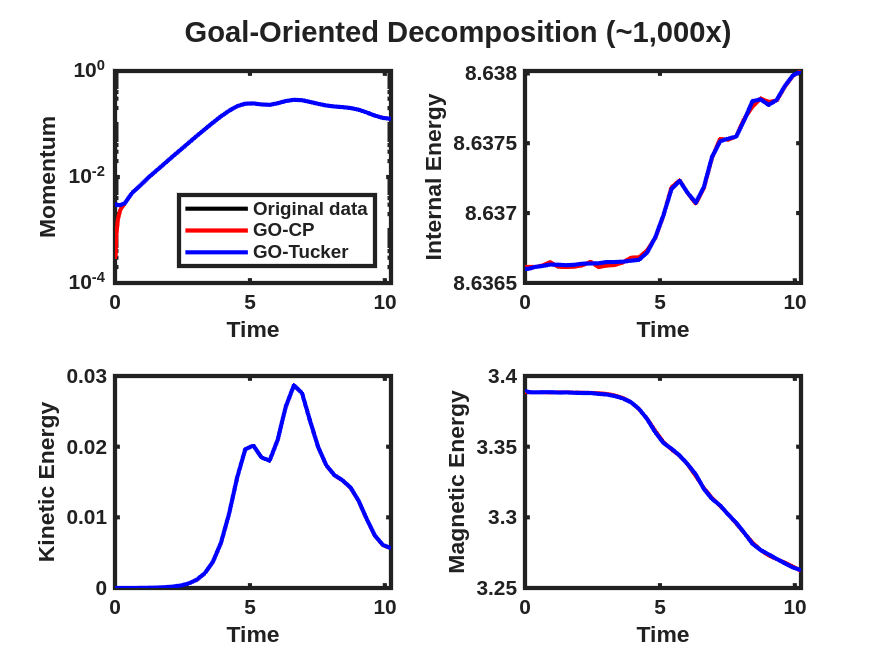}
		\caption{Relative tensor reconstruction errors of $0.011$/$0.053$ for unscaled/scaled GO-CP and $0.012$/$0.13$ for unscaled/scaled GO-Tucker.}
		\label{fig:island_coalescence_qois_final_1000}
	\end{subfigure}
	\hfill
	\begin{subfigure}[b]{0.49\textwidth}
		\includegraphics[width=\textwidth]{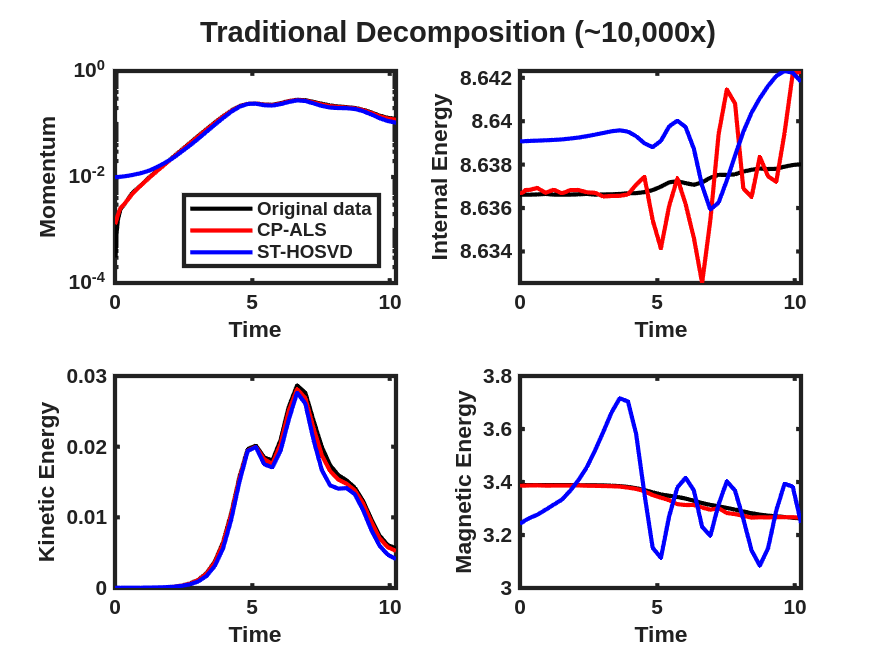}
		\caption{Relative tensor reconstruction errors of $0.024$/$0.18$ for unscaled/scaled CP-ALS and $0.027$/$0.29$ for unscaled/scaled ST-HOSVD.}
		\label{fig:island_coalescence_qois_initial_10000}
	\end{subfigure}
	\hfill
	\begin{subfigure}[b]{0.49\textwidth}
		\includegraphics[width=\textwidth]{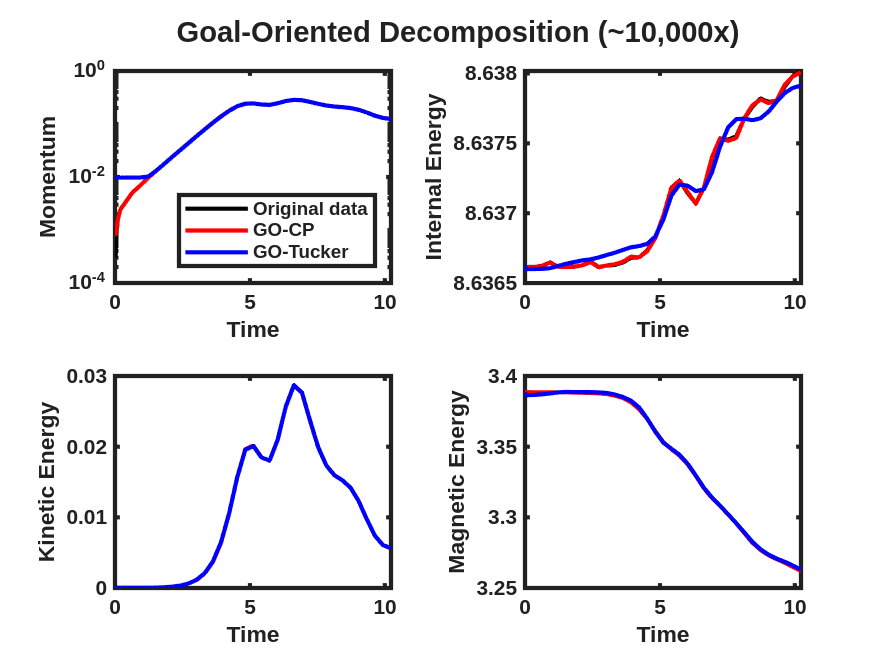}
		\caption{Relative tensor reconstruction errors of $0.024$/$0.18$ unscaled/scaled GO-CP and $0.025$/$0.29$ unscaled/scaled GO-Tucker.}
		\label{fig:island_coalescence_qois_final_10000}
	\end{subfigure}
	\hfill
	\caption{Traditional (left column) and goal-oriented (right column) CP and Tucker decomposition of the island coalescence data for with varying CP ranks (2294, 229) and scaled initial ST-HOSVD tolerances ($0.15$, $0.33$), yielding corresponding compression ratios of approximately 1,000x and 10,000x respectively.  Substantial errors in the QoIs are observed for both the original CP and Tucker decompositions due to the high compression, yet these errors are effectively removed for the goal-oriented approach.  As before,  the unscaled and unscaled relative tensor reconstruction errors are nearly the same between the two approaches.}
	\label{fig:island_coalescence_qois2}
\end{figure*}

\begin{figure*}[h!]
	\centering
	\begin{minipage}[b]{0.6\textwidth}
	\centering
		\includegraphics[width=\textwidth]{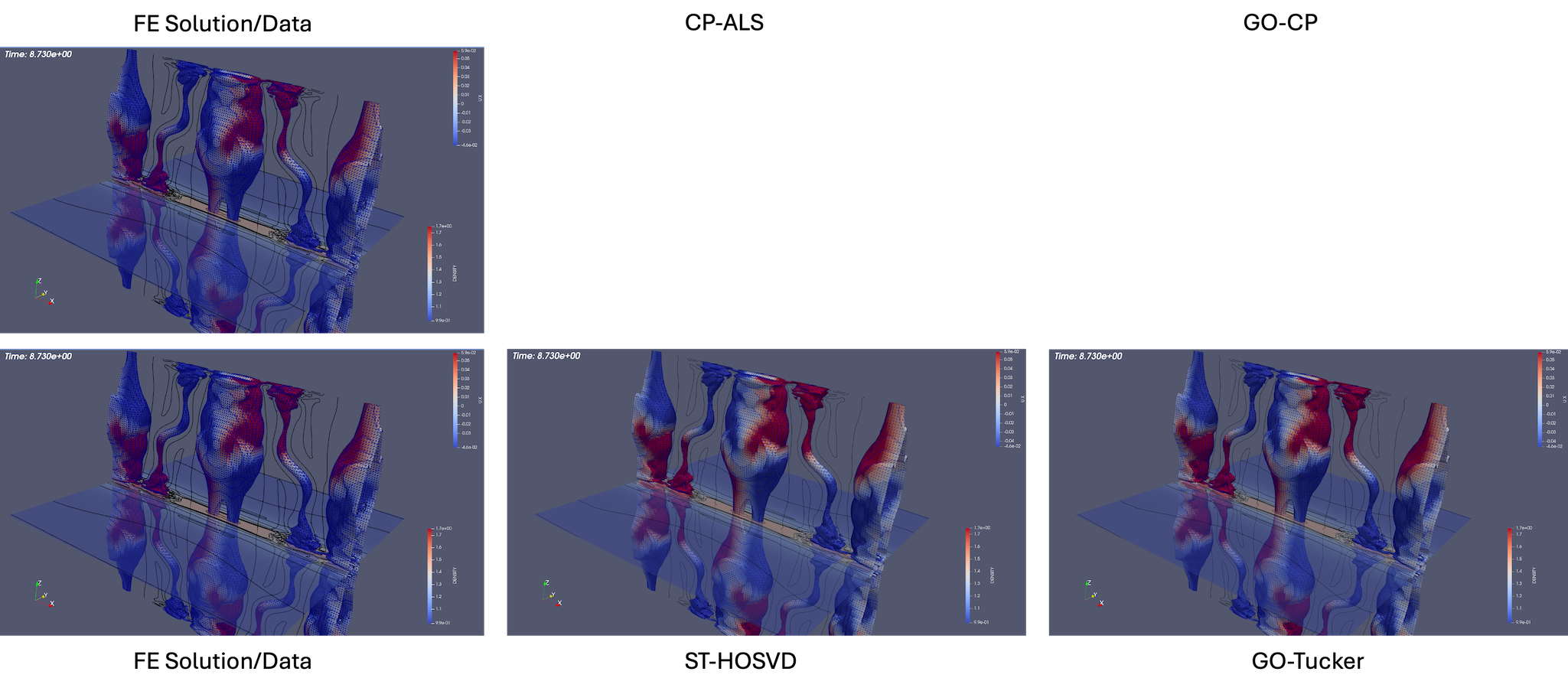}
		\caption*{10x}
			\label{fig:3D_IC_EffectCoarseDatabase_Reconstruction_10x_goal_oriented_t_8.730}
	\end{minipage}
	\hfill
	\begin{minipage}[b]{0.6\textwidth}
	\centering
		\includegraphics[width=\textwidth]{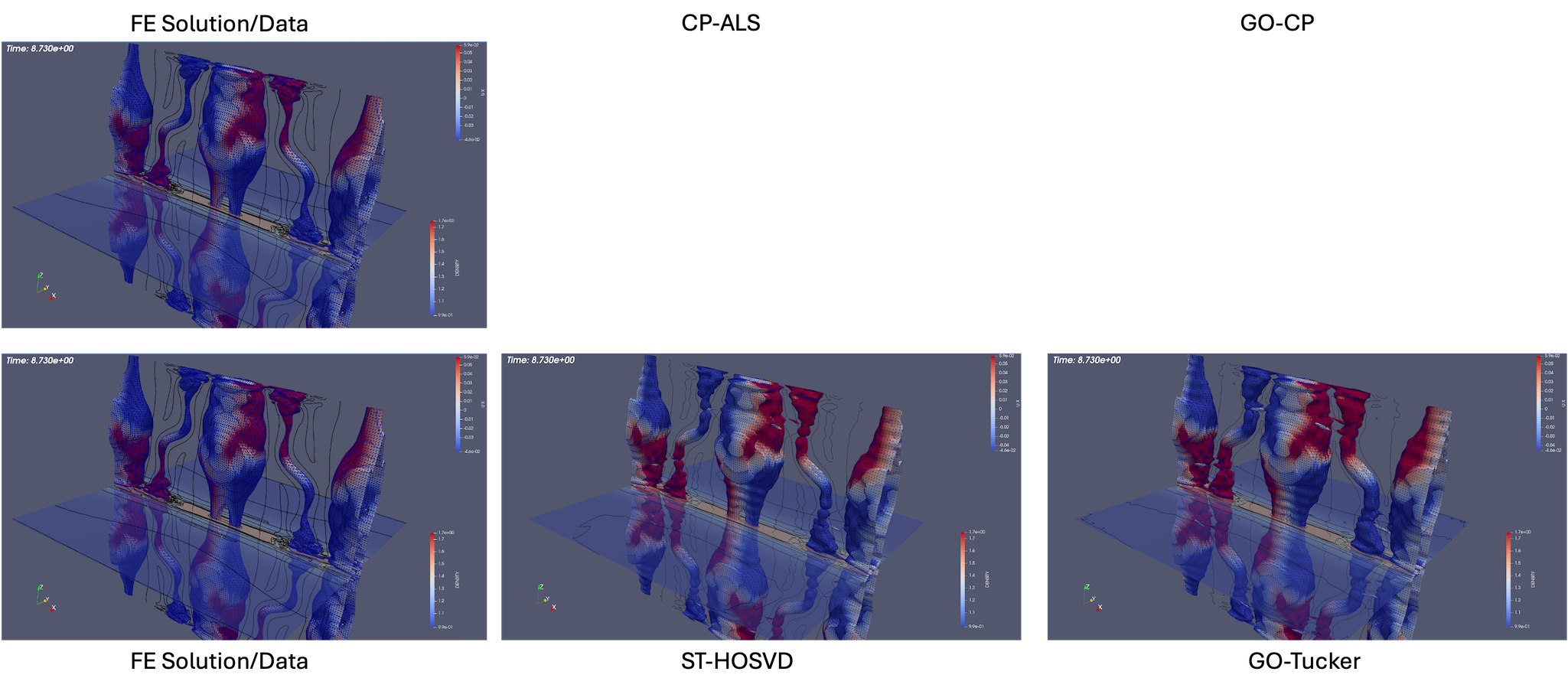}
		\caption*{100x}
			\label{fig:3D_IC_EffectCoarseDatabase_Reconstruction_100x_goal_oriented_t_8.730}
	\end{minipage}
	\begin{minipage}[b]{0.6\textwidth}
	\centering
		\includegraphics[width=\textwidth]{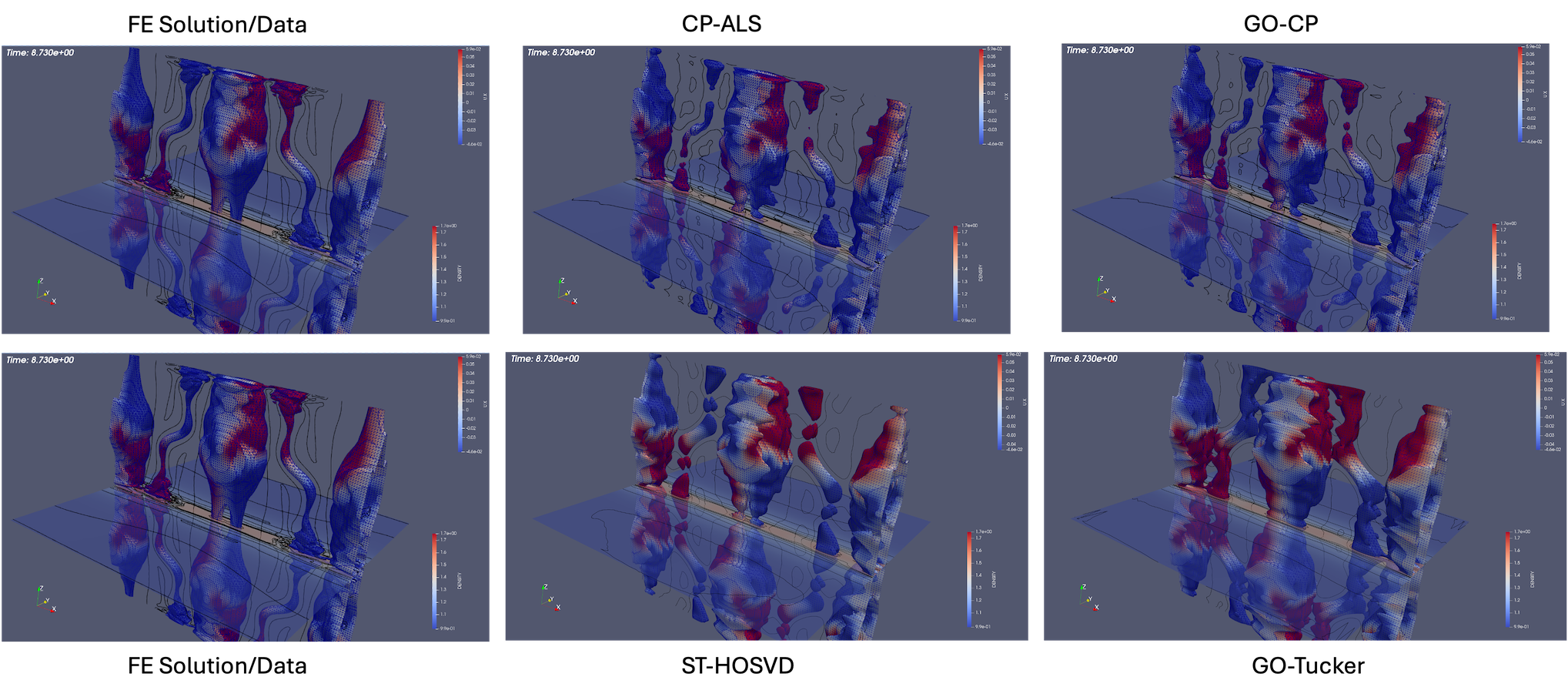}
		\caption*{1,000x}
			\label{fig:3D_IC_EffectCoarseDatabase_Reconstruction_1000x_goal_oriented_t_8.730}
	\end{minipage}
	\hfill
	\begin{minipage}[b]{0.6\textwidth}
	\centering
		\includegraphics[width=\textwidth]{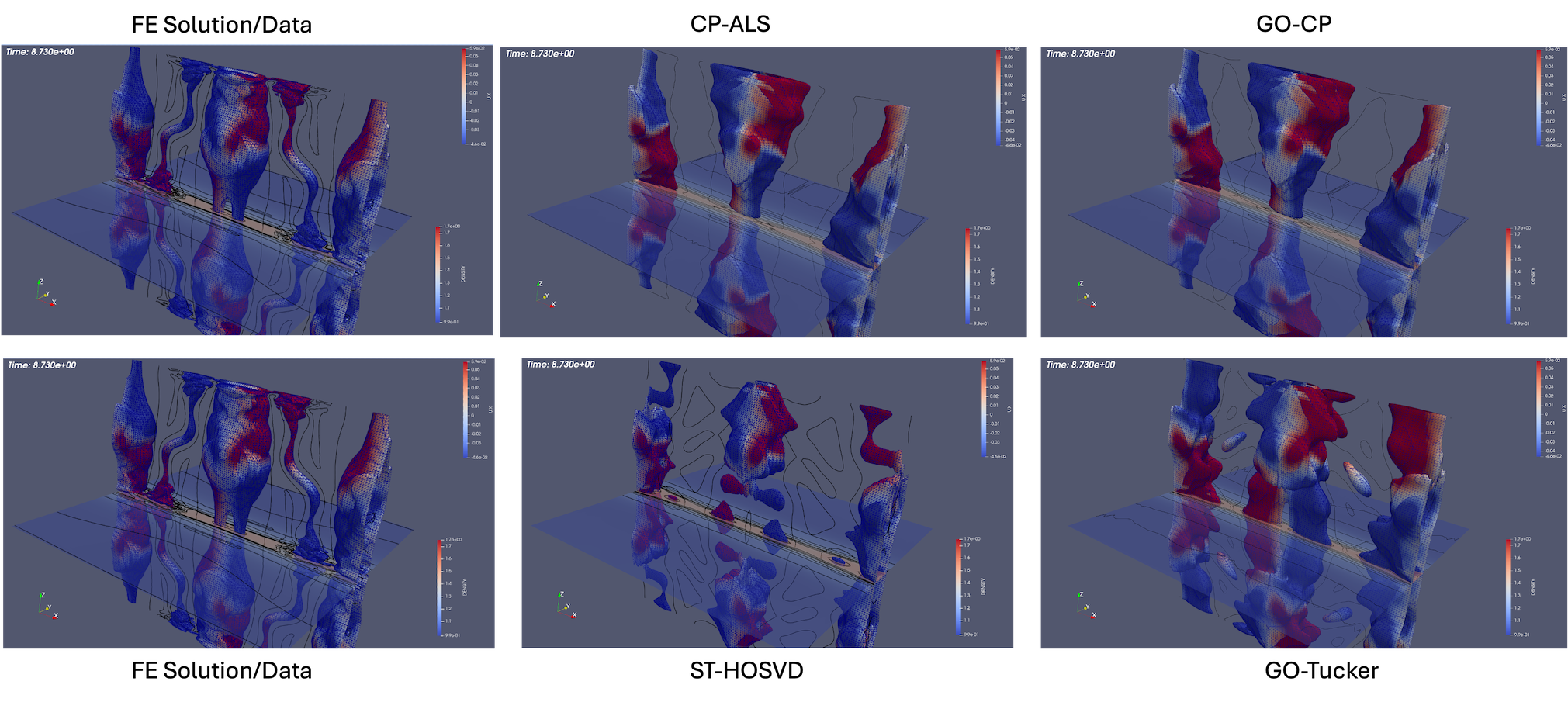}
		\caption*{10,000x}
			\label{fig:3D_IC_EffectCoarseDatabase_Reconstruction_10000x_goal_oriented_t_8.730}
	\end{minipage}
	\vskip -0.2in
	\caption{Comparison of FE method to the standard and goal-oriented tensor decompositions at various compressions for reconstructing the density iso-surface at 1.5 colored by $\bm{u}_x$, and slice plane colored by density at t = 8.730.}
		\label{fig:3D_IC_EffectCoarseDatabase_Reconstructions_goal_oriented_t_8.730}
\end{figure*}

\begin{figure*}[h!]
	\centering
	\begin{minipage}[b]{0.6\textwidth}
	\centering
		\includegraphics[width=\textwidth]{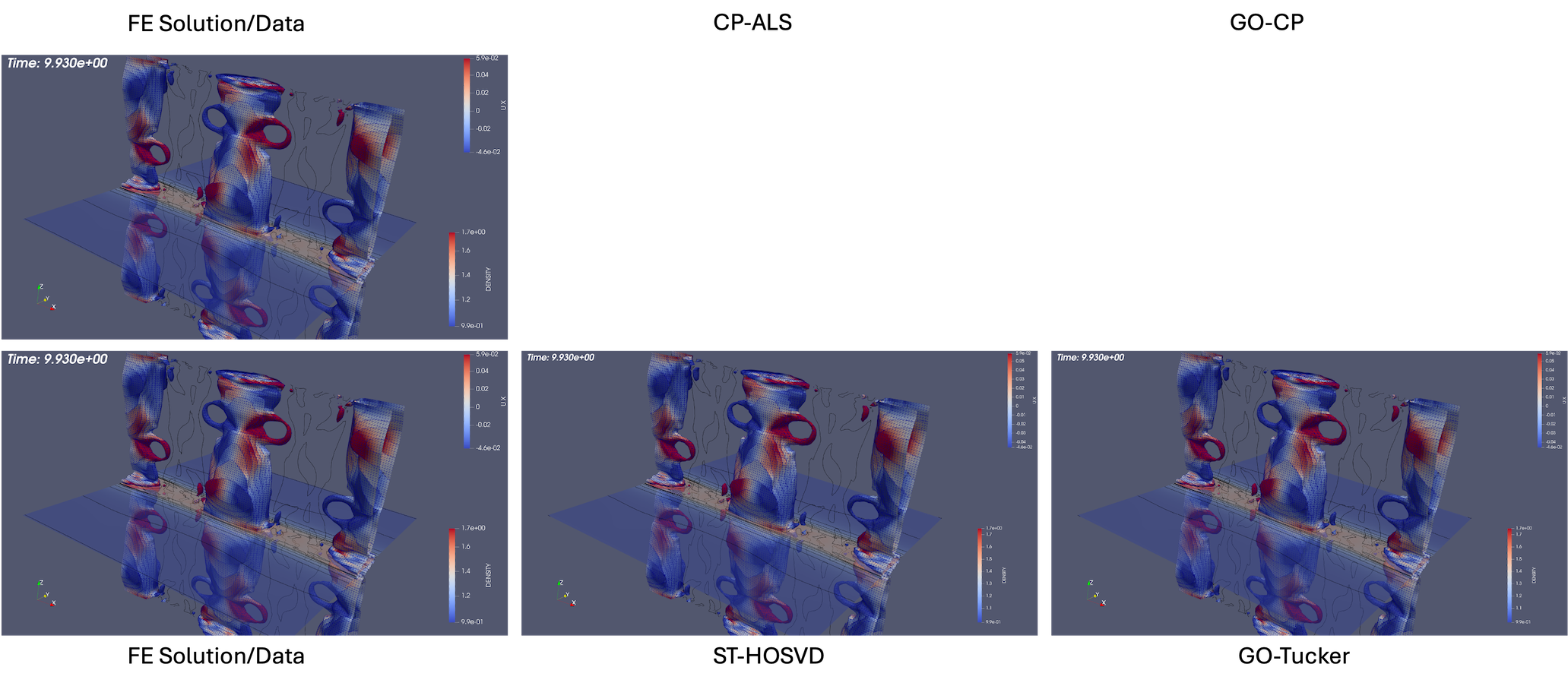}
		\caption*{10x}
			\label{fig:3D_IC_EffectCoarseDatabase_Reconstruction_10x_goal_oriented_t_9.930}
	\end{minipage}
	\hfill
	\begin{minipage}[b]{0.6\textwidth}
	\centering
		\includegraphics[width=\textwidth]{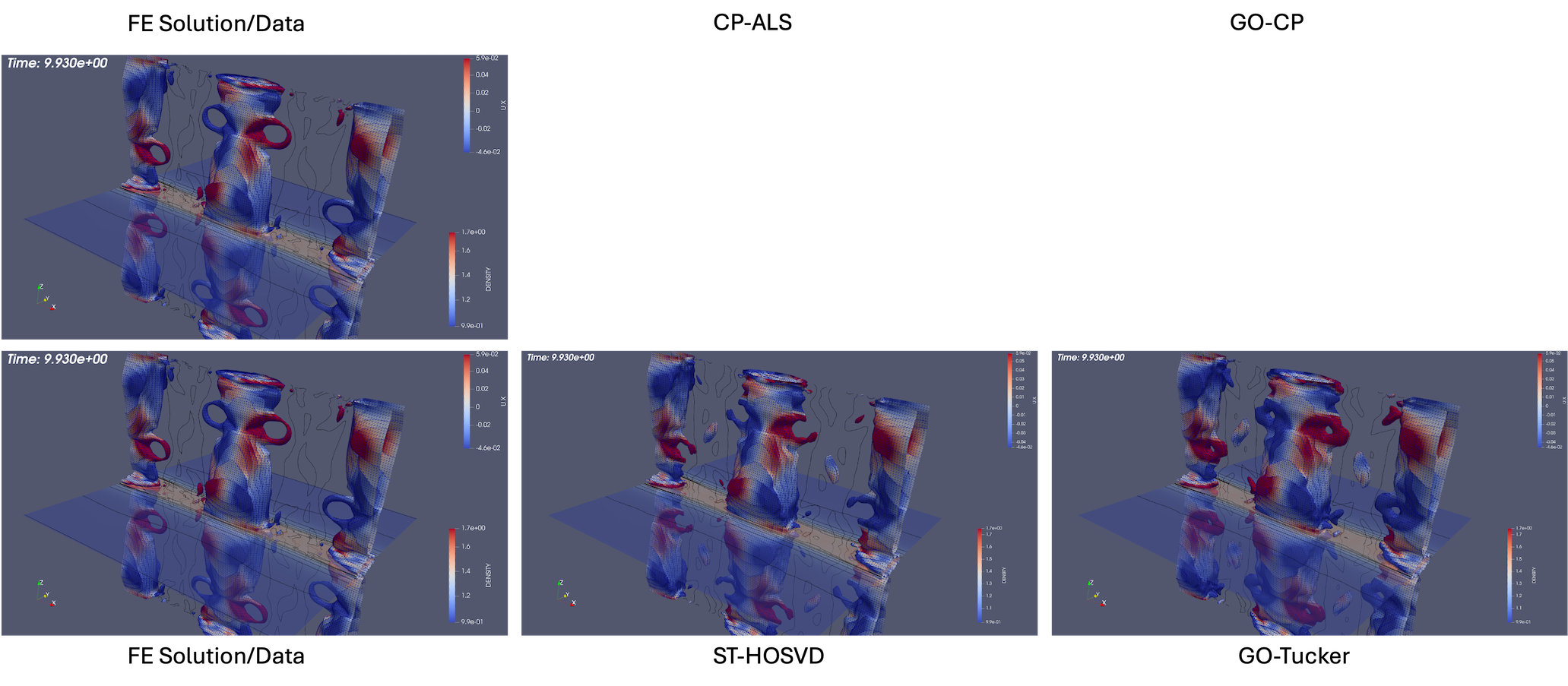}
		\caption*{100x}
			\label{fig:3D_IC_EffectCoarseDatabase_Reconstruction_100x_goal_oriented_t_9.930}
	\end{minipage}
	\begin{minipage}[b]{0.6\textwidth}
	\centering
		\includegraphics[width=\textwidth]{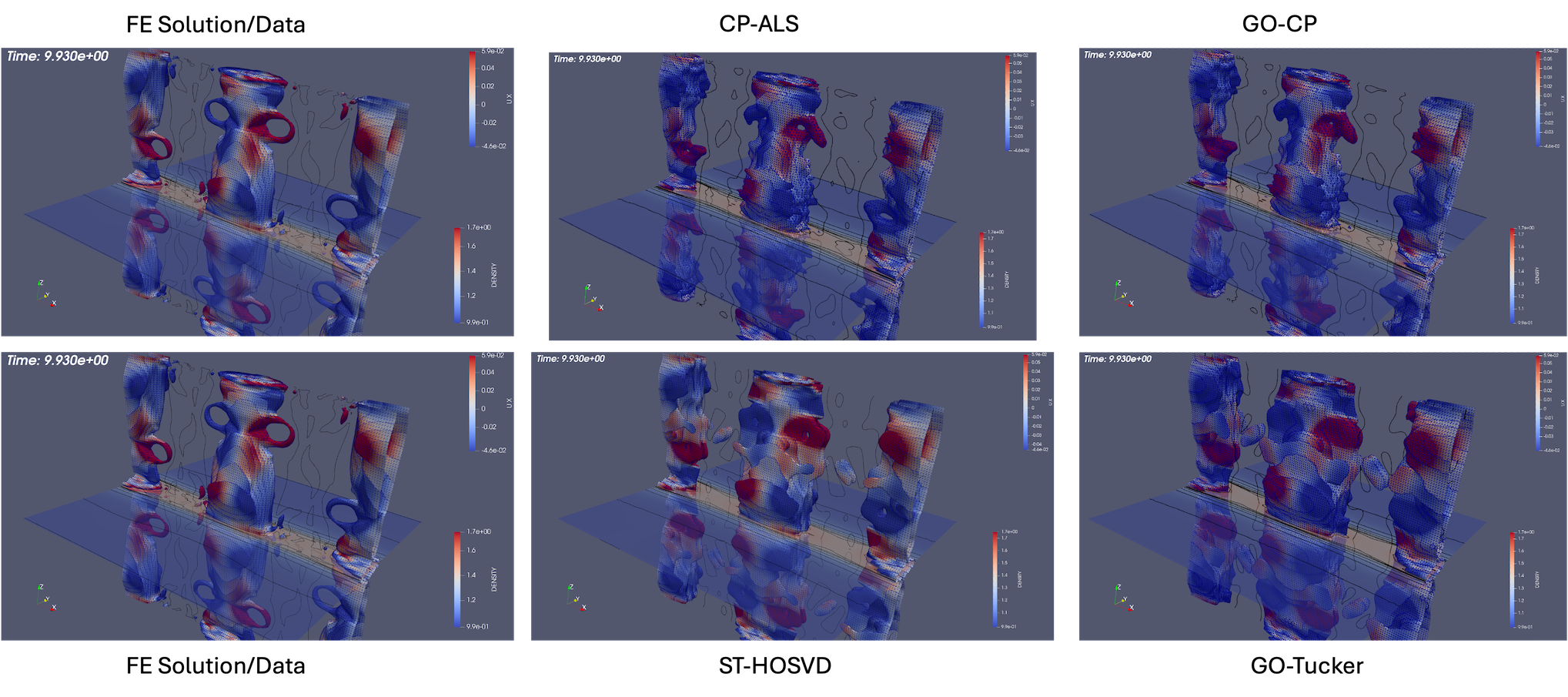}
		\caption*{1,000x}
			\label{fig:3D_IC_EffectCoarseDatabase_Reconstruction_1000x_goal_oriented_t_9.930}
	\end{minipage}
	\hfill
	\begin{minipage}[b]{0.6\textwidth}
	\centering
		\includegraphics[width=\textwidth]{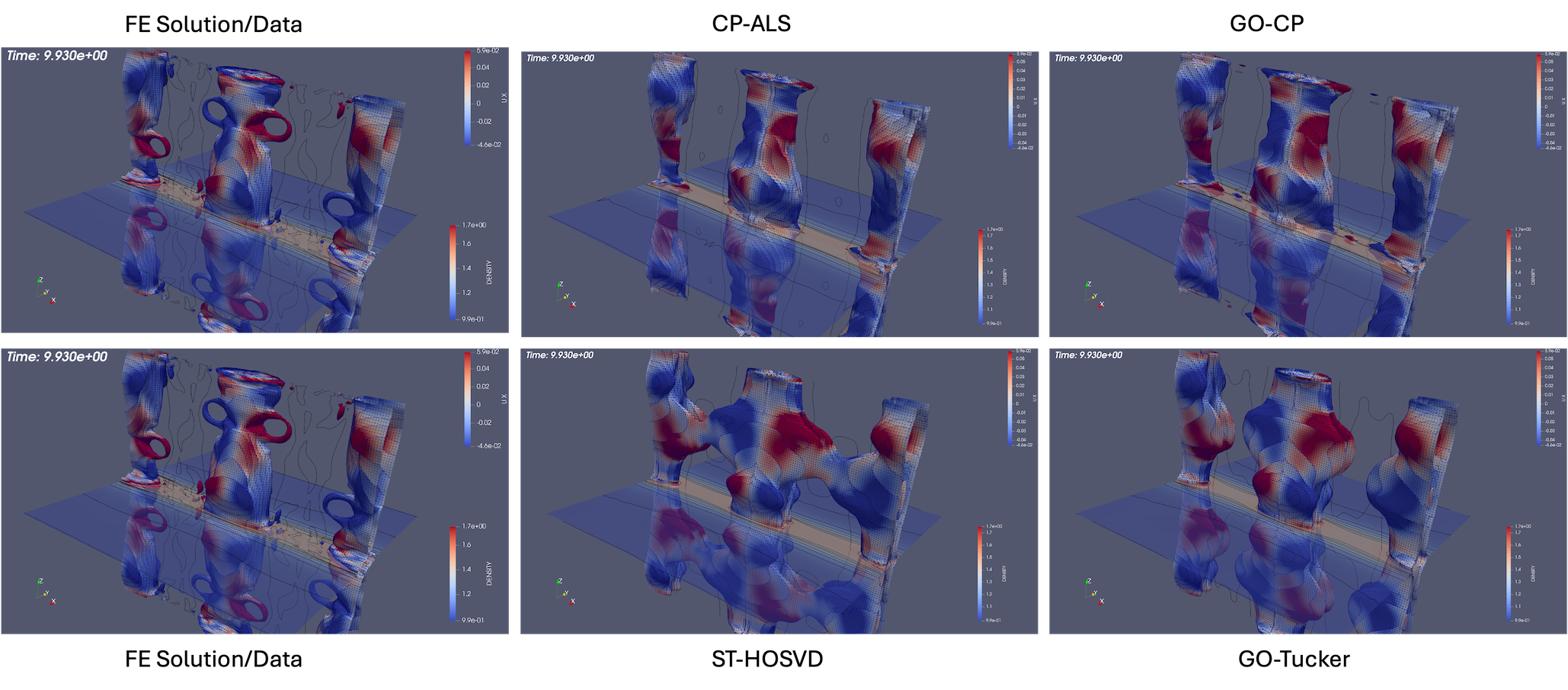}
		\caption*{10,000x}
			\label{fig:3D_IC_EffectCoarseDatabase_Reconstruction_10000x_goal_oriented_t_9.930}
	\end{minipage}
	\vskip -0.2in
	\caption{Comparison of FE method to the standard and goal-oriented tensor decompositions at various compressions for reconstructing the density iso-surface at 1.5 colored by $\bm{u}_x$, and slice plane colored by density at t = 9.930.}
		\label{fig:3D_IC_EffectCoarseDatabase_Reconstructions_goal_oriented_t_9.930}
\end{figure*}

\Cref{fig:island_coalescence_qois1} and  \cref{fig:island_coalescence_qois2} present the improvement in the satisfactions of the QoI evolution with the goal-oriented approaches, where the ranks/tolerances were chosen to yield compression ratios of approximately 10x, 100x, 1,000x, and 10,000x.  Here we see significant error in the QoIs for the 1,000x and 10,000x compressions for the traditional approach, which are substantially reduced for the goal-oriented approach.  As with the tearing mode data, essentially no change in the relative tensor reconstruction errors are observed.
\clearpage
As described above, the goal-oriented methods do strongly improve the QoI reconstruction.  In general these methods also appear to improve the iso-surface reconstructions in most cases but not uniformly. To illustrate this effect, \cref{fig:3D_IC_EffectCoarseDatabase_Reconstructions_goal_oriented_t_8.730} and \cref{fig:3D_IC_EffectCoarseDatabase_Reconstructions_goal_oriented_t_9.930} present results at times $t = 8.730$ and $t = 9.930$ respectively for reduction levels of 10x, 100x, 1,000x, and 10,000x. The results  look very good even on this coarse-database for the 10x, 100x compression ratios for the goal-oriented approaches. For higher compression ratios, it is evident for the $t = 8.730$ data that the goal-oriented methods appear to slightly improve the iso-surface reconstructions. This can most clearly be seen for the goal-oriented methods reconstructing some of the slender iso-surfaces at 1000x and for Tucker  at 10,000x reforming the continuous surfaces of the central iso-surface that has two islands that exist on the intersection with the ($x,y)$ plane, whereas the ST-HOSVD and CP-ALS-based  decompositions only show a single island outline. 

For the latter time of $t = 9.930$ the data is mixed. In~\cref{fig:3D_IC_EffectCoarseDatabase_Reconstructions_goal_oriented_t_9.930}, the 10x reductions look excellent, for 100x a slight improvement of the surface seems to appear with the loop structure reappearing in the goal-oriented Tucker reconstruction. Then for 1,000x a slight degradation in the loop structure seems evident in that the openings become closed surfaces compared to the standard ST-HOSVD. However, with a further increase in the compression rate to 10,000x the goal-oriented approach does isolate the three major iso-surfaces showing an improvement from the standard ST-HOSVD that merged these surfaces. 

\clearpage

%% file: sections/conclusion.tex

In this paper, we introduced a new approach for low-rank modeling of high-dimensional numerical simulation data based on tensor decompositions that augment the traditional optimization problems used to compute these tensor decompositions with additional functionals in the objective function that attempt to preserve application-specific quantities of interest and invariants.  The method was introduced specifically for CP and Tucker decompositions, but the approach is general and could likely be applied to any tensor format that can be computed through an optimization-based method.  Examples of applying the technique to large data sets arising from the numerical simulation of turbulent combustion and plasma physics relevant to fusion energy were presented, demonstrating several orders of magnitude reduction in the error for these QoIs while maintain the same level of accuracy in the overall tensor reconstruction.  When incorporating multiple QoIs in the problem formulation, we found trust region-Newton optimization methods based on Gauss-Newton Hessian approximations effective in solving the resulting optimization problems where other methods such as L-BFGS struggled.  

While our results compared the effectiveness of CP and Tucker decompositions for reduction of scientific data, our aim was not to claim that one these decompositions to be superior over the other.  Generally, we found the level of error for a given compression ratio to be comparable between the two methods, with comparable levels of reduction in QoI error for the goal-oriented approach.  We often found the Tucker method easier to deal with in practical settings however, particularly for high-accuracy decompositions.  With Tucker, setting a tight ST-HOSVD solver tolerance is straightforward and not terribly expensive computationally.  However, obtaining high-accuracy CP decompositions is challenging numerically due to the very high ranks (i.e., $O$(10,000)) required which existing CP decomposition algorithms (e.g., CP-ALS) and computational kernels supporting them (e.g., MTTKRP) are often not optimized to handle effectively.  Conversely, we found the optimization problem defining goal-oriented CP decompositions typically easier to solve numerically.  In fact, the challenges with solving the optimization problem for Tucker decompositions motivated the trust-region Newton optimization approach studied here.

While the approach was demonstrated to be successful in achieving our goal of reducing QoI error in reconstructions arising from the low-rank models, several challenges remain that could be studied through future work.  In particular, we found the goal-oriented approach to be substantially more expensive than the traditional solution methods the models are based on (i.e., CP-ALS and ST-HOSVD).  Part of this is due to our serial Matlab-based implementation where vectorizing the evaluation of the QoI functions and their derivatives is challenging.  We expect the approach to be much more competitive in terms of computational cost with a high-performance, parallel implementation.  Other potential approaches for reducing cost include randomized algorithms, incremental update/streaming algorithms, and improved preconditioning that incorporates the QoI terms in some fashion.  Effective stopping criteria beyond a fixed number of iterations were also not addressed in this work, and we suspect that extending the formulation to one on Riemannian manifolds to eliminate degeneracies  arising from non-uniqueness where traditional stopping criteria based on gradient norms could be used would likely be fruitful.  Finally, a theoretical study of how the goal-oriented approach can decrease QoI error without significantly increasing tensor reconstruction error, illustrating similarities and differences of the traditional and goal-oriented approaches, is needed.

%% file: sections/acknowledgement.tex
This work was supported by the Laboratory Directed Research and Development program at Sandia
National Laboratories, a multimission laboratory managed and operated by National Technology and
Engineering Solutions of Sandia LLC, a wholly owned subsidiary of Honeywell International Inc. for the U.S.
Department of Energy’s National Nuclear Security Administration under contract DE-NA0003525.

%% file: sections/go_grad_hess_prec.tex


In this section, we derive formulas for the efficient evaluation of the gradient and Hessian of the objective function $f_{go}$ needed for the optimization methods described above.  Calculation of these quantities for the tensor term (i.e., $f(\X,\M)$) is standard within the literature, so we focus solely on the goal-oriented terms (formulas for the tensor term are provided in \cref{sec:fro_grad_hess_prec} for reference).  For simplicity, we consider only a single goal term since the extension to multiple goals is straightforward, and therefore drop the $q$ subscript.  We also ignore variable scaling as updating the goal derivatives to include it is straightforward, and assume $\mathbb{T}=\{1,\dots,\tau\}$.  Define $G = \sum_{t=1}^\tau \left(g(\X_t)-g(\M_t)\right)^2$, which has the form of a nonlinear least-squares objective.  Given a CP or Tucker model $\M$, let $\Vc{v}$ be the vector containing the factor matrix coefficients (and core coefficients for Tucker). Define the residual function
\begin{equation}
  \Vc{F}(\Vc{v}) =
    \begin{bmatrix}
      g(\X_1)-g(\M_1) \\
      \vdots \\
      g(\X_\tau)-g(\M_\tau)
    \end{bmatrix}
 \end{equation}
so that $G = \Vc{F}^T \Vc{F}$.   Then $\Vc{\nabla}_{\Vc{v}} G = 2\Mx{J}^T \Vc{F}$ where $\Mx{J}$ is the Jacobian matrix $\Mx{D}_{\Vc{v}} \Vc{F}$.  Furthermore, the Gauss-Newton Hessian approximation $\Mx{H}\approx\Vc{\nabla}_{\Vc{v}}^2 G$ is given by $\Mx{H} = 2\Mx{J}^T \Mx{J}$.  For simplicity in what follows, we assume $\X\in\Real^{I_1 \times I_2 \times \tau}$ only contains a single spatial mode and therefore is a 3-way tensor.  The generalization to higher-dimensional $\X$ is straightforward.  Given a corresponding CP/Tucker model $\M_t$ for each $t$, define $\Tn{Z}\in\Real^{I_1 \times I_2 \times \tau}$ such that $z_{{i_1}{i_2}t} = \partial g/\partial m_{{i_1}{i_2}t}$ is the tensor of partial derivatives with respect to the reconstructed tensor model entries.
 
\leveldown{Goal-oriented CP}
\label{sec:go_cp}
 
Let $\M=\llbracket \Mx{A}, \Mx{B}, \Mx{C} \rrbracket$ be a given 3-way, rank-$R$ (unweighted) CP model.  Then $\Vc{v} = [ \vc(\Mx{A})^T \vc(\Mx{B})^T \vc(\Mx{C})^T ]^T$.  Note that $\M_t = \llbracket \Mx{A}, \Mx{B}, \Vc{c}_t \rrbracket$ where $\Vc{c}_t = \Mx{C}(t,:)$.  Let $\Vc{v}_t = \vc(\M_t)$.  Then we have
\begin{equation}
 \begin{aligned}
   \Vc{v}_t &= \vc(\M_t) \\
   &= \Mx{P}_1^T\vc((\M_t)_{(1)}) \\
   &= \Mx{P}_1^T\vc(\Mx{A}(\Vc{c}_t\odot\Mx{B})^T) \\
   &= \Mx{P}_1^T((\Vc{c}_t\odot\Mx{B})\otimes\Mx{I}_{I_1})\vc(\Mx{A})
 \end{aligned}
\end{equation}
where $\Mx{P}_k$ is the mode-$k$ perfect shuffle matrix such that $\Mx{P}_k\vc(\X) = \vc(\Mx{X}_{(k)})$ and $\Mx{I}_{I_1}$ is the $I_1\times I_1$ identity matrix.  Thus
\begin{equation}
   \frac{\partial\Vc{v}_t}{\partial\vc(\Mx{A})} = \Mx{P}_1^T((\Vc{c}_t\odot\Mx{B})\otimes\Mx{I}_{I_1})
\end{equation}
and through similar manipulations we obtain
\begin{equation}
 \begin{aligned}
   \frac{\partial\Vc{v}_t}{\partial\vc(\Mx{B})} &= \Mx{P}_2^T((\Vc{c}_t\odot\Mx{A})\otimes\Mx{I}_{I_2}), \\
   \frac{\partial\Vc{v}_t}{\partial\vc(\Vc{c}_t)} &= \Mx{P}_3^T((\Mx{B}\odot\Mx{A})\otimes\Mx{I}_{\tau}). \\
 \end{aligned}
\end{equation}
The Jacobian matrix $\Mx{J}=\Mx{D}_{\Vc{v}} \Vc{F}$ has the block structure 
\begin{equation}
   \Mx{J} = 
   \begin{bmatrix}
    \Mx{J}_{\Mx{A}}^1 & \Mx{J}_{\Mx{B}}^1 & \Mx{J}_{\Vc{c}_1}^1 & 0 & \hdots & 0 \\
    \Mx{J}_{\Mx{A}}^2 & \Mx{J}_{\Mx{B}}^2 & 0 & \Mx{J}_{\Vc{c}_2}^2 & \hdots & 0 \\
    \vdots & \vdots & \vdots & \vdots & \ddots & \vdots \\
    \Mx{J}_{\Mx{A}}^\tau & \Mx{J}_{\Mx{B}}^\tau & 0 & 0 & \hdots & \Mx{J}_{\Vc{c}_\tau}^\tau  \\
  \end{bmatrix}
\end{equation}
where
\begin{equation}
 \begin{aligned}
   \Mx{J}_{\Mx{A}}^t &= \frac{\partial f_t}{\partial\Vc{v}_t}\frac{\partial\Vc{v}_t}{\partial\vc(\Mx{A})} \\
   &= -\vc(\Tn{Z}^t)^T \frac{\partial\Vc{v}_t}{\partial\vc(\Mx{A})} \\
   &= -\vc(\Mx{Z}_{(1)}^t)^T \Mx{P}_1 \Mx{P}_1^T((\Vc{c}_t\odot\Mx{B})\otimes\Mx{I}_{I_1}) \\
   &= -\left(((\Vc{c}_t\odot\Mx{B})^T\otimes\Mx{I}_{I_1})\vc(\Mx{Z}_{(1)}^t) \right)^T \\
   &= -\vc(\Mx{Z}_{(1)}(\Vc{c}_t\odot\Mx{B})^T)^T
 \end{aligned}
\end{equation}
for $t\in\mathbb{T}$.  Similarly,
\begin{equation}
 \begin{aligned}
   \Mx{J}_{\Mx{B}}^t &= -\vc(\Mx{Z}_{(2)}(\Vc{c}_t\odot\Mx{A})^T)^T, \\
   \Mx{J}_{\Vc{c}_t}^t &= -\vc(\Mx{Z}_{(3)}(\Mx{B}\odot\Mx{A})^T)^T.
 \end{aligned}
\end{equation}
Thus the Jacobian blocks are given by the well-known Matricized Tensor Times Khatri-Rao Product (MTTKRP) operation with the gradient tensors $\Tn{Z}^t$.  Furthermore, one can efficiently compute Jacobian-vector and Jacobian-transpose-vector products using these formulas.  Given $\Vc{u}\in\Real^T$ one can easily show
\begin{equation}\label{eq:gocp_jac_trans_vec}
\Mx{J}^T\Vc{u} = 
\begin{bmatrix}
  \sum_{t=1}^\tau (\Mx{J}_{\Mx{A}}^t)^T u_t \\
  \sum_{t=1}^\tau (\Mx{J}_{\Mx{B}}^t)^T u_t \\
  (\Mx{J}_{\Vc{c}_1}^1)^T u_1 \\
  \vdots \\
  (\Mx{J}_{\Vc{c}_\tau}^\tau)^T u_\tau \\
\end{bmatrix}
=
\begin{bmatrix}
  \vc(\Mx{\tilde{Z}}_{(1)}(\Mx{C}\odot\Mx{B})) \\
  \vc(\Mx{\tilde{Z}}_{(2)}(\Mx{C}\odot\Mx{A})) \\
  \vc(\Mx{\tilde{Z}}_{(3)}(\Mx{B}\odot\Mx{A})) \\
\end{bmatrix}
\end{equation}
where $\Tn{\tilde{Z}}(:,:,t) = -u_t \Tn{Z}(:,:,t)$, which allows for efficient calculation of the gradient.  Moreover, given $\Vc{w} = [\Vc{w}_{\Mx{A}}^T \; \Vc{w}_{\Mx{B}}^T \; \Vc{w}_{\Mx{C}}^T ]^T = [\vc(\Mx{W}_{\Mx{A}})^T \vc(\Mx{W}_{\Mx{B}})^T \vc(\Mx{W}_{\Mx{C}})^T ]^T \in\Real^{R (I_1+I_2+\tau)}$,
\begin{equation}
 \Mx{J}\Vc{w} = 
\begin{bmatrix}
  \Mx{J}_{\Mx{A}}^1 \Vc{w}_{\Mx{A}} + \Mx{J}_{\Mx{B}}^1 \Vc{w}_{\Mx{A}} + \Mx{J}_{\Vc{c}_1}^1 \Vc{w}_{\Vc{c}_1} \\
  \vdots \\
  \Mx{J}_{\Mx{A}}^\tau \Vc{w}_{\Mx{A}} + \Mx{J}_{\Mx{B}}^\tau \Vc{w}_{\Mx{A}} + \Mx{J}_{\Vc{c}_\tau}^\tau \Vc{w}_{\Vc{c}_\tau} \\
\end{bmatrix}
\end{equation}
where, e.g.,
\begin{equation}
\begin{aligned}
  \Mx{J}_{\Mx{A}}^t \Vc{w}_{\Mx{A}} &= -\vc(\Mx{Z}_{(1)}^t)^T \Mx{P}_1 \Mx{P}_1^T((\Vc{c}_t\odot\Mx{B})\otimes\Mx{I}_{I_1})\vc(\Mx{W}_{\Mx{A}}) \\
  &= -\vc(\Mx{Z}_{(1)}^t)^T \Mx{P}_1 \Mx{P}_1^T\vc(\Mx{W}_{\Mx{A}}(\Vc{c}_t\odot\Mx{B})^T) \\
  &= -\vc(\Mx{Z}_{(1)}^t)^T \Mx{P}_1 \Mx{P}_1^T\vc(\llbracket \Mx{W}_{\Mx{A}}, \Mx{B}, \Vc{c}_t \rrbracket_{(1)}) \\
  &= -\vc(\Tn{Z}^t)^T \vc(\llbracket \Mx{W}_{\Mx{A}}, \Mx{B}, \Vc{c}_t \rrbracket). \\
\end{aligned}
\end{equation}
Thus
\begin{multline}\label{eq:gocp_jac_vec}
  \Mx{J}_{\Mx{A}}^t \Vc{w}_{\Mx{A}} + \Mx{J}_{\Mx{B}}^t \Vc{w}_{\Mx{A}} + \Mx{J}_{\Vc{c}_t}^t \Vc{w}_{\Vc{c}_t} = \\
  -\vc(\Tn{Z}^t)^T\vc\left( \llbracket \Mx{W}_{\Mx{A}}, \Mx{B}, \Vc{c}_t \rrbracket + \llbracket \Mx{A}, \Mx{W}_{\Mx{B}}, \Vc{c}_t \rrbracket + \llbracket \Mx{A}, \Mx{B}, \Vc{w}_{\Vc{c}_t} \rrbracket \right).
\end{multline}
Combining this with \cref{eq:gocp_jac_trans_vec}, one can efficiently compute Gauss-Newton Hessian-vector products.

\levelstay{Goal-oriented Tucker}
 \label{sec:go_tucker}
 
Let $\M=\llbracket \Tn{G}; \Mx{A}, \Mx{B}, \Mx{C} \rrbracket$ be a given 3-way, rank-$(R_1,R_2,R_3)$ Tucker model.  Then $\Vc{v} = [ \vc(\Tn{G})^T \vc(\Mx{A})^T \vc(\Mx{B})^T \vc(\Mx{C})^T ]^T$ and $\M_t = \llbracket \Tn{G}; \Mx{A}, \Mx{B}, \Vc{c}_t \rrbracket$ where $\Vc{c}_t = \Mx{C}(t,:)$.  With $\Vc{v}_t = \vc(\M_t)$ and using similar techniques as above, one can show the Jacobian matrix has a similar block structure
\begin{equation}
   \Mx{J} = 
   \begin{bmatrix}
    \Mx{J}_{\Tn{G}}^1 & \Mx{J}_{\Mx{A}}^1 & \Mx{J}_{\Mx{B}}^1 & \Mx{J}_{\Vc{c}_1}^1 & 0 & \hdots & 0 \\
    \Mx{J}_{\Tn{G}}^2 & \Mx{J}_{\Mx{A}}^2 & \Mx{J}_{\Mx{B}}^2 & 0 & \Mx{J}_{\Vc{c}_2}^2 & \hdots & 0 \\
    \vdots & \vdots & \vdots & \vdots & \vdots & \ddots & \vdots \\
    \Mx{J}_{\Tn{G}}^\tau & \Mx{J}_{\Mx{A}}^\tau & \Mx{J}_{\Mx{B}}^\tau & 0 & 0 & \hdots & \Mx{J}_{\Vc{c}_\tau}^\tau  \\
  \end{bmatrix}
\end{equation}
where
\begin{equation}
  \begin{aligned}
    \Mx{J}_{\Tn{G}}^t &= -\vc(\Tn{Z}^t \times_1 \Mx{A}^T \times_2 \Mx{B}^T \times_3 \Vc{c}_t^T)^T, \\
    \Mx{J}_{\Mx{A}}^t &= -\vc((\Tn{Z}^t \times_2 \Mx{B}^T \times_3 \Vc{c}_t^T)_{(1)} \Mx{G}_{(1)}^T)^T, \\
    \Mx{J}_{\Mx{B}}^t &= -\vc((\Tn{Z}^t \times_1 \Mx{A}^T \times_3 \Vc{c}_t^T)_{(2)} \Mx{G}_{(2)}^T)^T, \\
    \Mx{J}_{\Vc{c}_t}^t &= -\vc((\Tn{Z}^t \times_1 \Mx{A}^T \times_2 \Mx{B}^T)_{(3)} \Mx{G}_{(3)}^T)^T.
  \end{aligned}
\end{equation}
Thus, as before,
\begin{equation}\label{eq:gotucker_jac_trans_vec}
\Mx{J}^T\Vc{u} = 
\begin{bmatrix}
     \vc(\Tn{\tilde{Z}} \times_1 \Mx{A}^T \times_2 \Mx{B}^T \times_3 \Mx{C}^T) \\
     \vc((\Tn{\tilde{Z}} \times_2 \Mx{B}^T \times_3 \Mx{C}^T)_{(1)} \Mx{G}_{(1)}^T) \\
     \vc((\Tn{\tilde{Z}} \times_1 \Mx{A}^T \times_3 \Mx{C}^T)_{(2)} \Mx{G}_{(2)}^T) \\
     \vc((\Tn{\tilde{Z}} \times_1 \Mx{A}^T \times_2 \Mx{B}^T)_{(3)} \Mx{G}_{(3)}^T)
\end{bmatrix}
\end{equation}
with $\Tn{\tilde{Z}}$ defined as above.  Similarly,
\begin{equation}
 \Mx{J}\Vc{w} = 
\begin{bmatrix}
  \Mx{J}_{\Tn{G}}^1 \Vc{w}_{\Tn{G}} + \Mx{J}_{\Mx{A}}^1 \Vc{w}_{\Mx{A}} + \Mx{J}_{\Mx{B}}^1 \Vc{w}_{\Mx{A}} + \Mx{J}_{\Vc{c}_1}^1 \Vc{w}_{\Vc{c}_1} \\
  \vdots \\
  \Mx{J}_{\Tn{G}}^\tau \Vc{w}_{\Tn{G}} + \Mx{J}_{\Mx{A}}^\tau \Vc{w}_{\Mx{A}} + \Mx{J}_{\Mx{B}}^\tau \Vc{w}_{\Mx{A}} + \Mx{J}_{\Vc{c}_\tau}^\tau \Vc{w}_{\Vc{c}_\tau} \\
\end{bmatrix}
\end{equation}
with 
\begin{multline} \label{eq:gotucker_jac_vec}
  \Mx{J}_{\Tn{G}}^t \Vc{w}_{\Tn{G}} + \Mx{J}_{\Mx{A}}^t \Vc{w}_{\Mx{A}} + \Mx{J}_{\Mx{B}}^t \Vc{w}_{\Mx{A}} + \Mx{J}_{\Vc{c}_t}^t \Vc{w}_{\Vc{c}_t} = 
  -\vc(\Tn{Z}^t)^T\vc(\Tn{W}_{\Tn{G}}\times_1\Mx{A}\times_2\Mx{B}\times_3\Vc{c}_t \\
  + \Tn{G}\times_1\Mx{W}_{\Mx{A}}\times_2\Mx{B}\times_3\Vc{c}_t
  + \Tn{G}\times_1\Mx{A}\times_2\Mx{W}_{\Mx{B}}\times_3\Vc{c}_t 
  + \Tn{G}\times_1\Mx{A}\times_2\Mx{B}\times_3\Vc{w}_{\Vc{c}_t}).
\end{multline}

\levelstay{Frobenius Gradient, Hessian, and Preconditioner}
\label{sec:fro_grad_hess_prec}

While the goal-oriented formulation in \cref{eq:goal_oriented} is general and can be applied with any statistical loss function $f(\X,\M)$, our examples below focus specifically on Frobenius (sum-of-squares) loss with $f(\X,\M) = \frac{1}{2}\|\X-\M\|_F^2$.  While construction of the gradient, Gauss-Newton Hessian, and diagonal blocks for preconditioning is standard within the literature, we summarize these formulas here for completeness.  As above, we restrict to the 3-way case for simplicity of exposition as these formulas naturally generalize to high-order tensors.  First, in the CP case we have
\begin{equation} \label{eq:cp_grad}
\begin{aligned}
  \frac{\partial f}{\partial\Mx{A}} &= \Mx{A}(\Mx{C}^T\Mx{C}\ast\Mx{B}^T\Mx{B}) - \Mx{X}_{(1)}(\Mx{C}\odot\Mx{B}), \\
  \frac{\partial f}{\partial\Mx{B}} &= \Mx{B}(\Mx{C}^T\Mx{C}\ast\Mx{A}^T\Mx{A}) - \Mx{X}_{(2)}(\Mx{C}\odot\Mx{A}), \\
  \frac{\partial f}{\partial\Mx{C}} &= \Mx{C}(\Mx{B}^T\Mx{B}\ast\Mx{A}^T\Mx{A}) - \Mx{X}_{(3)}(\Mx{B}\odot\Mx{A}), \\
\end{aligned}
\end{equation}
where $\ast$ is the matrix Hadamard (element-wise) product.  Given $\Vc{w}$ as defined above, the corresponding Gauss-Newton Hessian-vector product is
\begin{equation}\label{eq:cp_hess_vec}
  \Mx{J}^T\Mx{J}\Vc{w} = 
  \begin{bmatrix}
    \vc(\Mx{W}_{\Mx{A}}(\Mx{B}^T\Mx{B}\ast\Mx{C}^T\Mx{C}) + \Mx{A}(\Mx{W}_{\Mx{B}}^T\Mx{B}\ast\Mx{C}^T\Mx{C}) + \Mx{A}(\Mx{B}^T\Mx{B}\ast\Mx{W}_{\Mx{C}}^T\Mx{C})) \\
    \vc(\Mx{B}(\Mx{W}_{\Mx{A}}^T\Mx{A}\ast\Mx{C}^T\Mx{C}) + \Mx{W}_{\Mx{B}}(\Mx{A}^T\Mx{A}\ast\Mx{C}^T\Mx{C}) + \Mx{B}(\Mx{A}^T\Mx{A}\ast\Mx{W}_{\Mx{C}}^T\Mx{C})) \\
    \vc(\Mx{C}(\Mx{W}_{\Mx{A}}^T\Mx{A}\ast\Mx{B}^T\Mx{B}) + \Mx{C}(\Mx{A}^T\Mx{A}\ast\Mx{W}_{\Mx{B}}^T\Mx{B}) + \Mx{W}_{\Mx{C}}(\Mx{A}^T\Mx{A}\ast\Mx{B}^T\Mx{B}))
  \end{bmatrix}.
\end{equation}
Finally, the diagonal blocks of the Gauss-Newton Hessian are
\begin{equation}
  \begin{aligned}
    \Mx{J}_{\Mx{A}}^T  \Mx{J}_{\Mx{A}} &= (\Mx{C}^T\Mx{C}\ast\Mx{B}^T\Mx{B})\otimes\Mx{I}_1, \\
    \Mx{J}_{\Mx{B}}^T  \Mx{J}_{\Mx{B}} &= (\Mx{C}^T\Mx{C}\ast\Mx{A}^T\Mx{A})\otimes\Mx{I}_2, \\
    \Mx{J}_{\Mx{C}}^T  \Mx{J}_{\Mx{C}} &= (\Mx{B}^T\Mx{B}\ast\Mx{A}^T\Mx{A})\otimes\Mx{I}_\tau, \\
  \end{aligned}
\end{equation}
whose inverses can be efficiently applied through Cholesky decompositions of the above Hadamard product matrices.

For Tucker, the Gauss-Newton Jacobian is $\Mx{J} = [\Mx{J}_{\Tn{G}} \; \Mx{J}_{\Mx{A}} \; \Mx{J}_{\Mx{B}} \; \Mx{J}_{\Mx{C}}]$ with
\begin{equation}
\begin{aligned}
  \Mx{J}_{\Tn{G}} &= -\Mx{C}\otimes\Mx{B}\otimes\Mx{A}, \\
  \Mx{J}_{\Mx{A}} &= -\Mx{P}_1^T ((\Mx{C}\otimes\Mx{B})\Mx{G}_{(1)}^T\otimes\Mx{I}_1), \\
  \Mx{J}_{\Mx{B}} &= -\Mx{P}_2^T ((\Mx{C}\otimes\Mx{A})\Mx{G}_{(2)}^T\otimes\Mx{I}_2), \\
  \Mx{J}_{\Mx{C}} &= -\Mx{P}_3^T ((\Mx{B}\otimes\Mx{A})\Mx{G}_{(3)}^T\otimes\Mx{I}_\tau). \\
\end{aligned}
\end{equation}
Thus one can compute Jacobian-transpose products as
\begin{equation}\label{eq:tucker_grad}
\Mx{J}^T\Vc{u} = -
\begin{bmatrix}
     \vc(\Tn{U} \times_1 \Mx{A}^T \times_2 \Mx{B}^T \times_3 \Mx{C}^T) \\
     \vc((\Tn{U} \times_2 \Mx{B}^T \times_3 \Mx{C}^T)_{(1)} \Mx{G}_{(1)}^T) \\
     \vc((\Tn{U} \times_1 \Mx{A}^T \times_3 \Mx{C}^T)_{(2)} \Mx{G}_{(2)}^T) \\
     \vc((\Tn{U} \times_1 \Mx{A}^T \times_2 \Mx{B}^T)_{(3)} \Mx{G}_{(3)}^T)
\end{bmatrix}
\end{equation}
where $\Vc{u} = \vc(\Tn{U})$. In particular, this allows calculation of the gradient where $\Tn{U} = \Tn{X}-\Tn{M}$.  Furthermore, Jacobian-vector products can be computed as
\begin{multline} \label{eq:tucker_jac_vec}
  \Mx{J}_{\Tn{G}} \Vc{w}_{\Tn{G}} + \Mx{J}_{\Mx{A}} \Vc{w}_{\Mx{A}} + \Mx{J}_{\Mx{B}} \Vc{w}_{\Mx{A}} + \Mx{J}_{\Mx{C}} \Vc{w}_{\Mx{C}} = 
  -\vc(\Tn{W}_{\Tn{G}}\times_1\Mx{A}\times_2\Mx{B}\times_3\Mx{C} \\
  + \Tn{G}\times_1\Mx{W}_{\Mx{A}}\times_2\Mx{B}\times_3\Mx{C}
  + \Tn{G}\times_1\Mx{A}\times_2\Mx{W}_{\Mx{B}}\times_3\Mx{C} 
  + \Tn{G}\times_1\Mx{A}\times_2\Mx{B}\times_3\Vc{w}_{\Mx{C}}).
\end{multline}
Finally, the Gauss-Newton Hessian diagonal blocks for preconditioning are
\begin{equation}
  \begin{aligned}
    \Mx{J}_{\Tn{G}}^T  \Mx{J}_{\Tn{G}} &= \Mx{C}^T\Mx{C}\otimes\Mx{B}^T\Mx{B}\otimes\Mx{A}^T\Mx{A}, \\
    \Mx{J}_{\Mx{A}}^T  \Mx{J}_{\Mx{A}} &= \Tn{G}_{(1)}(\Mx{C}^T\Mx{C}\otimes\Mx{B}^T\Mx{B})\Tn{G}_{(1)}^T\otimes\Mx{I}_1, \\
    \Mx{J}_{\Mx{B}}^T  \Mx{J}_{\Mx{B}} &= \Tn{G}_{(2)}(\Mx{C}^T\Mx{C}\otimes\Mx{A}^T\Mx{A})\Tn{G}_{(2)}^T\otimes\Mx{I}_2, \\
    \Mx{J}_{\Mx{C}}^T  \Mx{J}_{\Mx{C}} &= \Tn{G}_{(3)}(\Mx{B}^T\Mx{B}\otimes\Mx{A}^T\Mx{A})\Tn{G}_{(3)}^T\otimes\Mx{I}_\tau. \\
  \end{aligned}
\end{equation}
Since these blocks can be quite large, our goal-oriented Tucker experiments only explicitly construct and invert the diagonal of these blocks.

%% file: sections/plasma_qois.tex
The calculation of the QoI integrals for the plasma physics data are more complicated than simple sums as in the combustion data, owing to the finite element discretization approach taken in these simulations.  Here the values stored in the tensor do not necessarily correspond to values of the solution variables at points in the mesh, but rather are coefficients for the projection of the solution onto a finite element basis.  The QoI integrals are then given by integrating the discrete solution represented by these basis functions over the computational domain.  This is summarized in \cref{alg:plasma_qoi} for evaluation over a 3-dimensional mesh $\mathcal{N}(\Omega)$ for the finite element approximation to the integral
\begin{equation}
  g(t) = \int_{\Omega} f(\bm{u}(\bm{x},t)) d\bm{x},
\end{equation}
assuming nodal finite element basis functions, which takes as input a data tensor $\Tn{X}$, a set of variable indices $\Vc{v}$ indicating which variables are used in the QoI, a set of time indices $\Vc{t}$ indicating which time steps are used for the QoI, the coordinates $\Vc{x}$, $\Vc{y}$, and $\Vc{z}$ of the nodes in the mesh, a matrix $\Mx{A}\in\Real^{N_{qp}\times N_n}$ that interpolates each variable to the Gauss points in the cell from the nodal values, and a set $\Vc{w}\in\Real^N_{qp}$ quadrature weights.  For trilinear finite element basis functions using first-order Gaussian quadrature $N_n=N_{qp}=8$, $\Vc{w}$ is a vector of all ones, and $\Mx{A} = \Mx{A}_1\otimes\Mx{A}_1\otimes\Mx{A}_1$ where
\begin{equation}
  \Mx{A}_1 =
  \frac{1}{2}
  \begin{bmatrix}
    1+\frac{1}{\sqrt{3}} &  1-\frac{1}{\sqrt{3}} \\
    1-\frac{1}{\sqrt{3}}  & 1+\frac{1}{\sqrt{3}}
  \end{bmatrix}.
\end{equation}
It relies on several mappings/functions that must be initialized when reading the discrete mesh, including \texttt{node\_ind}, which provides a list of indices of mesh nodes for the given element,   \texttt{tensor\_ind}, which provides the $(x,y,z)$ spatial indices of those nodes in the tensor $\X$, and \texttt{det\_jac}, which computes the determinant of the Jacobian matrix of the transformation for a cell determined by its nodal coordinates to the reference cell.  In addition to approximating the above QoI integral, \cref{alg:plasma_qoi} also computes the QoI derivative needed for the Gauss-Newton gradient/Hessian computations described in \cref{sec:go_grad_hess} by approximating
\begin{equation}
  \Mx{Z}(t) = \int_{\Omega} \frac{\partial f}{\partial\bm{u}}(\bm{u}(\bm{x},t)) d\bm{x}.
\end{equation}

\begin{algorithm}
	\caption{Finite-element QoI and derivative evaluation for plasma physics data.}
	\label{alg:plasma_qoi}
	\begin{algorithmic}
		\Require $\Tn{X}$, $\Vc{v}$, $\Vc{t}$, $\Vc{x}$, $\Vc{y}$, $\Vc{z}$, $\Mx{A}$, $\Vc{w}$
		\Ensure $\Vc{g}$, $\Tn{Z}$
		\State $\Vc{g} = 0$
		\State $\Tn{Z} = 0$
		\For{$e\in\mathcal{N}(\Omega)$}
		  \State $\Vc{n} = \mbox{\texttt{node\_ind}}(e)$ \Comment{indices for nodes in mesh cell $e$ (length $N_n$)}
		  \State $\Vc{b} = \mbox{\texttt{det\_jac}}(\Vc{x}(\Vc{n}),\Vc{y}(\Vc{n}),\Vc{z}(\Vc{n}))$ \Comment{determinant of cell transformation Jacobian (length $N_{qp}$)}
		  \State $\Mx{i} = \mbox{\texttt{tensor\_ind}}(\Vc{n})$ \Comment{tensor indices for each node index ($N_n\times 3$)}
		  \For{$k=1,\dots,N_n$} 
		    \State $\Tn{U}(k,:,:) = \Tn{X}(i(k,1),i(k,2),i(k,3),\Vc{v},\Vc{t})$ \Comment{extract values corresponding to given nodes, variables, and time steps}
		  \EndFor
		  \State $\Mx{V}_{(1)} = \Mx{A} \Mx{U}_{(1)}$ \Comment{Interpolate variables to quadrature points}
		  \State $\Mx{F} = f(\Tn{V})$ \Comment{Evaluate QoI integrand at each quadrature point ($N_n\times N_T$)}
		  \State $\Tn{D} = \frac{\partial f}{\partial u}(\Tn{V})$ \Comment{Evaluate derivative of QoI integrand at each quadrature point ($N_{qp}\times N_v \times N_t$)}
		  \State $\Vc{g} \pluseq \sum_{k=1}^8 w(k) b(k) \Vc{f}(k,:)$ \Comment{Cell's contribution to integral}
		   \For{$k=1,\dots,N_n$} 
		     \For{$l=1,\dots,N_{qp}$} 
		       \State $\Tn{Z}(i(k,1),i(k,2),i(k,3),\Vc{v},\Vc{t}) \pluseq w(l) b(l) a(l,k) \Tn{D}(l,:,:)$ \Comment{Cell's contribution to derivative integral}
		     \EndFor
		   \EndFor
		\EndFor
	\end{algorithmic}
\end{algorithm}

%% file: main.bbl
\begin{thebibliography}{10}

\bibitem{genrtr}
{\sc P.-A. Absil, C.~G. Baker, and K.~A. Gallivan}, {\em Trust-region methods
  on {Riemannian} manifolds}, Foundations of Computational Mathematics, 7
  (2007), pp.~303--330, \url{https://doi.org/10.1007/s10208-005-0179-9}.

\bibitem{TTB_CPOPT}
{\sc E.~Acar, D.~M. Dunlavy, and T.~G. Kolda}, {\em A scalable optimization
  approach for fitting canonical tensor decompositions}, Journal of
  Chemometrics, 25 (2011), pp.~67--86, \url{https://doi.org/10.1002/cem.1335}.

\bibitem{acar_cmtf}
{\sc E.~Acar, T.~G. Kolda, and D.~M. Dunlavy}, {\em All-at-once optimization
  for coupled matrix and tensor factorizations}.
\newblock arXiv:1105.3422, 2011, \url{https://arxiv.org/abs/1105.3422}.

\bibitem{AiTuWhKl19}
{\sc M.~Ainsworth, O.~Tugluk, B.~Whitney, and S.~Klasky}, {\em Multilevel
  techniques for compression and reduction of scientific data-quantitative
  control of accuracy in derived quantities}, SIAM Journal on Scientific
  Computing, 41 (2019), pp.~A2146--A2171,
  \url{https://doi.org/10.1137/18M1208885}.

\bibitem{AuBaKo16}
{\sc W.~Austin, G.~Ballard, and T.~G. Kolda}, {\em Parallel tensor compression
  for large-scale scientific data}, in IPDPS'16: Proceedings of the 30th IEEE
  International Parallel and Distributed Processing Symposium, May 2016,
  pp.~912--922, \url{https://doi.org/10.1109/IPDPS.2016.67}.

\bibitem{matlab_ttb_dense}
{\sc B.~W. Bader and T.~G. Kolda}, {\em Algorithm 862: Matlab tensor classes
  for fast algorithm prototyping}, ACM Trans. Math. Softw., 32 (2006),
  p.~635–653, \url{https://doi.org/10.1145/1186785.1186794}.

\bibitem{BaKo07}
{\sc B.~W. Bader and T.~G. Kolda}, {\em Efficient {MATLAB} computations with
  sparse and factored tensors}, SIAM Journal on Scientific Computing, 30
  (2007), pp.~205--231, \url{https://doi.org/10.1137/060676489}.

\bibitem{BHK18}
{\sc G.~Ballard, K.~Hayashi, and R.~Kannan}, {\em Parallel nonnegative {CP}
  decomposition of dense tensors}, in 25th IEEE International Conference on
  High Performance Computing (HiPC), Dec 2018, pp.~22--31,
  \url{https://ieeexplore.ieee.org/document/8638076}.

\bibitem{BaKlKo20}
{\sc G.~Ballard, A.~Klinvex, and T.~G. Kolda}, {\em {TuckerMPI}: A parallel
  {C++/MPI} software package for large-scale data compression via the {Tucker}
  tensor decomposition}, ACM Transactions on Mathematical Software, 46 (2020),
  13 (31~pages), \url{https://doi.org/10.1145/3378445}.

\bibitem{TensorBook}
{\sc G.~Ballard and T.~G. Kolda}, {\em Tensor Decompositions for Data Science},
  Cambridge University Press, 2025,
  \url{https://www.mathsci.ai/post/tensor-textbook/}.

\bibitem{BhagatwalaCL2014}
{\sc A.~Bhagatwala, J.~H. Chen, and T.~Lu}, {\em Direct numerical simulations
  of hcci/saci with ethanol}, Combustion and Flame, 161 (2014), pp.~1826--1841,
  \url{https://doi.org/https://doi.org/10.1016/j.combustflame.2013.12.027},
  \url{https://www.sciencedirect.com/science/article/pii/S0010218014000030}.

\bibitem{biskamp1994magnetic}
{\sc D.~Biskamp}, {\em Magnetic reconnection}, Physics Reports, 237 (1994),
  pp.~179--247.

\bibitem{bittencourt2013fundamentals}
{\sc J.~A. Bittencourt}, {\em Fundamentals of plasma physics}, Springer Science
  \& Business Media, 2013.

\bibitem{BonillaShadidetalVMSMHDMCF2023}
{\sc J.~Bonilla, J.N.Shadid, X.~Tang, M.~Crockatt, P.~Ohm, E.~G. Phillips,
  R.~Pawlowski, and S.~Conde}, {\em On a fully-implicit {VMS}-stabilized {FE}
  formulation for low {M}ach number compressible resistive {MHD} with
  application to {MCF}}, Computer Methods in Applied Mechanics and Engineering,
  417 (2023), p.~116359.

\bibitem{borggaard2016goal}
{\sc J.~Borggaard, Z.~Wang, and L.~Zietsman}, {\em A goal-oriented
  reduced-order modeling approach for nonlinear systems}, Computers \&
  Mathematics with Applications, 71 (2016), pp.~2155--2169.

\bibitem{manopt}
{\sc N.~Boumal, B.~Mishra, P.-A. Absil, and R.~Sepulchre}, {\em {M}anopt, a
  {M}atlab toolbox for optimization on manifolds}, Journal of Machine Learning
  Research, 15 (2014), pp.~1455--1459, \url{https://www.manopt.org}.

\bibitem{bui2007goal}
{\sc T.~Bui-Thanh, K.~Willcox, O.~Ghattas, and B.~van Bloemen~Waanders}, {\em
  Goal-oriented, model-constrained optimization for reduction of large-scale
  systems}, Journal of Computational Physics, 224 (2007), pp.~880--896.

\bibitem{CaCh70}
{\sc J.~D. Carroll and J.~J. Chang}, {\em Analysis of individual differences in
  multidimensional scaling via an {N}-way generalization of ``{Eckart-Young}''
  decomposition}, Psychometrika, 35 (1970), pp.~283--319,
  \url{https://doi.org/10.1007/BF02310791}.

\bibitem{Ca44}
{\sc R.~B. Cattell}, {\em Parallel proportional profiles and other principles
  for determining the choice of factors by rotation}, Psychometrika, 9 (1944),
  pp.~267--283, \url{https://doi.org/10.1007/BF02288739}.

\bibitem{Ca52}
{\sc R.~B. Cattell}, {\em The three basic factor-analytic research designs ---
  their interrelations and derivatives}, Psychological Bulletin, 49 (1952),
  pp.~499--452, \url{https://doi.org/10.1037/h0054245}.

\bibitem{chacon-pop-08-3dmhd}
{\sc L.~Chac\'on}, {\em An optimal, parallel, fully implicit newton-krylov
  solver for three-dimensional visco-resistive magnetohydrodynamics}, Phys.
  Plasmas, 15 (2008), p.~056103.

\bibitem{Chen_S3D}
{\sc J.~H. Chen, A.~Choudhary, B.~de~Supinski, M.~DeVries, E.~R. Hawkes,
  S.~Klasky, W.~K. Liao, K.~L. Ma, J.~Mellor-Crummey, N.~Podhorszki,
  R.~Sankaran, S.~Shende, and C.~S. Yoo}, {\em Terascale direct numerical
  simulations of turbulent combustion using s3d}, Computational Science \&
  Discovery, 2 (2009), p.~015001,
  \url{https://doi.org/10.1088/1749-4699/2/1/015001},
  \url{https://dx.doi.org/10.1088/1749-4699/2/1/015001}.

\bibitem{cheng2016semi}
{\sc L.~Cheng, S.~Mattei, P.~W. Fick, and S.~J. Hulshoff}, {\em A
  semi-continuous formulation for goal-oriented reduced-order models: 1d
  problems}, International Journal for Numerical Methods in Engineering, 105
  (2016), pp.~464--480.

\bibitem{snect}
{\sc D.~Choi and L.~Sael}, {\em {SNeCT}: {Scalable} network constrained
  {Tucker} decomposition for multi-platform data profiling}, IEEE/ACM
  Transactions on Computational Biology and Bioinformatics, 17 (2020),
  pp.~1785--1796, \url{https://doi.org/10.1109/TCBB.2019.2906205}.

\bibitem{cp_nn_constraints}
{\sc J.~Cohen, R.~C. Farias, and P.~Comon}, {\em Fast decomposition of large
  nonnegative tensors}, IEEE Signal Processing Letters, 22 (2015),
  pp.~862--866, \url{https://doi.org/10.1109/LSP.2014.2374838}.

\bibitem{tour_constraints}
{\sc J.~Cohen, K.~Usevich, and P.~Comon}, {\em A tour of constrained tensor
  canonical polyadic decomposition}, May 2016.
\newblock \url{lhttps://hal.archives-ouvertes.fr/hal-01311795/}. Accessed
  05/24/2021.

\bibitem{CoVa25}
{\sc J.~E. Cohen and V.~Leplat}, {\em Efficient algorithms for regularized
  nonnegative scale-invariant low-rank approximation models}, SIAM Journal on
  Mathematics of Data Science, 7 (2025), pp.~468--494,
  \url{https://doi.org/10.1137/24M1657948}.

\bibitem{DeCJV2024}
{\sc S.~De, E.~Corona, P.~Jayakumar, and S.~Veerapaneni}, {\em Tensor-train
  compression of discrete element method simulation data}, Journal of
  Terramechanics, 113-114 (2024), p.~100967,
  \url{https://doi.org/https://doi.org/10.1016/j.jterra.2024.100967},
  \url{https://www.sciencedirect.com/science/article/pii/S0022489824000090}.

\bibitem{DeDeVa00}
{\sc L.~{De Lathauwer}, B.~{De Moor}, and J.~Vandewalle}, {\em A multilinear
  singular value decomposition}, SIAM Journal on Matrix Analysis and
  Applications, 21 (2000), pp.~1253--1278,
  \url{https://doi.org/10.1137/S0895479896305696}.

\bibitem{DeDeVa00a}
{\sc L.~{De Lathauwer}, B.~{De Moor}, and J.~Vandewalle}, {\em On the best
  rank-1 and rank-{$(R_1, R_2, \dots, R_N)$} approximation of higher-order
  tensors}, SIAM Journal on Matrix Analysis and Applications, 21 (2000),
  pp.~1324--1342, \url{https://doi.org/10.1137/S0895479898346995}.

\bibitem{coupled}
{\sc L.~De~Lathauwer and E.~Kofidis}, {\em Coupled matrix-tensor
  factorizations---the case of partially shared factors}, in 2017 51st Asilomar
  Conference on Signals, Systems, and Computers, 2017, pp.~711--715,
  \url{https://doi.org/10.1109/ACSSC.2017.8335436}.

\bibitem{DeHo16}
{\sc H.~De~Sterck and A.~Howse}, {\em Nonlinearly preconditioned optimization
  on grassmann manifolds for computing approximate tucker tensor
  decompositions}, SIAM Journal on Scientific Computing, 38 (2016),
  pp.~A997--A1018, \url{https://doi.org/10.1137/15M1037288}.

\bibitem{ElSa09}
{\sc L.~Eld\'{e}n and B.~Savas}, {\em A {Newton–Grassmann} method for
  computing the best multilinear rank-$(r_1, r_2, r_3)$ approximation of a
  tensor}, SIAM Journal on Matrix Analysis and Applications, 31 (2009),
  pp.~248--271, \url{https://doi.org/10.1137/070688316}.

\bibitem{FanaeetG2016}
{\sc H.~Fanaee-T and J.~Gama}, {\em Tensor-based anomaly detection: An
  interdisciplinary survey}, Knowledge-Based Systems, 98 (2016), pp.~130--147,
  \url{https://doi.org/https://doi.org/10.1016/j.knosys.2016.01.027},
  \url{https://www.sciencedirect.com/science/article/pii/S0950705116000472}.

\bibitem{fang2017efficient}
{\sc F.~Fang, C.~Pain, I.~M. Navon, and D.~Xiao}, {\em An efficient goal-based
  reduced order model approach for targeted adaptive observations},
  International Journal for Numerical Methods in Fluids, 83 (2017),
  pp.~263--275.

\bibitem{farrando2022goal}
{\sc O.~T. Farrando}, {\em A goal-oriented, reduced-order modeling framework of
  wall-bounded shear flows}, PhD thesis, Master’s thesis]. Delft University
  of Technology, 2022.

\bibitem{parafac_constrained}
{\sc G.~Favier and A.~de~Almeida}, {\em Overview of constrained parafac
  models}, EURASIP Journal on Advances in Signal Processing, 2014 (2014),
  p.~142, \url{https://doi.org/10.1186/1687-6180-2014-142}.

\bibitem{GhahremaniB2024}
{\sc B.~Ghahremani and H.~Babaee}, {\em Cross interpolation for solving
  high-dimensional dynamical systems on low-rank tucker and tensor train
  manifolds}, Computer Methods in Applied Mechanics and Engineering, 432
  (2024), p.~117385,
  \url{https://doi.org/https://doi.org/10.1016/j.cma.2024.117385},
  \url{https://www.sciencedirect.com/science/article/pii/S0045782524006406}.

\bibitem{GoedbloedPoedts2004}
{\sc H.~Goedbloed and S.~Poedts}, {\em Principles of Magnetohydrodynamics with
  Applications to Laboratory and Astrophysical Plasmas}, Cambridge Univ. Press,
  2004.

\bibitem{Hackbusch12}
{\sc W.~Hackbusch}, {\em Tensor Spaces and Numerical Tensor Calculus}, vol.~42
  of Springer Series in Computational Mathematics, Springer-Verlag Berlin
  Heidelberg, 2012, \url{https://doi.org/10.1007/978-3-642-28027-6}.

\bibitem{Ha70}
{\sc R.~A. Harshman}, {\em Foundations of the {PARAFAC} procedure: Models and
  conditions for an ``explanatory" multi-modal factor analysis}, UCLA working
  papers in phonetics, 16 (1970), pp.~1--84.
\newblock Available at
  \url{http://www.psychology.uwo.ca/faculty/harshman/wpppfac0.pdf}.

\bibitem{HeChLi23}
{\sc Q.~Heng, E.~C. Chi, and Y.~Liu}, {\em Robust low-rank tensor decomposition
  with the {L2} criterion}, Technometrics, 65 (2023), pp.~537--552,
  \url{https://doi.org/10.1080/00401706.2023.2200541}.

\bibitem{Hi27a}
{\sc F.~L. Hitchcock}, {\em The expression of a tensor or a polyadic as a sum
  of products}, Journal of Mathematics and Physics, 6 (1927), pp.~164--189,
  \url{https://doi.org/10.1002/sapm192761164}.

\bibitem{Hi27}
{\sc F.~L. Hitchcock}, {\em Multiple invariants and generalized rank of a p-way
  matrix or tensor}, Journal of Mathematics and Physics, 7 (1927), pp.~39--79,
  \url{https://doi.org/10.1002/sapm19287139}.

\bibitem{HoKoDu2020}
{\sc D.~Hong, T.~G. Kolda, and J.~A. Duersch}, {\em Generalized canonical
  polyadic tensor decomposition}, SIAM Review, 62 (2020), pp.~133--163,
  \url{https://doi.org/10.1137/18M1203626}.

\bibitem{HuSiLi}
{\sc K.~Huang, N.~D. Sidiropoulos, and A.~P. Liavas}, {\em A flexible and
  efficient algorithmic framework for constrained matrix and tensor
  factorization}, IEEE Transactions on Signal Processing, 64 (2016),
  pp.~5052--5065, \url{https://doi.org/10.1109/TSP.2016.2576427}.

\bibitem{Ion_Phd_thesis_2024}
{\sc I.~G. Ion}, {\em Low-rank tensor decompositions for surrogate modeling in
  forward and inverse problems}, PhD thesis, Technische Universit{\"a}t
  Darmstadt, Darmstadt, March 2024,
  \url{https://doi.org/https://doi.org/10.26083/tuprints-00026678},
  \url{http://tuprints.ulb.tu-darmstadt.de/26678/}.

\bibitem{KizhakkinanDLVRR2023}
{\sc U.~Kizhakkinan, P.~L.~T. Duong, R.~Laskowski, G.~Vastola, D.~W. Rosen, and
  N.~Raghavan}, {\em Development of a surrogate model for high-fidelity laser
  powder-bed fusion using tensor train and gaussian process regression},
  Journal of Intelligent Manufacturing, 34 (2023), pp.~369--385,
  \url{https://doi.org/https://doi.org/10.1007/s10845-022-02038-4}.

\bibitem{KoBa09}
{\sc T.~G. Kolda and B.~W. Bader}, {\em Tensor decompositions and
  applications}, SIAM Review, 51 (2009), pp.~455--500,
  \url{https://doi.org/10.1137/07070111X}.

\bibitem{Kolla2020}
{\sc H.~Kolla, K.~Aditya, and J.~H. Chen}, {\em Higher order tensors for dns
  data analysis and compression}, in Data Analysis for Direct Numerical
  Simulations of Turbulent Combustion: From Equation-Based Analysis to Machine
  Learning, H.~Pitsch and A.~Attili, eds., Springer International Publishing,
  Cham, 2020, pp.~109--134, \url{https://doi.org/10.1007/978-3-030-44718-2_6},
  \url{https://doi.org/10.1007/978-3-030-44718-2_6}.

\bibitem{LeGoCh22}
{\sc J.~Lee, Q.~Gong, J.~Choi, T.~Banerjee, S.~Klasky, S.~Ranka, and
  A.~Rangarajan}, {\em Error-bounded learned scientific data compression with
  preservation of derived quantities}, Applied Sciences, 12 (2022),
  \url{https://doi.org/10.3390/app12136718}.

\bibitem{LiYuan2023}
{\sc X.~Li and Y.~Yuan}, {\em A tensor-based hyperspectral anomaly detection
  method under prior physical constraints}, IEEE Transactions on Geoscience and
  Remote Sensing, 61 (2023), pp.~1--12,
  \url{https://doi.org/10.1109/TGRS.2023.3324147}.

\bibitem{Li_etal_PINO_2024}
{\sc Z.~Li, H.~Zheng, N.~Kovachki, D.~Jin, H.~Chen, B.~Liu, K.~Azizzadenesheli,
  and A.~Anandkumar}, {\em Physics-informed neural operator for learning
  partial differential equations}, ACM / IMS J. Data Sci., 1 (2024),
  \url{https://doi.org/10.1145/3648506}, \url{https://doi.org/10.1145/3648506}.

\bibitem{NocedalWright}
{\sc J.~Nocedal and S.~J. Wright}, {\em Numerical Optimization}, Springer
  Series in Operations Research, Springer-Verlag, New York, 2006.

\bibitem{Pa97}
{\sc P.~Paatero}, {\em A weighted non-negative least squares algorithm for
  three-way ``{PARAFAC}'' factor analysis}, Chemometrics and Intelligent
  Laboratory Systems, 38 (1997), pp.~223--242,
  \url{https://doi.org/10.1016/S0169-7439(97)00031-2}.

\bibitem{PhTiCi13a}
{\sc A.-H. Phan, P.~Tichavsk\'{y}, and A.~Cichocki}, {\em Low complexity damped
  {Gauss}--{Newton} algorithms for {CANDECOMP/PARAFAC}}, SIAM Journal on Matrix
  Analysis and Applications, 34 (2013), pp.~126--147,
  \url{https://doi.org/10.1137/100808034}.

\bibitem{PhTiCi2011}
{\sc A.~H. Phan, P.~Tichavský, and A.~Cichocki}, {\em Damped {Gauss-Newton}
  algorithm for nonnegative {Tucker} decomposition}, in 2011 IEEE Statistical
  Signal Processing Workshop (SSP), 2011, pp.~665--668,
  \url{https://doi.org/10.1109/SSP.2011.5967789}.

\bibitem{cp_rank_constraints}
{\sc A.-H. Phan, P.~Tichavský, K.~Sobolev, K.~Sozykin, D.~Ermilov, and
  A.~Cichocki}, {\em Canonical polyadic tensor decomposition with low-rank
  factor matrices}, in ICASSP 2021 - 2021 IEEE International Conference on
  Acoustics, Speech and Signal Processing (ICASSP), 2021, pp.~4690--4694,
  \url{https://doi.org/10.1109/ICASSP39728.2021.9414606}.

\bibitem{QiLu17}
{\sc L.~Qi and Z.~Luo}, {\em Tensor Analysis: Spectral Theory and Special
  Tensors}, Society for Industrial and Applied Mathematics, Philadelphia, PA,
  2017, \url{https://doi.org/10.1137/1.9781611974751}.

\bibitem{RaissiPK_PINN_2019}
{\sc M.~Raissi, P.~Perdikaris, and G.~Karniadakis}, {\em Physics-informed
  neural networks: A deep learning framework for solving forward and inverse
  problems involving nonlinear partial differential equations}, Journal of
  Computational Physics, 378 (2019), pp.~686--707,
  \url{https://doi.org/https://doi.org/10.1016/j.jcp.2018.10.045},
  \url{https://www.sciencedirect.com/science/article/pii/S0021999118307125}.

\bibitem{RaBa20}
{\sc T.~M. Ranadive and M.~M. Baskaran}, {\em Large-scale sparse tensor
  decomposition using a damped {Gauss-Newton} method}, in 2020 IEEE High
  Performance Extreme Computing Conference (HPEC), 2020, pp.~1--8,
  \url{https://doi.org/10.1109/HPEC43674.2020.9286202}.

\bibitem{RoScCaAd22}
{\sc M.~Roald, C.~Schenker, V.~D. Calhoun, T.~Adali, R.~Bro, J.~E. Cohen, and
  E.~Acar}, {\em An ao-admm approach to constraining parafac2 on all modes},
  SIAM Journal on Mathematics of Data Science, 4 (2022), pp.~1191--1222,
  \url{https://doi.org/10.1137/21M1450033}.

\bibitem{ShadidetalResistiveMHD2016}
{\sc J.~N. Shadid, R.~P. Pawlowski, E.~C. Cyr, R.~S. Tuminaro, L.~Chacon, and
  P.~D. Weber}, {\em {Scalable Implicit Incompressible Resistive {MHD} with
  Stabilized FE and Fully-coupled {N}ewton-{K}rylov-{AMG}}}, Comput. Methods
  Appl. Mech. Engrg., 304 (2016), pp.~1--25.

\bibitem{singh2020comparison}
{\sc N.~Singh, L.~Ma, H.~Yang, and E.~Solomonik}, {\em Comparison of accuracy
  and scalability of {Gauss-Newton} and alternating least squares for {CP}
  decomposition}, 2020, \url{https://arxiv.org/abs/1910.12331}.

\bibitem{Tu66}
{\sc L.~R. Tucker}, {\em Some mathematical notes on three-mode factor
  analysis}, Psychometrika, 31 (1966), pp.~279--311,
  \url{https://doi.org/10.1007/BF02289464}.

\bibitem{VaVeLa21}
{\sc M.~Vandecappelle, N.~Vervliet, and L.~D. Lathauwer}, {\em Inexact
  generalized {Gauss–Newton} for scaling the canonical polyadic decomposition
  with non-least-squares cost functions}, IEEE Journal of Selected Topics in
  Signal Processing, 15 (2021), pp.~491--505,
  \url{https://doi.org/10.1109/JSTSP.2020.3045911}.

\bibitem{VaVaMe12}
{\sc N.~Vannieuwenhoven, R.~Vandebril, and K.~Meerbergen}, {\em A new
  truncation strategy for the higher-order singular value decomposition}, SIAM
  Journal on Scientific Computing, 34 (2012), pp.~A1027--A1052,
  \url{https://doi.org/10.1137/110836067}.

\bibitem{ZhByLuNo1997}
{\sc C.~Zhu, R.~H. Byrd, P.~Lu, and J.~Nocedal}, {\em Algorithm 778:
  {L-BFGS-B}: {Fortran} subroutines for large-scale bound-constrained
  optimization}, ACM Trans. Math. Softw., 23 (1997), p.~550–560,
  \url{https://doi.org/10.1145/279232.279236}.

\end{thebibliography}
